\documentclass[a4paper,11pt]{article}

\usepackage{amssymb}
\usepackage{amsmath}
\usepackage{amsfonts}
\usepackage{stmaryrd}
\usepackage{theorem}
\usepackage[dvips]{graphicx}
\usepackage[all]{xy}
\usepackage{pstricks}
\usepackage{fancybox}
\usepackage{tabls}
\usepackage{array}
\usepackage{boxedminipage}

	\usepackage{fancyhdr}

\newtheorem{Thm}{Theorem}[section]
{\theoremstyle{break} }
\newtheorem{Prop}[Thm]{Proposition}
{\theoremstyle{break} }
\newtheorem{Lemma}[Thm]{Lemma}
{\theoremstyle{break} }
\newtheorem{Cor}[Thm]{Corollary}
{\theoremstyle{break} }
\newtheorem{Def}[Thm]{Definition}
{\theoremstyle{break} }

{\theoremstyle{break} }
{\theorembodyfont{\rm} \newtheorem{Rem}[Thm]{Remark}}
{\theoremstyle{break} {\theorembodyfont{\rm} \newtheorem{Remb}[Thm]{Remark}}}
{\theorembodyfont{\rm} \newtheorem{Ex}[Thm]{Example}} 
{\theoremstyle{break} {\theorembodyfont{\rm} }}
{\theorembodyfont{\rm} } 
{\theoremstyle{break} {\theorembodyfont{\rm} }}

\newenvironment{Proof}[1]{\textsc{Proof} \textsc{#1}\quad}{ \hspace*{\fill} $\Box$ \\}

\newcommand{\R}{\mathbb{R}} \newcommand{\C}{\mathbb{C}} 
 \newcommand{\N}{\mathbb{N}} \newcommand{\Z}{\mathbb{Z}}
  \newcommand{\T}{\mathcal{T}}
 \newcommand{\E}{\mathcal{E}} \newcommand{\J}{\mathcal{J}}
\newcommand{\cR}{\mathcal{R}}  
 
\newcommand{\map}{\longrightarrow} 
 \newcommand{\smap}{\twoheadrightarrow}

\newcommand{\Sym}{\mathrm{Sym}}

  \newcommand{\LC}{\nabla}

\newcommand{\qq}[1]{\underline{\smash{#1}}}

\newcommand{\PR}{\mathcal{P}}

\newcommand{\M}{\mathcal{M}}
\newcommand{\del}{\partial}

\newcommand{\supp}{\text{supp}}

\newcommand{\tsum}{\textstyle\sum}
\newcommand{\F}{\mathcal{F}}

\newcommand{\Diff}{\mathrm{Diff}}
\newcommand{\Der}{\mathrm{Der}}

\newcommand{\sign}{\mathrm{sign}} 
\newcommand{\tL}{\textstyle\bigwedge}

\newcommand{\Hom}{\mathrm{Hom}} \newcommand{\Set}{\mathrm{Set}}\newcommand{\End}{\mathrm{End}}

\newcommand{\Top}{\mathrm{Top}}
\newcommand{\TGr}{\mathrm{Top}^{\mathrm{Gr}}}
\newcommand{\Man}{\mathrm{Man}}
\newcommand{\MGr}{\mathrm{Man}^{\mathrm{Gr}}}
\newcommand{\VBun}{\mathrm{VBun}}
\newcommand{\Gr}{\mathrm{Gr}}
\newcommand{\BKL}{\mathrm{BKL}}
\newcommand{\0}{{\overline{0}}} 
\newcommand{\pr}{\mathrm{pr}}

\newcommand{\odd}{\mathrm{odd}} \newcommand{\ev}{\mathrm{ev}}

\newcommand{\SMan}{\mathrm{SMan}}
\newcommand{\Jet}{\mathrm{Jet}}\newcommand{\jet}{\mathrm{jet}}

\newcommand{\sss}{\scriptscriptstyle}
\newcommand{\nil}{\mspace{-4mu}\mathrm{nil}}
\newcommand{\SC}{SC}
\newcommand{\op}{\mathrm{op}}
\newcommand{\symb}{\mathrm{sb}^{\! \scriptscriptstyle \nabla}}
\newcommand{\Id}{\mathrm{Id}}

\renewcommand{\O}{\mathcal{O}}
\renewcommand{\emptyset}{\varnothing}
\renewcommand{\Vec}{\mathrm{Vec}}
\renewcommand{\1}{{\overline{1}}} \renewcommand{\2}{{\overline{2}}}

\newcounter{Z0}\newcounter{Z1}
\newcounter{Z2}\newcounter{Z3}
\newcounter{Z4}\newcounter{Z5}
\newcounter{Z6}\newcounter{Z7}
\newcounter{Z8}\newcounter{Z9}
\numberwithin{equation}{section}

\title{A Supermanifold structure on Spaces of Morphisms between Supermanifolds}
\author{Florian Hanisch}

\author{Florian Hanisch \vspace{4mm}\\
{\small Institut f\"ur Mathematik} \\
{\small Universit\"at Potsdam, Am Neuen Palais 10, 14469 Potsdam, Germany.}\vspace{2mm}\\
{\footnotesize fhanisch@math.uni-potsdam.de~}
}

\begin{document}

\maketitle

\begin{abstract}
The aim of this work is the construction of a ``supermanifold of morphisms $X \rightarrow Y$'', given two finite-dimensional supermanifolds $X$ and $Y$. More precisely, we will define an object $\qq{\SC}^\infty(X,Y)$ in the category of supermanifolds proposed by Molotkov and Sachse. Initially, it is given by the set-valued functor characterised by the adjunction formula $\Hom(P \times X,Y) \cong \Hom(P,\qq{SC}^\infty(X,Y))$ where $P$ ranges over all superpoints. We determine the structure of this functor in purely geometric terms: We show that it takes values in the set of certain differential operators and establish a bijective correspondence to the set of sections in certain vector bundles associated to $X$ and $Y$. Equipping these spaces of sections with infinite-dimensional manifold structures using the convenient setting by Kriegl and Michor, we obtain at a supersmooth structure on $\qq{\SC}^\infty(X,Y)$, i.e. a supermanifold of all morphisms $X \rightarrow Y$.      

\vspace{2cm}
\textit{Key words:} supergeometry, infinite-dimensional geometry, mapping spaces, fermionic geometry 
\end{abstract}

\pagestyle{fancy}
\pagenumbering{arabic}

\newpage

\section{Introduction} \label{ChapIntro}

Superanalysis and supermanifolds were introduced in theoretical physics in the 1960s. The need for spaces which incorporate both anticommuting (or fermionic) and commuting (or bosonic) degrees of freedom on the same footing first arose in quantum theory. Especially in the path integral approach, these concepts provided suitable classical configuration spaces for setting up that theory. Later on, supermanifolds also served as ``superspacetimes'' providing a geometric framework to write down certain field theories in an intrinsically supersymmetric way.\\
The literature from mathematics and mathematical physics contains several (partially equivalent) definitions for supermanifolds. Taking into account the need to include anticommuting (and hence nilpotent) elements in their rings of functions, a natural approach uses the language of ringed spaces. It is mainly due to Berezin (\cite{Ber}), Leites (\cite{Lei}) and Kostant (\cite{Kos}). A different approach, due to Rogers (see \cite{Rog1} for an overview), DeWitt (\cite{DeWitt}) and others, replaces $\R^n$ by modules of the form $(\Lambda_\0)^p \times (\Lambda_\1)^q$ for some Grassmann algebra $\Lambda = \Lambda^{\mathrm{ev}} \oplus \Lambda^{\mathrm{odd}}$ and provides an appropriate scheme for glueing them to obtain $p|q$-dimensional supermanifolds. A comparison and overview can e.g. be found in \cite{BBHR}.\\
All these concepts describe \emph{finite-dimensional} supermanifolds. However, the spaces of greatest physical interest are infinite-dimensional: Configuration spaces of bosonic and fermionic fields in quantum theory or spaces of superfields and their component fields in applications to supersymmetric theories. Again, there are different approaches to define infinite-dimensional supermanifolds. Batchelor (\cite{Bat2}) and Jaskolski (\cite{Jas}) used a coalgebraic approach extending ideas from Kostant's work. Schmitt extended the ringed space approach by a generalisation of Kostant's finitely supported distributions (\cite{Schmitt}, \cite{Schmitt23}, vol. I, p.25 ff and section 3.1, 3.4) to include infinite-dimensional analytic supermanifolds. To the authors knowledge, all these proposals eventually do not provide a construction of an infinite-dimensional ``supermanifold of all morphisms $X \rightarrow Y$'', where $X$ and $Y$ are arbitrary, finite-dimensional supermanifolds. A different approach has been suggested by Molotkov (\cite{Mol}) which was later extended by Alldridge and Laubinger (\cite{AL}) to a larger class of base fields. Our construction of mapping spaces is based on these ideas. Finally, Alldridge (\cite{All}) more recently defined a category of supermanifolds based on Douady's notion of ``functored spaces'' and proved (under certain assumptions) the existence of inner Hom objects. The comparison of these results and approaches with our construction is beyond the scope of this work but it would be interesting to address this task in the future.\\   

The aim of the present work is the construction of a ``supermanifold of all morphisms $X \rightarrow Y$'', denoted by $\qq{SC}^\infty(X,Y)$, in the framework by Molotkov. The latter was worked out in great detail by Sachse (\cite{Sa1}, \cite{Sa2}) and has been already used in \cite{SaWo} to construct the infinite-dimensional supermanifold of diffeomorphisms of a compact supergroup. The basic idea (which will be explained in chapter \ref{ChapSuper} in more detail) is as follows: Supermanifolds are defined to be certain functors from a sufficiently large category of ``test objects'' (here: all supermanifolds of dimension $0|n$) into a suitable category of (possibly infinite-dimensional) manifolds. Molotkov's work defines the notion of Banach-supermanifolds which clearly needs to be enlarged to discuss spaces of smooth mappings. We will choose the category of manifolds defined by Kriegl and Michor (\cite{KM}), which provides manifold structures for spaces of smooth mappings $C^\infty(M,N)$, where $M,N$ are finite-dimensional smooth manifolds. In particular, $M$ need not be compact which is crucial for physical applications: Spacetime should be modelled by globally hyperbolic Lorentzian manifolds (see \cite{BGP}, \cite{QFTPotsdam}) and these are never compact.\\
The rough idea of our construction is as follows: The functorial framework in fact suggests a candidate for the space of morphisms: The functor $\qq{SC}^\infty(X,Y)$ defined in \ref{DefSCInfty} by an adjunction formula taking values in the category of sets. These sets can be described in terms of differential operators and we arrive at a geometric identification of the values of $\qq{SC}^\infty(X,Y)$ with sections in certain bundles (Theorem \ref{ThmPointsSCInfty}, \ref{CorPointsSCInfty}). Under this identification, there is an obvious manifold structure on the values of $\qq{SC}^\infty(X,Y)$ (see Theorem \ref{ThmBundleStruct}) and an appropriate supersmooth atlas can finally be constructed (Theorem \ref{ThmSMfSCInfty}).\\
The motivation to study these particular spaces is twofold: They provide a domain of definition for supersymmetric action functionals (see e.g. \cite{AMS2}) and can hopefully be a starting point to obtain a precise geometric description of the corresponding spaces of critical points.  Secondly, these spaces (or subspaces thereof) can possibly serve as configuration spaces for fermionic fields in (not necessarily supersymmetric) classical field theories, similar to the bosonic case in \cite{BFLR}.  

This work is organized as follows: In the the second chapter, we provide some prerequisites which are mostly known but sometimes difficult to find: First, we give a brief (and very incomplete) overview of the convenient calculus on infinite dimensional vector spaces and its use to equip $C^\infty(M,N)$ with the structure of a manifold. Secondly, we discuss the Riemannian geometry of the total space of a vector bundle, its local product structure and the theory of its jet bundles. Here, we make use of a connection to identify jets with symmetric algebras and derive some formulas. The third chapter contains an introduction to the ringed space picture of supergeometry and the functorial Molotkov-Sachse approach mentioned above. At the end, we slightly generalize this approach to locally affine model spaces in order to include those supermanifolds constructed in the following chapters. Chapter 4 contains the structural results describing $\qq{SC}^\infty(X,Y)$ as a set-valued functor; these are summarized in the fundamental Theorems \ref{ThmStructComp}, \ref{ThmStructComp2} and \ref{ThmPointsSCInfty}. More precisely, we obtain a characterization of this functor in terms of the tangent bundles associated to $X$ and $Y$; a similar result has been obtained independently in \cite{BK} for a special case. We will prove these results using an algebraic formulation of the theory of differential operators on supermanifolds which will be briefly reviewed at the begin of the chapter. It should be mentioned that a Chapter 5 is a brief detour back to finite-dimensional supermanifolds. Based on the results from chapter 4, we arrive at an explicit geometric description of the functor of points associated to finite-dimensional supermanifolds. Related results have been obtained in a different fashion in \cite{BCF1}, \cite{Sa2}, \cite{AL}. In particular, we recover the basic concepts familiar from the Rogers-DeWitt approach in a geometric language. Motivated by these finite-dimensional results, we eventually construct a supersmooth structure on $\qq{SC}^\infty(X,Y)$ in chapter 6; the result being given in the main Theorem \ref{ThmSMfSCInfty}. To this end, we first construct an infinite-dimensional bundle structure on $C^\infty(M,S)$ for a vector bundle $S \rightarrow N$.  The sets $\qq{SC}^\infty(X,Y)(\Lambda_n)$ are eventually identified with certain submanifolds of such bundles.\\

{\bf Acknowledgements:} Special thanks to Gregor Weingart for many helpful discussions concerning jet formalism, to Peter Michor for providing the bundle construction used in Remark \ref{RemMichor} on MathOverFlow and to the anonymous author ``TaQ" on MathOverFlow, who indicated the argument in Remark \ref{RemCInftyTop}, iii). I am also grateful to Tillmann Wurzbacher and Christoph Wockel for interesting discussions on infinite-dimensional- and supergeometry and the geometry group in Potsdam (in particular C.~B\"{a}r, C.~Becker and M.~Ludewig) for comments and suggestions. Finally, thanks also go to SFB 647 funded by Deutsche Forschungsgemeinschaft and Max Planck Institut f\"{u}r Gravitationsphysik for financial support.\\

{\bf Notations and Conventions}: 
If $\mathcal{C}, \mathcal{D}$ are two categories, $\mathcal{D}^\mathcal{C}$ denotes the category whose objects are functors $\mathcal{C} \rightarrow \mathcal{D}$ and whose morphisms are natural transformations among them. Moreover, $\mathcal{C}^{opp}$ denotes the opposite category, i.e. $\mathcal{C}$ with arrows reversed. We will use the following categories:
\begin{itemize}
 \item $\mathrm{SALg}$ the category of (real) superalgebras. Morphisms are assumed to preserve the $\Z_2$ parity. $\Gr \subset \mathrm{SAlg}$ the full subcategory of real, finite-dimensional Grassmann 
       algebras. Generic objects of $\Gr$ will be denoted by $\Lambda, \Lambda'$ and morphisms by $\rho \in \Hom(\Lambda,\Lambda')$. Starting from chapter \ref{ChapAlg}, we will restrict to a skeleton of $\Gr$ and simply work with the algebras $\Lambda_n := \tL \R^n$.
 \item $\Set$ the category of sets, $\Top$ the category of topological spaces, $\mathrm{cVec}$ the category of convenient vector spaces (cf. chapter \ref{ChapPre}) 
 \item $\Man$ the category of manifolds according to \cite{KM}; $\Man_F$ and $\mathrm{fMan}$ its subcategories of Fr\'{e}chet- and finite-dimensional manifolds, respectively. Morphisms are always given by the 
       corresponding smooth maps.  
 \item $\BKL$ the category of supermanifolds as defined by Berezin-Kostant-Leites, $\SMan$ the category of supermanifolds in the Molotkov-
       Sachse-framework build upon $\Man$ (cf. chapter \ref{ChapSuper} for details).
\end{itemize}
When working with $\Z_2$-graded objects, we will use $\0$ and $\1$ to denote this grading, e.g. $M = M_\0 \oplus M_1$ for a supermodule $M$. We will furthermore decompose real Grassmann algebras $\Lambda$ as $\Lambda = \R \oplus \Lambda^{\nil} = \R \oplus \Lambda^{\ev \geq 2} \oplus \Lambda^{\odd}$. Here, $\Lambda^{\nil}$ is the space of all nilpotent elements, $\Lambda^{\odd}$ the space of all elements of odd degree and $\Lambda^{\ev \geq 2} = \Lambda^{\ev} \cap \Lambda^{\nil}$ the space of all elements of even degree greater or equal than 2. If $V$ is some vector space, we will denote the resulting  decomposition of $\mu \in V \otimes \Lambda$ by $\mu = \tilde{\mu} + \mu^{\nil} = \tilde{\mu} + \mu_\2 + \mu_\1 = \mu_\0 + \mu_\1$, where $\mu_\2$ is the component in $\Lambda^{\ev \geq 2}$.\\

A multi index $I= (i_1,\ldots,i_k)$ will denote an element of $\N_0^k$ for some $k \in \N$. As usual, we set $|I| := i_1 + \cdots + i_k$ and $I! = i_1! \cdots i_k!$. Let $\N_0^{k,l} := \{I \in \N_0^k \mid |I|\leq l \}$ denote the set of multi indices $I$ with length at most $l$. An analogous convention is used for indices $J \in \Z_2^q$. \\


\section{Geometric Prerequisites} \label{ChapPre}

We start with a brief sketch of convenient calculus and its applications to construct infinite-dimensional manifolds. We will follow the approach explained in great detail in \cite{KM}.\\

Let $V,W$ be locally convex vector spaces. The fundamental concept for the discussion of smooth maps $V \rightarrow W$ is that of smooth curves in $V,W$. Defining the derivative by $c'(t) := \lim_{s\rightarrow 0} (c(t+s)-c(s))/s$, a curve $c : \R \rightarrow V$ is called smooth if all its iterated derivatives exist. $V$ can be equipped with the $c^\infty$-topology (\cite{KM} 2.12), which is defined to be the final topology w.r.t. $C^\infty(\R,V)$, i.e. the finest topology on $V$ s.t. all smooth curves are continuous. It is finer than the original locally convex topology and the vector space $V$, equipped with the $c^\infty$-topology, will be denoted by $c^\infty V$. As usual, one needs some sort of completeness assumption on $V$. Following \cite{KM}, we will be mainly interested in so called \emph{convenient vector spaces}, which may be characterized by the property that smoothness of curves can be tested by composing it with continuous linear functionals (see \cite{KM} 2.14, this theorem also provides many other equivalent characterizations). Smoothness of maps $V \rightarrow W$ is now defined on arbitrary $c^\infty$-open subsets of $V$ by testing with smooth curves:

\begin{Def}[\cite{KM}, 3.11] \label{DefSmoothMap}
A map $f : U \subset V \rightarrow W$, defined on the $c^\infty$-open set $U$, is called smooth, if $f_\ast$ maps $C^\infty(\R,U)$ to $C^\infty(\R,W)$. 
\end{Def}

Let $\mathrm{cVec}$ denote the category of convenient vector spaces with morphisms given by smooth linear maps. If $U$ is equipped with the $c^\infty$-topology, it now follows immediately from the properties of the final topology that smooth maps $U \rightarrow W$ are continuous. Similarly, a diffeomorphism $f : U \subset V \rightarrow U' \subset W$ is a homeomorphism, if $U$ and $U'$ are equipped with the respective $c^\infty$-topology.\\

In the sequel, we will mainly be interested in the following examples, which provide the model spaces for manifolds of mappings: Let $\pi: S \rightarrow M$ be a vector bundle of finite rank over the finite-dimensional manifold $M$ (which can be compact or not). Consider the spaces of sections
\begin{align} \label{EqModell}
\Gamma(S) &:= \{\sigma \in C^\infty(M,S) \mid \pi \sigma = \Id_M\} \\
\Gamma_c(S) &:= \{\sigma \in C^\infty(M,S) \mid \pi \sigma = \Id_M, \ \mathrm{supp}(\sigma) \text{ is compact }\},   \notag
\end{align}
equipped with the usual locally convex, complete, nuclear topology. They turn out to be convenient vector spaces; several equivalent characterizations of   
their convenient structure can be found in chapter 30 of \cite{KM}, in particular 30.1, 30.3 and 30.4. In view of Definition \ref{DefSmoothMap}, we will not discuss these details and simply quote the relevant characterization of smooth of curves:

\begin{Lemma}[\cite{KM} 30.8, 30.9] \label{LemSmoothCurvVS}
A curve $c : \R \rightarrow \Gamma(S)$ is smooth iff the induced map \linebreak $c^\wedge : \R \times M \rightarrow S, c^\wedge(t,x) := c(t)(x)$ is smooth in the ordinary sense. A curve $c : \R \rightarrow \Gamma_c(S)$ is smooth, iff $c^\wedge$ is smooth and in addition, for each interval $[a,b] \subset \R$, there exists $K \subset M$ compact s.t. $c^\wedge(\cdot,x)$ is constant on $[a,b]$ for every $x \in M\setminus K$.
\end{Lemma}

Before discussing manifolds, we list a few properties of $c^\infty$-topologies, especially with regard to the space $\Gamma_c(S)$: 

\begin{Remb} \label{RemCInftyTop}
\begin{enumerate}
 \item The $c^\infty$-topology on $V$ is finer than the original locally convex one, it agrees with it under certain circumstances (e.g. if $V$ is metrizable,  see \cite{KM} 4.11 for more details)
 \item $c^\infty V$ is not always a topological vector space. In particular, this is not the case for $\Gamma_c(S)$ unless $M$ is compact since addition fails to be continuous (see \cite{KM} 4.26, Remark).   
 \item $c^\infty(V \times W) \rightarrow c^\infty(V) \times c^\infty(W)$ is continuous bijection and in case $W \cong \R^n$, it is a homeomorphism (\cite{KM} 4.16). In general, the topology of $c^\infty(V 
       \times W)$ may be strictly finer than the product topology on $c^\infty(V) \times c^\infty(W)$. Again, this happens for $V = W = \Gamma_c(S)$ unless $M$ is compact. The following argument was hinted to me by TaQ on MathOverFlow: Addition, considered as a map $c^\infty(\Gamma_c(S) \times \Gamma_c(S)) \rightarrow c^\infty\Gamma_c(S)$ is continuous (even smooth). If $c^\infty\Gamma_c(S) \times c^\infty\Gamma_c(S) \rightarrow c^\infty(\Gamma_c(S) \times \Gamma_c(S))$ was continuous, so would be addition as a map $c^\infty\Gamma_c(S) \times c^\infty\Gamma_c(S) \rightarrow c^\infty\Gamma_c(S)$. This contradicts (b).
\end{enumerate}
\end{Remb}

Manifolds are defined in the usual way, i.e. by considering charts $u_i : U_i \subset M \rightarrow u_i(U_i) \subset V_i$ from a set $M$ on a $c^\infty$-open subset $u_i(U_i)$ of a convenient vector space $V_i$. Atlases are obtained by requiring chart changes to be smooth in the sense of Definition \ref{DefSmoothMap} and it can be shown that maps between manifolds are smooth iff they map smooth curves to smooth curves. A detailed discussion can be found in \cite{KM} (27.1 to 27.4). We only note that the topology on the set $M$ is chosen to be the one induced by the $c^\infty$-topologies of the charts $V$, i.e. $U \subset M$ is open iff $u_i(U \cap U_i) \subset V_i$ is open for all $i$, which is required to be smoothly Hausdorff (\cite{KM}, 27.4). A slightly different definition is used in \cite{Michoralt}, Def. 9.1. We will not discuss tangent structures but refer to chapter 28 of \cite{KM}.\\

A very simple example of a manifold, which will be used subsequently, is given by a locally affine space, taken from \cite{Michoralt}, 12.1 and 4.10. 

\begin{Ex} \label{ExLocAffine}
Assume that $M$ is not compact. Consider the subspace $\Gamma_c(S) \subset \Gamma(S)$ and the algebraic quotient $\Gamma(S) / \Gamma_c(S)$. Then clearly, two sections $s_1,s_2 \in \Gamma(S)$ are in the same  equivalence class of the quotient ($s_1 \sim s_2$) iff $s_1 - s_2$ has compact support. We may now put a new, finer topology on $\Gamma(S)$ (called $\mathrm{F}\mathcal{D}$-topology in \cite{Michoralt}) by putting the topology induced from $\Gamma_c(S)$ on each equivalence class and taking their disjoint union. This is clearly no longer a vector space topology since for $s \in \Gamma(S)\setminus \Gamma_c(S)$, $t\cdot s \notin \Gamma_c(S)$ for $t \neq 0$ and hence, $t\cdot s$ does not converge to zero for $t\searrow 0$.\\
However, $\Gamma(S)$ is a manifold modelled on $\Gamma_c(S)$, called a \emph{locally affine space}. In fact, a chart near $s_0 \in \Gamma(S)$ is given by $u_{s_0} : s_0 + \Gamma_c(S) \rightarrow \Gamma_c(S), s \mapsto s-s_0$. $\Gamma(S)$ is simply the disconnected union of components diffeomorphic to $\Gamma_c(S)$. Their union does not form a topological vector space even though it is a vector space. As also remarked in \cite{Michoralt}, 12.1, the tangent bundle of the locally affine space is simply $T\Gamma(S) = \Gamma(S) \times \Gamma_c(S)$ and differentiation works as if $\Gamma(S)$ was a topological vector space.\\
Interestingly, this type of structure has also have been considered recently in the physics literature in the context of spaces of field configurations in classical field theory, see \cite{BFLR} section 2.2.
\end{Ex}

We now briefly recall the construction of the infinite-dimensional structure on $C^\infty(M,N)$, described in chapter 42 of \cite{KM}, for arbitrary smooth, finite-dimensional manifolds $M,N$. The model space near $f \in C^\infty(M,N)$ is given by the convenient vector space $\Gamma_c(f^\ast TN)$. It induces a very fine topology on $C^\infty(M,N)$; maps being in the same connected component necessarily coincide outside of a compact set $K \subset M$. However, in contrast to many other topologies, this choice makes $C^\infty(M,N)$ locally contractible and hence allows for the construction of charts, even if $M$ is not compact. A detailed discussion of this issue can be found in \cite{KM} (introduction to chapter IX), \cite{Michoralt} (4.9 to 4.11) or \cite{Michor1} (section 2 and Rem. 3.5). To define charts for $C^\infty(M,N)$, we choose a Riemannian metric on $N$ and open sets $V \subset TN$, $U \subset N \times N$ containing the diagonal s.t. the exponential map induced by the metric defines a diffeomorphism 
\begin{align} \label{EqCharts0}
 (\pi^{TN},\exp) : TN \supset V &\map U \subset N \times N
\end{align}
A chart near $f \in C^\infty(M,N)$ may then be defined by
\begin{align} \label{EqCharts1}
 u_f : U_f &\rightarrow V_f, \quad u_f(g) := (x \mapsto (x, (\pi^{TN},\exp)^{-1}\circ (f,g)(x))) 
\end{align} 
where
\begin{align} \label{EqCharts2}
 U_f &:= \{g \in C^\infty(M,N) \mid (f,g)(M) \subset U, f|_{M\setminus K} =  g|_{M\setminus K} \text{ for some compact } K \subset M\} \\
 V_f &:= \{s \in \Gamma_c(f^\ast TN) \mid s(M) \subset (\pi^\ast_{TN} f)^{-1}(V)\}    \notag
\end{align}
where $\pi^\ast_{TN}f: f^\ast TN \rightarrow TN$ denotes the canonical fibrewise isomorphism. The following theorem is proven in \cite{KM} (42.1 for the manifold part, 42.17 and its proof for the statements concerning the tangent bundle):

\begin{Thm} \label{ThmCInftyMN}
$(U_f,u_f)_{f\in C^\infty(M,N)}$ defines a smooth atlas for $C^\infty(M,N)$ and the resulting smooth structure does not depend on the choice of the Riemannian metric. The tangent bundle of this manifold can be canonically identified with $C^\infty_c(M,TN)$ by
\begin{align}\label{EqIdentTangentBundle}
 \Phi : TC^\infty(M,N) &\rightarrow C^\infty_c(M,TN) & \Phi(\dot{c}(0)) = (x \mapsto \tfrac{\del}{\del t} c^\wedge(t,x)|_{t=0})
\end{align}
\end{Thm}

Here, $C^\infty_c(M,TN)$ denotes the space of smooth maps $M \rightarrow TN$ taking the value zero (in the fibres of $TN$) outside of compact sets in $M$. Loosely speaking, \eqref{EqIdentTangentBundle} says that instead of obtaining a vector tangent to $C^\infty(M,N)$ by differentiating a smooth curve $c : \R \rightarrow C^\infty(M,N)$, one may compute the derivative of the induced curve $c^\wedge : \R \times M \rightarrow N$ at each $x \in M$ individually and identify the resulting map with this tangent vector.

\begin{Rem} \label{RemFrechetMf}
In case $M$ is compact, the spaces $\Gamma_c(f^\ast TN) = \Gamma(f^\ast TN)$ are Fr\'{e}chet spaces and $C^\infty(M,N)$ becomes a Fr\'{e}chet-manifold (see \cite{KM} 42.3 or \cite{Ham}, 4.1.3). By the previous theorem, the tangent bundle can now be identified with $C^\infty(M,TN)$.
\end{Rem}

\begin{Rem} \label{RemAlgPropMapping}
In finite-dimensional differential geometry, it is possible to reformulate many concepts using algebraic notions associated with algebras of smooth functions. This is still true for the spaces $C^\infty(M,N)$ ($M$ compact or not). Taking into account that they are Lindel\"{o}f and smoothly paracompact by \cite{KM} 42.3, these spaces enjoy the following properties: 
\begin{enumerate}
 \item Elements of $C^\infty(M_i,N_i)$ ($i=1,2$) bijectively correspond to unital algebra homomorphisms $C^\infty(C^\infty(M_i,N_i),\R) 
       \rightarrow \R$ (use \cite{KM} 17.2 to show that $C^\infty(M_i,N_i)$ is    
       realcompact, then apply 17.6. See also \cite{KMS} 3.3).
 \item Smooth maps $C^\infty(M_1,N_1) \rightarrow C^\infty(M_2,N_2)$ bijectively correspond to unital algebra homomorphisms $C^\infty(C^\infty(M_2,N_2),C^\infty(M_1,N_1)) \rightarrow \R$ (use (a) and \cite{KM} 
       31.4).
 \item $T_fC^\infty(M_i,N_i)$ bijectively correspond to derivations $C^\infty(C^\infty(M_i,N_i),\R)_f \rightarrow \R$ of germs of functions at $f$ (Use \cite{KM} 28.7. The relevant model spaces have the 
       approx. property by \cite{KM} 6.14. They are reflexive in the bornological sense by 6.4 (4), since they are Montel and hence reflexive in the locally convex sense).
\end{enumerate} 
\end{Rem}

We conclude this short description of manifolds of mappings by restating Lemma \ref{LemSmoothCurvVS} for $C^\infty(M,N)$ (see \cite{KM}, 42.5):

\begin{Lemma} \label{LemSmoothCurvCInftyMN}
A curve $c : \R \rightarrow C^\infty(M,N)$ is smooth iff $c^\wedge \in C^\infty(\R \times M,N)$ and in addition, for each interval $[a,b] \subset \R$, there exists $K \subset M$ compact s.t. $c^\wedge(\cdot,x)$ is constant on $[a,b]$ for every $x \in M\setminus K$. The additional condition is clearly void if $M$ is compact.
\end{Lemma}

\vspace{2 \baselineskip}

In chapter \ref{ChapSmooth}, we will look at the special case when the manifold in the target is the total space of a vector bundle. The next paragraphs provide some basic geometry on these spaces. Let $\pi : A \rightarrow M$ be a vector bundle with bundle metric $\langle,\rangle$ and a metric connection $\LC$. $\LC$ gives rise to a decomposition of the tangent bundle of $A$ (see \cite{Lee}, chapter 12 or \cite{KM}, chapter 37 for a detailed exposition),
\begin{align} \label{EqDecompTS}
TA &= \mathcal{H} A \oplus \mathcal{V} A = \mathcal{H} A \oplus \ker(d\pi)
\end{align}
into the subbundles of horizontal and vertical vectors respectively. The corresponding projectors will also be denoted by $\mathcal{H}$ and $\mathcal{V}$. It is well known (\cite{Lee} chapter 12.4 or \cite{KM} 29.9 and 37.2) that these subbundles can be identified in the following way:
\begin{align} \label{EqIdentHV}
\mathcal{H} A &\cong \pi^\ast TM & \mathcal{V} A &\cong \pi^\ast A 
\end{align}
These identifications induce horizontal and vertical lifts
\begin{align*}
h : \Gamma(TM) &\rightarrow \Gamma(\mathcal{H} A) & v : \Gamma(A) &\rightarrow \Gamma(\mathcal{V} A)
\end{align*}
and the horizontal lift $v$ is a right inverse to $d\pi$. Note that the vertical bundle and the corresponding lift are canonically defined whereas their horizontal counterparts as well as both projections $v,h$ depend on the choice of $\LC$. Let now $g$ be a Riemannian metric on $M$ (which may be assumed to be complete). We can define an Riemannian metric on $A$ by declaring the horizontal and the vertical subbundles to be orthogonal and using $g$ and $\langle,\rangle$ on horizontal and vertical tangent spaces, respectively. More precisely, for $\xi_1,\xi_2 \in T_a A$ 
\begin{align} \label{EqDefMetA}
G(\xi_1,\xi_2) &:= g(d\pi(\xi_1), d\pi(\xi_2)) + \langle v^{-1}(\mathcal{V}\xi_1), v^{-1}(\mathcal{V}\xi_2) \rangle 
\end{align}
Let $U_p$ be a geodesically convex neighbourhood of $p \in M$. For each $q \in U$, let $P_{p,q} : A_p \rightarrow A_q$ denote the parallel transport determined by $\LC$ along the unique geodesic joining $p$ and $q$. This map induces the following local trivializations of $A$: 
\begin{align} \label{EqLocIsom}
\Psi_p : A|_{U_p} &\rightarrow U_p \times A_p, \qquad e \mapsto (\pi(e), P_{\pi(e),p}^{-1}(e)) 
\end{align}
It is straightforward to check that its differential at $a \in A|_{U_p}$ is given by 
\begin{align}\label{EqDLocIsom}
d\Psi_p|_a(\xi) &= d\pi|_{\pi(a)}(\mathcal{H}\xi) + P_{\pi(a),p}(\mathcal{V}\xi)   
\end{align}
where $\mathcal{H}\xi \in \mathcal{H}_a A \cong T_{\pi(a)}M$ and $\mathcal{V}\xi \in \mathcal{V}_a A \cong A_{\pi(a)}$. Equipping $U_p \times A_p$ with the product metric $g|_{U_p} + \langle,\rangle_p$ and taking into account that parallel transport preserves the bundle metric, it follows from \eqref{EqDLocIsom} and \eqref{EqDefMetA} that \eqref{EqLocIsom} is a local isometry. Thus, the Riemannian manifold $(A,G)$ is locally isometric to $(M,g) \times (A_p, \langle,\rangle_p)$. Since the Riemann exponential map on $U_p \times A_p$ at $(q,a)$ is given by
\begin{align*}
\exp_{(q,a)}: T_{(q,a)}(U_p \times A_p) \cong T_qM \oplus A_p &\rightarrow U_p \times A_p, &
\exp_{(q,a)}(X,Y) &= (\exp^g_q(X), a + Y),
\end{align*}
the exponential map on $(A,G)$ and its differential is readily obtained by taking compositions with \eqref{EqLocIsom} and  \eqref{EqDLocIsom}:
\begin{align} \label{EqExpMapBun}
\exp_a^G(\xi) &= P_{\pi(a),\exp^g_{\pi(a)}(d\pi(\xi))}(a+\mathcal{V}\xi)  \\
d\exp_a^G|_{\xi}(\zeta) &= h\bigl(d\exp_{\pi(a)}|_{d\pi(\xi)}(d\pi(\zeta))\bigr) + P_{\pi(a),\exp^g_{\pi(a)}(d\pi(\xi))}(\mathcal{V}\zeta) \notag
\end{align}
where we identified $T_{\xi}T_a A \cong T_aA$ for $\xi \in T_a A$. Using this explicit form of the exponential map, we obtain:

\begin{Lemma} \label{LemConvNghdA}
Let $V \subset TM$ ($U \subset M \times M$) be an open neighbourhood of the zero section (the diagonal) s.t. $(\pi^{TM},\exp^g) : V \rightarrow U$ is a diffeomorphism. Then, $(\pi^{TA},\exp^G) : (d\pi)^{-1}(V) \rightarrow (\pi \times \pi)^{-1}(U)$ is a diffeomorphism.   
\end{Lemma}

In particular, we may choose the neighbourhoods appearing in \eqref{EqCharts0} for $(A,G)$ to be products of open sets of the base and the entire fibres.\\

Next, we consider the bundle $\Jet^k(A)$ of $k$-jets of section in $A$. A detailed exposition of the theory of jets can be found in \cite{Sau1} or \cite{Pom}. Here, we only recall that $\Jet^k(A) = \{\jet^k_p(\omega) \mid  p \in M, \ \omega \text{ section in $A$ near } p \}$, where the $k$-th jet $\jet^k_p(\omega)$ of $\omega$ at $p$ is the equivalence class of sections, whose derivatives at $p$ up to order $k$ coincide with those of $\sigma$. It can be shown (see \cite{Sau1}, 6.2.7) that the obvious projection $\Jet^k(A) \smap M$ defines again a smooth vector bundle over $M$. The map $\jet^k : \Gamma(A) \rightarrow \Gamma(\Jet^k(A))$ then assigns to each section its $k$-th jet. Moreover, we clearly have a canonical projection $\Jet^k(A) \rightarrow \Jet^{k-1}_p(A)$ and it is well known (see \cite{Pal} p.58 bottom), that it fits into the following exact sequence of vector bundles:
\begin{align} \label{EqJetExSeq}
\xymatrix{ 0 \ar[r] & \Sym^k T^\ast M \otimes A \ar[r] & \Jet^k(A) \ar[r] & \Jet^{k-1}(A) \ar[r] & 0} 
\end{align}
Since there is no canonical splitting of this exact sequence, there is no canonical way of identifying $\Jet^k(A)$ with $\Sym^{\leq k} T^\ast M  \otimes A$. Thus, there is no canonical way to identify the abstract object $\jet^k_p(\omega) \in \Jet^k(A)$ with a Taylor polynomial of $\omega$ of order $\leq k$, which would be an element of $\Sym^{\leq k} T^\ast M \otimes A$.\\

We will need jet bundles because they give rise to a geometric characterisation of differential operators. If $B \rightarrow M$ is another vector bundle, a linear differential operator $A \rightarrow B$ of order $\leq k$ is an $\R$-linear map $D : \Gamma(A) \rightarrow \Gamma(B)$ s.t the following diagram commutes:
\begin{equation} \label{EqDOJetDiag}
\xymatrix{ \Gamma(A) \ar[r]^{\jet^k} \ar[rd]^{D} & \Gamma(\Jet^k(A)) \ar[d]^{\mathrm{sb}(D)} \\
                                                & \Gamma(B)
         } 
\end{equation}
The $C^\infty_M$-linear map $\mathrm{sb}(D)$ is called the (total) symbol of the differential operator $D$. Thus, denoting the space of differential operators $A \rightarrow B$ up to order $k$ by $\Diff^k(\Gamma(A),\Gamma(B))$, we have bijections
\begin{align} \label{EqJetFactorProp}
\Diff^k(\Gamma(A),\Gamma(B)) &\cong  \Hom_{C^\infty_M}(\Gamma(\Jet^k(A)),\Gamma(B)) \cong \Gamma(\Hom_\R(\Jet^k A, B))
\end{align}
Here, the first bijection is given by taking the symbol, the second one follows from the fact that any $C^\infty_M$-linear map on the level of sections is induced by a morphism of bundles and vice versa.\\

In chapter \ref{ChapAlg}, we will need a $k-th$ order differential for $A$, i.e. a splitting map for the sequence \eqref{EqJetExSeq}, to analyse certain maps between jet spaces in more detail. This map can be obtained as follows: Fixing connections $\LC$ on $A$ and $TM$ (the latter one assumed to be torsion-free) and denoting the induced connections on the various tensor bundles by $\LC$ again, we can form the symmetrized covariant derivative
\begin{align*}
D : \Gamma(\Sym^k T^\ast M \otimes A) &\map \Gamma(\Sym^{k+1} T^\ast M \otimes A), & D(\sigma) := \tsum_{i=1}^n v^i \cdot \LC^{\Sym^k T^\ast M \otimes A}_{v_i} \sigma   
\end{align*}
Here, $k \in \N_0$, $\{v_i\}$ is a local frame of $TM$ with dual frame $\{v^i\}$ and $\cdot$ denotes the product in $\Sym(T^\ast M)$. The $m$-fold iterated derivative $D^m : \Gamma(A) \rightarrow \Gamma(\Sym^m T^\ast M \otimes A)$ then clearly depends only on derivatives of a section $\sigma \in \Gamma(A)$ up to order $m$. Hence, we obtain a map 
\begin{align} \label{EqMthCovDer}
S^k : \Jet^k(A) &\map \Sym^{\leq k} T^\ast M \otimes A, 
& S^k(\jet^k_x(\sigma)) &:= \tsum_{m=0}^k \tfrac{1}{m!} (D^m\sigma)(x)   
\end{align}

By a slight abuse of notation, we will also use the symbol $S^k$ to denote the map \linebreak $\Gamma(A) \rightarrow \Gamma(\Sym^{\leq k} T^\ast M \otimes A), \ \sigma \mapsto \tsum_{m=0}^k \tfrac{1}{m!} D^m\sigma$. We have the following important statement (see \cite{Pal}, p.90):
\begin{Prop} \label{PropTotalDer}
$D^m / m!$ is a $m$-th order total derivative for $A$. In particular, $S^k$ is an isomorphism of vector bundles which preserves the filtration (given by the degree $m$ and the order of the jet, respectively) and we can identify $\jet^k(\sigma)$ with the $k$-th order Taylor polynomial $S^k(\sigma)$.
\end{Prop}

Now assume that $A$ is equipped with the structure of an algebra bundle. Then, for each $k \in \N$, $\Sym^{\leq k} T^\ast M \otimes A$ also carries the structure of an algebra bundle. Its fibrewise multiplication $\cdot_k$ is induced by that of $A$ and truncated multiplication on $\Sym^{\leq k} T^\ast M $, i.e. for $\sigma_i \otimes a_i \in \Sym^{\leq k} T^\ast_x M \otimes A_x$ ($i=1,2$), we have
\begin{align} \label{EqAlgStruct}
 (\sigma_1 \otimes a_1) \cdot_k (\sigma_1 \otimes a_2) &= \pr^{\Sym^{\leq k}} (\sigma_1 \cdot \sigma_2) \otimes (a_1 \cdot a_2)
\end{align}
Assuming that the connection on $A$ is compatible with the algebra structure (i.e. it satisfies the Leibniz rule), Proposition \ref{PropTotalDer} can be improved as follows

\begin{Cor} \label{CorIdentJet} 
Let $A$ be an algebra bundle equipped with a compatible connection. Then, for each $k \in \N$, the map $S^k : \Jet^k(A) \rightarrow \Sym^{\leq k} T^\ast M \otimes A$ defined in \eqref{EqMthCovDer} is an isomorphism of algebra bundles w.r.t. the algebra structure given in \eqref{EqAlgStruct}.
\end{Cor}

\begin{Proof}{}
It remains to show that $S^k$ is multiplicative. $D$ inherits the Leibniz rule from $\LC$. Thus, we have the generalized Leibniz rule for expressions of the form $D^l(\sigma\tau)$ ($\sigma,\tau \in \Gamma(A)$) and obtain the required identity: 
\begin{align*}
\tsum_{l=0}^k \tfrac{1}{l!} D^l(\sigma\tau) 
&= \tsum_{l=0}^k \tfrac{1}{l!} \tsum_{j=0}^l \tbinom{l}{j} D^j(\sigma) \cdot_k D^{l-j}(\tau) 
 = \tsum_{l_1 + l_2 \leq k} \tfrac{1}{l_1!} D^{l_1}(\sigma)  \cdot_k \tfrac{1}{l_2!} D^{l_2}(\tau) \\
&= \bigl( \tsum_{l_1=0}^k \tfrac{1}{l_1!} D^{l_1}(\sigma) \bigr) \cdot_k \bigl( \tsum_{l_2=0}^k \tfrac{1}{l_2!} D^{l_2}(\sigma) \bigr)
\end{align*}
In the last step, we used \eqref{EqAlgStruct}: The product of $\tsum_{l_1=0}^k \tfrac{1}{l_1!} D^{l_1}(\sigma)$ and $\tsum_{l_2=0}^k \tfrac{1}{l_2!} D^{l_2}(\tau)$ in $\Sym^{\leq k} T^\ast M$ is obtained by computing the product in $\Sym(T^\ast M)$ and dropping all terms of order $> k$.  
\end{Proof}

In the remaining part of this section, we discuss chain rule-type formulas for $S^k$ under composition with smooth maps. Let $A \rightarrow M$ and $B \rightarrow N$ be vector bundles equipped with connections. For $\varphi \in C^\infty(M,N)$, the bundles $\varphi^\ast B$ and $\Hom(\varphi^\ast B,A)$ over $M$ carry induced connections. All connections are denoted by $\LC$; occasionally we write $\LC^\varphi$ to indicate that a connection has been pulled back along $\varphi$. Given $\phi \in \Gamma(\Hom(\varphi^\ast B, A))$ (or, equivalently, a vector bundle morphism $\phi: \varphi^\ast B \rightarrow A$), we have induced maps 
\begin{align} \label{EqConcatMor}
\xymatrix{ \Gamma(N,B) \ar[r]^{\varphi^\ast} & \Gamma(M,\varphi^\ast B) \ar[r]^{\phi_\ast} & \Gamma(M,A)} 
\end{align}
The definition of $S^k$ in \eqref{EqMthCovDer} clearly extends to sections of arbitrary vector bundles. In particular, we have
\begin{align*}
S^k(\phi) &= \tsum_{l=0}^k \tfrac{1}{l!} D^l \phi \quad \in  \Gamma(\Sym^{\leq k} T^\ast M \otimes \Hom(\varphi^\ast B,A))
\end{align*}

\vspace{\baselineskip}

The composition $\circ: \Hom(\varphi^\ast B,A) \otimes \varphi^\ast B \rightarrow A$, viewed as a morphism of vector bundles, naturally extends to
\begin{align} \label{EqCompProd1}
\circ_k: \Bigl( \Sym^{\leq k} T^\ast M \otimes \Hom(\varphi^\ast B,A) \Bigr) \otimes \Bigl( \Sym^{\leq k} T^\ast M \otimes \varphi^\ast B  \Bigr) & \rightarrow \Sym^{\leq k} T^\ast M \otimes A \\
(\sigma_1 \otimes b) \otimes (\sigma_2 \otimes \psi) &\mapsto \pr^{\leq k}(\sigma_1 \sigma_2) \otimes \psi(b). \notag
\end{align}
Hence, we may compose $S^k(\phi)$ with elements or sections of $\Sym^{\leq k} T^\ast M \otimes \varphi^\ast B$. It is easy to verify that the induced connections satisfy the Leibniz rule w.r.t. the product $\circ_k$. By using the argument from the proof of Corollary \ref{CorIdentJet}, we obtain

\begin{Lemma} \label{LemVBHomChainRule}
For $k \in \N$, $\phi \in \Gamma(\Hom(\varphi^\ast B,A))$ and $b \in \Gamma(\varphi^\ast B)$, we have 
\begin{align} \label{EqVBHomChainRule}
S^k(\phi \circ b) &= S^k(\phi) \circ_k S^k(b), 
\end{align}
and we may view $S^k(\phi)$ as a map $\Sym^{\leq k} T^\ast M \otimes \varphi^\ast B  \rightarrow \Sym^{\leq k} T^\ast M \otimes A$.
\end{Lemma}

\begin{Rem}
Since $\Hom(\varphi^\ast B,A) \cong \varphi^\ast B^\ast \otimes A$, we may view $S^k(\phi)(x) \in \Sym^{\leq k} T^\ast_x M \otimes B^\ast_{\varphi(x)} \otimes A_x$ as a polynomial on $T_x M \oplus B_{\varphi(x)}$ with values in $A_x$. It has only contributions of degree one in the directions of the fibres of $B$ which reflects the fact that $\phi$, being a morphism of vector bundles, is fibrewise linear. Thus, we may interpret formula \eqref{EqVBHomChainRule} as the \emph{chain rule} for the $k$-th derivative of the composition of $b$ and $\phi$.
\end{Rem}

\begin{Rem}\label{RemFunkSK}
Given a third bundle $B' \rightarrow N'$ and a map $\varphi' : N \rightarrow N'$, we may compose bundle morphisms $\phi \in \Hom(\varphi^\ast B, A)$ and $\phi' \in \Hom(\varphi'{}^\ast B', B)$. In analogy to Lemma \ref{LemVBHomChainRule}, we find $S^k(\phi' \circ \phi) = S^k(\phi') \circ_k S^k(\phi)$. Moreover, $S^k(\mathrm{Id}_A) = \mathrm{Id}_A$.
\end{Rem}

It is slightly more complicated to derive a similar result for the nonlinear map $\varphi \in C^\infty(M,N)$. We define symmetrized derivatives recursively (for $k>1$) by 
\begin{align*}
D^0(\varphi) &:= \varphi &
D^1(\varphi) &:= d\varphi &
D^k(\varphi) &:= D^{k-1}d\varphi 
\end{align*}
where $D^{k-1}$ on the right hand side simply denotes the $(k-1)$-fold symmetrized covariant derivative on the bundle $T^\ast M \otimes \varphi^\ast TN$. The $k$-th order Taylor polynomial may be defined by 
\begin{align} \label{EqTaylorVarphi}
S^k(\varphi) &:= \tsum_{l=1}^k \tfrac{1}{l!} D^l\varphi \in \Gamma(\bigoplus_{l=1}^k \Sym^l T^\ast M \otimes \varphi^\ast TN) 
\end{align}

\begin{Rem} \label{RemFlatTaylor}
For $x \in M$, $S^k(\varphi)(x)$ is a polynomial $T_{x}M \rightarrow T_{\varphi(x)}N$ of order $\leq k$ without constant term. In case that $M$ and $N$ are vector spaces equipped with the canonical flat connection, we may choose a basis $\{v_i\}_{i=1,\ldots,m}$ on $M$ with dual basis $\{v^i\}$  and write 
\begin{align} \label{EqTaylorVarphiFlat}
S^k(\varphi) &= \tsum_{0 \neq I \in \N_0^{m,k}} \tfrac{1}{I!} v^I D_I \varphi  
\end{align}
where $D_I \varphi := \langle v_I, D^{|I|} \varphi \rangle / |I| !$ coincides with the ordinary coordinate derivative. Thus, \eqref{EqTaylorVarphiFlat} is just the ordinary Taylor expansion (with variable reference point $x \in M$) whose $0$-th order contribution $\varphi(x)$ at $x$ is encoded in the base point of the target $T_{\varphi(x)}N$. We will use \eqref{EqTaylorVarphiFlat} locally on arbitrary manifolds $M,N$ with non-flat connections; the expression represents the full Taylor polynomial of $\varphi \in C^\infty(M,N)$ up to order $k$.
\end{Rem}

Note that we can not expect to obtain an analogue to Proposition \ref{PropTotalDer} for $\Jet^k(M,N)$ since the bundle in \eqref{EqTaylorVarphi} depends on $\varphi$ and in general, $\Jet^k(M,N)$ is not even a vector bundle. Since $S^k(\varphi)$ is meaningful for each $\varphi$ individually, we may yet compose sections of $B$ with $\varphi$ as indicated in \eqref{EqConcatMor} and derive a chain rule. Using $C^\infty(M,B) = \bigsqcup_{\varphi \in C^\infty(M,N)}\Gamma(\varphi^\ast B)$, the composition $\circ : \Gamma(B) \times C^\infty(M,N) \rightarrow C^\infty(M,B)$ induces for each $\varphi \in C^\infty(M,N)$ a map 
\begin{align*}
\circ_k :  (\varphi^\ast \Sym^{\leq k} T^\ast N \otimes \varphi^\ast B) \times_M (\Sym^{\leq k} T^\ast M \otimes \varphi^\ast TN)  
&\rightarrow \Sym^{\leq k} T^\ast M \otimes \varphi^\ast B
\end{align*}
which is defined by truncated composition of polynomial maps $TM \rightarrow \varphi^\ast TN$, $\varphi^\ast TN \rightarrow \varphi^\ast B$ ((see \cite{D4}, A.21 or \cite{Car} 7.4)). Note that this map can only be defined on the fibred product but not on the tensor product of bundles, because it is not bilinear on the fibres. To formulate a chain rule for higher derivatives of $b \circ \varphi$ where $b \in \Gamma(B)$, we define for $k \in \N$ the following set of multi-indices:
\begin{align*}
\M_k &:= \{ \alpha \in \N_0^k \mid \tsum_{i=1}^k i\alpha_i = k\}, 
\end{align*}
We have the following Faa-di-Bruno-type result (see \cite{Ron} or \cite{Frae} for a discussion of such formulas on vector spaces):

\begin{Prop} \label{PropFaaDiBruno}
For $\varphi \in C^\infty(M,N)$, $b \in \Gamma(B)$ and each $m \in \N$, we have
\begin{align} \label{EqFaaDiBruno}
D^m (b \circ \varphi) &= \sum_{\alpha \in \M_m} \tfrac{m!}{\alpha!} \varphi^\ast (D^{|\alpha|}(b)) \circ \prod_{j=1}^m \bigl( \tfrac{1}{j!} D^j \varphi \bigr)^{\alpha_j}
\end{align}
\end{Prop}

Note that for $\alpha \in \M_m$ fixed, the product appearing on the r.h.s. of \eqref{EqFaaDiBruno} is simply the product in $\Gamma(\Sym^{\leq m} T^\ast M \otimes \varphi^\ast \Sym^{\leq |\alpha|} TN)$. The composition of $\varphi^\ast (D^{|\alpha|}(b))  \in \Gamma(\varphi^\ast \Sym^{\leq |\alpha|} T^\ast N \otimes \varphi^\ast B)$ and the aforementioned product is now simply given by the pairing of $\varphi^\ast \Sym^{\leq |\alpha|} TN$ and its dual. We will only sketch a proof of Proposition \ref{PropFaaDiBruno}: $D$ satisfies the Leibniz rule and the chain rule $D^\varphi(b \circ \varphi) = D(b) \circ d\varphi$; the latter one is inherited from the induced connection $\LC^\varphi$. Hence, we may proceed as in one of the proofs of Faa-di-Bruno's formula for maps between vector spaces (see e.g. \cite{Ron}, section 2) which only relies on these two facts and some combinatorics.\\

A formula similar to \eqref{EqFaaDiBruno} can be obtained for $D^m(\varphi_2 \circ \varphi_1)$  where $\varphi_1 \in C^\infty(M,N)$, \linebreak $\varphi_2 \in C^\infty(N,P)$ (e.g. by viewing maps $N \rightarrow P$ as sections in the trivial fibre bundle $N \times P \twoheadrightarrow N$ and extending the previous statement). We do not give the details but just note, that it allows us to conclude the following statement: 

\begin{Lemma} \label{LemmaFuncSK}
For $\varphi_1,\varphi_2$ as above, we have $S^k(\varphi_1 \circ \varphi_2) = S^k(\varphi_1) \circ_k S^k(\varphi_2)$. Moreover, $S^k(\Id_M) = \Id_{\Gamma(\End(TM))}$ where we identified $\End(TM) \cong T^\ast M \otimes TM \subset \Sym^{\leq k} T^\ast M \otimes TM $.  
\end{Lemma}

In fact, it is well known that this statement is true for ordinary Taylor polynomials (see \cite{D4}, 19.5 or \cite{Car} 7.4), i.e. for the case that the manifolds are vector spaces equipped with the canonical flat connections. But since the ordinary $k$-th order differential as well as $D^k$ satisfy the same decomposition formula \eqref{EqFaaDiBruno}, we may replace the differentials in the formula for ordinary Taylor polynomials by the $D^k$ defined before \eqref{EqMthCovDer} and obtain the lemma.   

\begin{Cor}
A completely analogous statement holds for expressions of the form $b \circ \varphi$ : $S^k(b \circ \varphi) = S^k(b) \circ_k S^k(\varphi)$. In particular, the Taylor polynomial $S^k(\varphi)$ induces a map $\varphi^\ast \Sym^{\leq k} T^\ast N \otimes \varphi^\ast B \rightarrow \Sym^{\leq k} T^\ast M \otimes \varphi^\ast B$.
\end{Cor}


\section{Supergeometry} \label{ChapSuper}

This chapter provides a brief introduction into some concepts of superalgebra and -geometry, covering the ``classical'' ringed space definition as well as the functorial formulation.\\ 

Superalgebra deals with categories of $\Z_2$-graded objects (e.g. rings, modules and their tensor products, \ldots), s.t. algebraic structures as well as morphisms are compatible with the grading. We will not discuss these structures in detail but refer to the literature (e.g. \cite{Ma1}, \cite{Var}, \cite{AMS1}, \cite{CCF} and many others). Even though morphisms in all categories respect the $\Z_2$-grading by definition (e.g. a morphism of supermodules $M_\0 \oplus M_\1 \rightarrow N_\0 \oplus M_\1$ maps $M_i$ to $N_i$, where $i=\0,\1$ denotes the $\Z_2$-grading), we will also need the set of maps respecting the algebraic structures but not the grading. To distinguish them from the ``honest'' morphisms $\Hom(X,Y)$, they will be denoted by $\qq{\Hom}(X,Y)$. In fact, the sets $\qq{\Hom}(X,Y)$ naturally arise as inner Hom functors/objects which means that they are characterized by an adjunction formula (see \cite{GM1} II.4.23 or \cite{Var} 3.7 for the use in superalgebra). In the category of supermodules, this e.g. simply means that $\qq{\Hom}(M,N)$ is characterized by natural bijections $\Hom(M,\qq{\Hom}(N,P)) \cong \Hom(M \otimes N,P)$.\\

We briefly recall the construction of changing the base ring of a module, since this is not discussed in the standard texts. Let $M$ be a $R$-supermodule and $\Phi : R \rightarrow S$ be a morphism of superrings. Then $S$ becomes a $R$-supermodule denoted $S_\Phi$ by $r \cdot s := \Phi(r)s$. We can define the $S$-supermodule, obtained from $M$ by \emph{extension of scalars from $R$ to $S$} by
\begin{align}\label{EqExtScalars}
S \otimes_\Phi M &:= S_\Phi \otimes_R M &  s'\cdot (s \otimes m) := (s's) \otimes m
\end{align}
A number of properties of this construction for ordinary modules is discussed in \cite{BouAlg1}, III.5.3. The proofs can easily be adapted to the case of supermodules.\\

We will use the approach due to Berezin, Kostant and Leites (see \cite{Ber}, \cite{Kos}, \cite{Lei}) as fundamental definition for finite-dimensional supermanifolds:

\begin{Def} \label{DefSMf}
A supermanifold $X$ of dimension $p|q$ is a locally ringed space $(\tilde{X},\O)$ where $\tilde{X}$ is a topological manifold and $\O$ is a sheaf of supercommutative algebras locally isomorphic to $\C^\infty(\R^p) \otimes \tL \R^q$. A morphism of supermanifolds is a morphism of the corresponding ringed spaces. The resulting category will be denoted by $\BKL$. 
\end{Def}

Note that, by definition, a morphism $\Phi: (X,\O) \rightarrow (Y,\cR)$ is given by a pair $(\Phi^\ast,\varphi)$, where $\varphi : \tilde{X} \rightarrow \tilde{Y}$ is continuous, $\Phi^\ast: \cR \rightarrow \varphi_\ast\O$ is a morphism of sheaves of supercommutative algebras over $\tilde{Y}$ and $\varphi_\ast\O$ denotes the direct image sheaf (see e.g. \cite{Ten} 3.7.1). Being morphisms of super algebras, the maps on sections have to preserve parity. As is true for ringed spaces in general, we may equivalently use morphisms of sheaves over $\tilde{X}$ to describe the morphisms in $\BKL$ (see \cite{Ten}, 3.7.11). More important for this work, $\Hom_{\BKL}(X,Y)$ can be characterized in terms of homomorphisms acting on global sections (see \cite{BBHR}, Corollary III.1.5):

\begin{Prop} \label{PropMorGlob}
For all supermanifolds $X,Y \in \BKL$, the following map is a bijection:
\begin{align*}
\Hom_{\BKL}(X,Y) &\rightarrow \Hom_{\mathrm{SAlg}}(\cR,\O), \qquad \Phi \mapsto \Phi^\ast{\tilde{Y}} 
\end{align*}
\end{Prop}

This statement relies on an argument using partition of unity and is not true on holomorphic supermanifolds. The sections of $\O$ and $\cR$ are interpreted as functions on $X$ and $Y$, respectively and often called superfunctions.\\

It is well known that supermanifolds can be described using the language of differential geometry. Denoting the sheaf of nilpotent elements in $\O$ by $\mathcal{J}$, it is not difficult to prove that the quotient sheaf $\O / J$ is locally isomorphic to the sheaf $C^\infty_{\R^p}$. Thus, $\O$ induces on $\widetilde{X}$ the structure of a smooth manifold of dimension $p$. There is a canonical embedding $\iota_X : \tilde{X} \hookrightarrow X$ which is given by the quotient map $\iota_X^\ast :\O \twoheadrightarrow \O / J \cong C^\infty(\tilde{X})$ on the level of sheaves of functions. It can be shown that each morphism $\Phi : X \rightarrow Y$ projects down to a smooth map $\varphi : \tilde{X} \rightarrow \tilde{Y}$. The sheaf of superfunctions $\O$ can always be identified in the following way: 

\begin{Thm}[Batchelor \cite{Bat}, Main Theorem] \label{ThmBachelor}
For each supermanifold $X$ of dimension $p|q$, there is a vector bundle $E \rightarrow \tilde{X}$ of rank $q$, s.t. $X \cong (\tilde{X},\Gamma(\tL E^\ast))$. The vector bundle is uniquely defined up to isomorphism. 
\end{Thm}

We will refer to $E$ as a Batchelor bundle for $X$. Again, the theorem is not valid in the holomorphic category. Note that even though $E$ is uniquely determined up to isomorphism, there is no canonical choice for the isomorphism $X \cong (\tilde{X},\Gamma(\tL E^\ast))$. In fact, the $\Z_2$-grading of the sheaf $\Gamma(\tL E^\ast)$ is induced by its natural $\Z$-grading but $\O$ does not carry a corresponding natural $\Z$-grading. In fact, it can be shown that a choice of an isomorphism is equivalent to the choice of a function factor or the choice of the even part of a coordinate system (\cite{Kos}, Theorem 2.7).   

\begin{Rem}
In \ref{ThmBachelor}, we adopted the convention to identify $\O$ with $\Gamma(\tL E^\ast)$ instead of $\Gamma(\tL E)$. It allows us to embed the category $\VBun$ of vector bundles (over different base spaces) into $\BKL$. Objects of $\VBun$ are just finite-dimensional vector bundles $S\rightarrow M$ and a morphism $\phi: (S_1 \rightarrow M_1) \rightarrow (S_2 \rightarrow M_2)$ is simply a smooth, fibre-preserving and fibrewise linear map $\phi : S_1 \rightarrow S_2$ (which then always projects down to a smooth map $\varphi: M_1 \rightarrow M_2$). The embedding functor $\Gr : \VBun \rightarrow \BKL$ is then defined by
\begin{align*}
(S \rightarrow M) &\mapsto (M,\Gamma(\tL S^\ast)) \\
(\phi : S_1 \rightarrow S_2) &\mapsto (\Gamma(\tL S_2^\ast) \rightarrow \Gamma(\tL S_1^\ast), \sigma \mapsto \sigma \circ \wedge \phi)
\end{align*}
This functor is faithful but not full, because a morphism not preserving the $\Z$-grading is not in its image. By Batchelor's theorem, it is moreover essentially surjective.  
\end{Rem}

\begin{Rem} \label{RemProducts}
It is well known that the category $\BKL$ admits finite products (see e.g. \cite{Ber} 4.4.2 or \cite{CCF} 4.5.3). A geometric description using Batchelor bundles can be given as follows (cf. \cite{CondeGro} 5.21): If $E_i \rightarrow \tilde{X}_i$ are Batchelor bundles for $X_1$ and $X_2$, then a Batchelor bundle of the product is given by $\pr_1^\ast E_1 \oplus \pr_2^\ast E_2 \rightarrow \tilde{X}_1 \times \tilde{X}_2$ where $\pr_i : \tilde{X}_1 \times \tilde{X}_2 \twoheadrightarrow \tilde{X}_i$ denote the canonical projections.
\end{Rem}

$\BKL$-supermanifolds come along with the usual tangent structure. Vector fields are (local) sections of the tangent sheaf $\T_X$ which is a sheaf of $\O$-modules over $\tilde{X}$. For $U \subset \tilde{X}$ open, it is defined by $\T_X(U) := \Der(\O(U))$, the set of super derivations of the ring $\O(U)$ of superfunctions. It is locally free of rank $\dim(X) = p|q$ and a local basis is given by coordinate derivatives $\tfrac{\del}{\del x^i}$, $\tfrac{\del}{\del \theta_\alpha}$ (cf. \cite{AMS1} \textsection 3.3). General super vector bundles over $X$ may be defined as locally free $\O$-modules over $\tilde{X}$ (cf. \cite{AMS1} \textsection 3.2). Associated to $\T_X$, there is the sheaf $\Gamma(TX)$ over $\tilde{X}$ defined by $\Gamma(U,TX) := \T_X(U)/ (\J(U)\T_X(U))$, i.e. its sections are given by sections of $\T_X$, composed with the ``evaluation map'' $\iota^\ast_X : \O \rightarrow C^\infty(\tilde{X})$. Due to the quotient construction, this is only a sheaf of $C^\infty(\tilde{X})$-modules. We have (see \cite{Schmitt1}, p.128, and \cite{Kos}, section 2.10)

\begin{Prop}
The sheaf $\Gamma(TX)$ is locally free of rank $p|q$, i.e. the sheaf of sections of a real, $\Z_2$-graded vector bundle $TX = TX_\0 \oplus TX_\1$ over $\tilde{X}$. We have $TX_\0 \cong T\tilde{X}$ and, after the choice of a Batchelor bundle $E$, also $TX_\1 \cong E$.  
\end{Prop}

We close the discussion of the BKL picture with a brief discussion of pullbacks. Let $\Phi=(\varphi,\Phi^\ast) \in \Hom_\BKL(X,Y)$ and consider the tangent sheaf $\T_Y$ (the construction works for arbitrary finite-dimensional super vector bundles). Following \cite{Ten} 4.4.13, we may form the sheaf $\Phi^\ast \T_Y$ over $\tilde{X}$, defined by $\Phi^\ast\T_Y(U) := \O(U) \otimes_\Phi (\varphi^\ast \T_Y)(U)$ (cf. \eqref{EqExtScalars} for the notion $\otimes_\Phi$). It is not hard to show that it is locally free, too. If $\{\xi^\nu\}$ are local coordinates on $Y$, a local basis for $\Phi^\ast \T_Y$ is given by $\{\Phi^\ast \tfrac{\del}{\del \xi^\nu} = 1 \otimes_\Phi \tfrac{\del}{\del \xi^\nu} \}$. By 
\begin{align} \label{EqPullBack}
\bigl( \Phi^\ast \tfrac{\del}{\del \xi^\nu} \bigr) g &= (1 \otimes_\Phi \tfrac{\del}{\del \xi^\nu}) g = 1 \otimes_\Phi \tfrac{\del g}{\del \xi^\nu} \cong \Phi^\ast(\tfrac{\del g}{\del \xi^\nu})
\end{align}
these sections act on superfunctions $g$ on $Y$ (or more general, on sections of $\varphi^\ast \cR$) as derivations. In the last step, we used the identification $\O(U) \otimes_\Phi (\varphi^\ast \cR)(U) \cong \O(U)$ induced by \eqref{EqExtScalars}.

\vspace{2 \baselineskip}
 
We can not apply Definition \ref{DefSMf} to construct a ''supermanifold of morphisms`` $\qq{\SC}^\infty(X,Y)$ since $\BKL$ does not contain infinite-dimensional objects. Following ideas of Molotkov (\cite{Mol}) and others, we make use of the ``functor of points'' from algebraic geometry (see \cite{EH1} for the general ideas and \cite{AMS1} 2.8-2.9, \cite{Var} 4.5 or \cite{CCF} 3.4 \& 4.6 for its applications in supergeometry ). It roughly says that instead of studying a supermanifold $X$ directly, one may study the functor $S \mapsto \Hom_\BKL(S,X)$ which assigns to an arbitrary finite-dimensional ``test''-supermanifold $S \in \BKL$ the set  $\Hom_{\BKL}(S,X)$ of so called $S$-points. Similarly, morphisms are replaced by natural transformations among these functors and Yoneda's lemma allows us to reconstruct finite-dimensional supermanifolds and their morphisms. The infinite-dimensional object $\qq{\SC}^\infty(X,Y)$ can be \emph{defined} this way by equipping the sets 
$\Hom_\BKL(S,\qq{\SC}^\infty(X,Y))$ with suitable extra structures.\\
It is sometimes said that the ringed space approach is not appropriate for the construction of infinite dimensional spaces. We already observed in Remark \ref{RemAlgPropMapping} that it is in fact possible to characterise important geometric objects associated to $C^\infty(M,N)$ in terms of its algebra of smooth functions. Recent work by Egeileh and Wurzbacher \cite{EW} (4.12 - 4.14) shows that for a large class of model spaces (which unfortunately does not contain $\Gamma_c(S)$ from \eqref{EqModell}), the notion of a reduced infinite-dimensional manifold can be defined entirely in the language of ringed spaces. They moreover extend their results in 4.24, 4.25 to Douady's concept (cf. \cite{Dou}, section 3.2) of functored spaces, which is closer to the ringed space language than the concept of the functor of points. Functored spaces were also used by Alldridge (\cite{All}) to define a category of supermanifolds containing infinite-dimensional objects. We use the functorial picture since we will see in chapter \ref{ChapAlg} that the sets of higher points of mapping spaces can be interpreted in a very geometric way leading to a natural choice of a smooth structure on them.  Nevertheless, it is an  interesting task to reformulate (parts of) the results presented here in the language of \cite{All} and \cite{EW}; it will be addressed elsewhere.\\    

Molotkov's approach (\cite{Mol}), which was worked out in detail by Sachse (\cite{Sa1}, \cite{Sa2}), relies on a smaller class of ``test objects'' $S$ than all of $\BKL$ (cf. Remark \ref{RemAdvMSApproach} for further discussion). To illustrate the idea, note that for each finite-dimensional smooth manifold $M$, we have an obvious bijection $M \cong C^\infty(\R^0,M)$ of sets. It is clear that the $C^\infty$-structure of $M$ is not determined by the set $C^\infty(\R^0,M)$. However, if this set of maps is equipped with the structure of a (finite-dimensional) smooth manifold, then we trivially obtain a smooth manifold structure on $M$. For $X \in \BKL$, it follows from Proposition \ref{PropMorGlob} that $\Hom_\BKL(\R^{0},X) \cong C^\infty(\R^0,\tilde{X})$ since homomorphisms map nilpotents to nilpotents. Thus, studying the set of morphism $\R^{0} \rightarrow X$ is insufficient to study $X$ since only information about the underlying manifold $\tilde{X}$ is retained. It was proposed by several authors (\cite{Lei} who quotes J. Bernstein, \cite{Shv1}, \cite{Shv2}, \cite{Vor} and others) to study sets $\Hom_{\BKL}(\PR_q,X)$ for arbitrary $q \in \N_0$, where $\PR_q \cong \R^{0|q}:= (\{0\},\tL \R^q)$ is a superpoint, i.e. an element of $\BKL$ whose underlying manifold consists of one point. This formally leads to consider functors from the category $\mathrm{SPoint}{\mathrm{opp}}$ of superpoints into a suitable category of manifolds (sets, rings, \ldots) as supermanifolds (supersets, superrings, \ldots). Molotkov used this approach to give an alternative definition for supermanifolds which coincides with \ref{DefSMf} above in the finite-dimensional case (cf. discussion before \eqref{EqAequiCat}) but which is also meaningful in the infinite-dimensional setting. We briefly summarize those parts of his theory relevant for this article.\\ 

Let $\Gr$ denote the category of finite-dimensional Grassmann algebras. It follows from Definition \ref{DefSMf} (see also \cite{Sa2}, Proposition 2.8) that the functor $\PR : \Gr \rightarrow \mathrm{SPoint}^{\op}, \ \Lambda \mapsto (\{pt\}, \Lambda)$ defines a natural equivalence. Thus, for any category $\mathcal{C}$, we can replace $\mathcal{C}^{\mathrm{SPoint}^{\op}}$ by $\mathcal{C}^\Gr$. Given $F \in \mathcal{C}^{\Gr}$, $F(\Lambda)$ will also be referred to as the $\Lambda$-points (or $\PR(\Lambda)$-points) of $F$. For any super vector space $W$, we consider the functor $\qq{W} \in \Vec^\Gr$ defined by 
\begin{align}\label{EqRepSMod}
 \qq{W}(\Lambda) &:= (\Lambda \otimes W)_\0 &
 \qq{W} (\rho : \Lambda \rightarrow \Lambda') &:= \rho \otimes \Id_W 
\end{align}
$\qq{W}$ is a $\qq{\R}$-module, i.e. each $\qq{W}(\Lambda)$ is a module over $\qq{\R}(\Lambda) = \Lambda_\0$ and the maps $\qq{W}(\rho)$ are module homomorphisms. Not each $\qq{\R}$-module is of this special form; those arising from super vector spaces as in \eqref{EqRepSMod} are called \emph{superrepresentable} and serve as model spaces. 

To discuss the required coordinate neighbourhoods for supermanifolds in the functorial framework, recall that $F'\in \Top^\Gr$ is called an open subfunctor of $F \in \Top^\Gr$ (denoted $F' \subset F$), if $F'(\Lambda) \subset F(\Lambda)$ is an open subset for all $\Lambda \in \Gr$ and the inclusions $i_\Lambda : F'(\Lambda) \hookrightarrow F(\Lambda)$ form a morphism in $\Top^\Gr$. The relevant notion (which we will slightly extend in Definition \ref{DefLocAffDom} below) is now given by

\begin{Def}[Molotkov-Sachse]
Let $W$ be a super vector space s.t. $\qq{W} \in \mathrm{cVec}^\Gr$. Then a subfunctor $D \subset \qq{W}$ is called a (convenient) superdomain, if it is an open subfunctor w.r.t the $c^\infty$-topology on each $\qq{W}(\Lambda)$.
\end{Def}
Note that by \cite{KM}, 2.15, the condition $\qq{W} \in \mathrm{cVec}^\Gr$ is equivalent to the requirement $W = W_\0 \oplus W_\1 \in \mathrm{cVec}$. Using the functorial properties, it can be shown that superdomains have a very special form (see \cite{Sa2} Proposition 4.8): They are always restrictions, i.e. functors of the form 
\begin{align} \label{EqRestriction}
 \qq{W}|_V(\Lambda) &:= \qq{W}(\pr^{\Lambda}_\R)^{-1}(V) 
\end{align}
where $V \subset \qq{W}(\R) = W_\0$ is some open set and $\pr^{\Lambda}_\R$ the canonical projection $\Lambda \twoheadrightarrow \R$. Supersmooth maps between two superdomains are now defined as follows:

\begin{Def}[Molotkov-Sachse] \label{DefSuperSmoothDom}
Let $D_i \subset \qq{V}_i$ be two convenient superdomains. A morphism $f : D_1 \rightarrow D_2$ in $\TGr$ is called supersmooth if
\begin{enumerate}
 \item For each $\Lambda \in \Gr$, the map $f_\Lambda : D_1(\Lambda) \rightarrow D_2(\Lambda)$ is a smooth map between $c^\infty$-open sets of $\qq{V}_1(\Lambda)$ and $\qq{V}_2(\Lambda)$ in the sense of Definition \ref{DefSmoothMap}.
 \item For each $\Lambda \in \Gr$, the differential $df_\Lambda : D_1(\Lambda) \times \qq{V}_1(\Lambda) \rightarrow \qq{V}_2(\Lambda)$  is $\Lambda^{\ev}$-linear in the second entry. 
\end{enumerate}
\end{Def}

Clearly, one may replace $\mathrm{cVec}$ by other categories of topological vector spaces allowing for a suitable concept of differential. To be able to give a precise meaning to notions like chart and atlas, we need to consider functors $G,F \in \Top^\Gr$ more general than superdomains and morphisms between them. Recall that that a morphism $\alpha : G \rightarrow F$ of functors is called open, if it factorises as $\alpha = i \circ \alpha'$, where $\alpha' : G \overset{\sim}{\rightarrow} F'$ is an isomorphism of functors onto an open subfunctor $F' \subset F$ and $i$ is the inclusion $F' \hookrightarrow F$. 

\begin{Def}[Molotkov-Sachse] \label{DefSMfFunc}
Let $X$ be a functor in $\Man^\Gr$. A collection of open morphisms $\mathcal{A} := \{u_i: V_i \overset{\sim}{\rightarrow} U_i \subset X \}_{i \in I}$ is called an atlas on $X$ provided that
\begin{enumerate}
 \item $\mathcal{A}$ is an open covering of $X$ by convenient superdomains $V_i$, i.e. the open sets $U_i(\Lambda)$ form an open covering of $X(\Lambda)$ in the usual sense for each $\Lambda \in \Gr$.
 \item For all $i,j \in I$, the fibre product $V_{ij} := V_i \times_X V_j \in \Man^\Gr$ carries the structure of a convenient superdomain, s.t. the canonical projections $\Pi_i : V_{ij} \rightarrow V_i$, $\Pi_j: V_{ij} \rightarrow V_j$ are supersmooth in the sense of Definition \ref{DefSuperSmoothDom}.   
\end{enumerate}
Two atlases $\mathcal{A}_1,\mathcal{A}_2$ on $X$ are called equivalent, if $\mathcal{A}_1 \cup \mathcal{A}_2$ is again a supersmooth atlas. The functor $X$, equipped with such an equivalence class of smooth atlases, is called a (convenient) supermanifold.
\end{Def}

\begin{Rem} \label{RemSuperLin}
The fibred product $V_i \times_X V_j$ in the second condition exists as an element of $\Man^\Gr$; it is given by $V_{ij}(\Lambda) = ((u_i)_\Lambda \times (u_j)_\Lambda)^{-1}(U_i(\Lambda) \cap U_j(\Lambda)) \subset V_i(\Lambda) \times V_j(\Lambda)$ (see e.g.\cite{GM1} II.3.5 for the general definition of pullbacks) with $\Pi_i,\Pi_j$  given by the projections on the left and right factor. The commutative diagram
\begin{align*}
\xymatrix{ & U_{ij}(\Lambda) \ar[ld]_{(\Pi_i)_\Lambda} \ar[rd]^{(\Pi_i)_\Lambda} & \\ 
          (u_i)_\Lambda^{-1}(U_i(\Lambda) \cap U_j(\Lambda)) \ar[rr]^{(u_j)_\Lambda^{-1}(u_i)_\Lambda} & & (u_j)_\Lambda^{-1}(U_i(\Lambda) \cap U_j(\Lambda)) }
\end{align*}
now shows that the second condition of the previous definition is satisfied if and only if all chart transitions $(u_j)_\Lambda^{-1}(u_i)_\Lambda$ are supersmooth in the sense of Definition \ref{DefSuperSmoothDom} for each $\Lambda \in \Gr$.
\end{Rem}

\begin{Rem}
Note that not all morphisms in $\Man^\Gr$ are supersmooth (see e.g. \cite{BCF1} 3.15 for a counterexample). It is important to require (b) of Definition \ref{DefSMfFunc} since it implies that differentials of chart transitions are morphisms of $\qq{\R}$-modules, which is essential for the construction of a tangent bundle with fibres of the form \eqref{EqRepSMod} (cf. \cite{Sa1} section 3.8.2). In particular, a suitable notion of chart which allows for the formulation of (b) is crucial in the finite- and infinite-dimensional situation. This is an important reason for requiring the functors to take values in $\Man$ rather than in a more general category of smooth spaces.  
\end{Rem}

\begin{Rem} \label{RemSMfMorphism}
Given two supermanifolds $X,Y$ in the sense of Definition \ref{DefSMfFunc}, the notion of supersmoothness for a morphism $\Phi : X \rightarrow Y$ in $\Man^\Gr$ is defined using charts $u: U \rightarrow X$ and $u' : U' \rightarrow Y$: The fibred product $U \times_Y U' \in \Man^\Gr$ defined w.r.t the morphisms $\Phi \circ u \circ \Pi_X$ and $u' \circ \Pi_Y$ is required to be a superdomain (see \cite{Sa2} Def. 4.14 for details). Similar to the discussion in Remark \ref{RemSuperLin}, this condition is satisfied iff all local representations $(u')_\Lambda^{-1}\Phi_\Lambda u_\Lambda$  are supersmooth in the sense of Definition \ref{DefSuperSmoothDom} for each $\Lambda \in \Gr$ w.r.t. arbitrary charts $u$ on $X$ and $v$ on $Y$. The set of these morphisms will be denoted by $\SC^\infty(X,Y)$. 
\end{Rem}

Using the notion of supermanifold Definition \ref{DefSMfFunc} and morphisms as discussed in Remark \ref{RemSMfMorphism}, we obtain the category $\SMan$ of supermanifolds. It contains the full subcategory $\mathrm{fSMan}$ of finite-dimensional objects. $\mathrm{fSMan}$ is a priori different from the category $\BKL$. It can be shown (see \cite{Sa2} Theorem 5.1, see also \cite{AL} section 3.3 for a different formulation) that the assignments
\begin{align} \label{EqAequiCat}
X \in \mathrm{fSMan} &\mapsto (X(\R), \SC^\infty(X,\R^{1|1})) \in \BKL \\
X \in \BKL &\mapsto \Hom_\BKL(-,X) \in \mathrm{fSMan} \subset \Man^\Gr \notag
\end{align}
define an equivalence of categories. Thus, both definitions coincide in the finite-dimensional case.\\

In the Molotkov-Sachse framework, there is a natural candidate for the space of morphisms, given by the common adjunction formula for inner $\Hom$ objects (see \cite{GM1}, II.4.23 and also the discussion of ``exponential law'' in \cite{KM}, e.g. I.3, p.76, 23.2, 42.14):
\begin{Def}[Sachse, \cite{Sa1} 7.1.2]\label{DefSCInfty}
For finite-dimensional supermanifolds $X,Y$, the functor $\qq{\SC}^\infty(X,Y): \Gr \rightarrow \Set$ is defined by
\begin{align*}
\Lambda                               &\mapsto \qq{\SC}^\infty(X,Y)(\Lambda) := \SC^\infty(\PR(\Lambda) \times X, Y) \\
(\rho : \Lambda \rightarrow \Lambda') &\mapsto \qq{\SC}^\infty(\rho) := (\sigma \mapsto \sigma \circ (\PR(\rho) \times \mathrm{Id}_X))
\end{align*}
\end{Def}
As it is defined, it takes values in the category of sets. We will turn it into a functor $\Gr \rightarrow \Man$ in the course of this work and eventually equip it with a supersmooth structure. However, it should be pointed out that by part (b) of the following remark, the supermanifold is not an inner $\Hom$ object in the category $\SMan$.

\begin{Rem} \label{RemAdvMSApproach}
The Molotkov-Sachse-approach uses functors $\mathrm{SPoint}^{\op} \rightarrow \Man$ to describe supermanifolds. As mentioned above, some authors work with the larger category $\BKL$ instead of $\mathrm{SPoint}$ (see e.g. \cite{AMS1} or \cite{Var}). In view of the construction of mapping spaces based on Definition \ref{DefSCInfty}, choosing $\BKL$ has two disadvantages: 
\begin{enumerate}
 \item In case $M$ is compact manifold, it is well known that $C^\infty(M,N)$ is a Fr\'{e}chet manifold (\cite{KM} 42.3). Hence, we expect $\qq{\SC}^\infty(X,Y)$ to take values in the category of Fr\'{e}chet-
       manifolds in case $\tilde{X}$ is compact. However, if we allow for arbitrary test objects $S \in \BKL$, the manifolds $\widetilde{S \times X} \cong \tilde{S} \times \tilde{X}$ will in general not be compact and consequently, $\qq{\SC}^\infty(S \times X, Y)$ will not be Fr\'{e}chet. For $S \in \mathrm{SPoint}$, we have $\tilde{S} = \{\mathrm{point}\}$ and the problem does not occur. 
 \item Let $M,N,S$ be smooth manifolds. Using the manifold structure on $C^\infty(M,N)$ sketched in chapter \ref{ChapPre}, it is known that in 
       general (for $M$ not compact) $C^\infty(S,C^\infty(M,N)) \neq C^\infty(S \times M,N)$ (see \cite{KM}, 42.14 for details of this failure of cartesian closedness). However, we need equality for Definition \ref{DefSCInfty} to make sense. If $S$ is just a point (i.e. the underlying manifold of a superpoint) the problem again disappears. Moreover, this issue also indicates that the supermanifold $\qq{SC}^\infty(X,Y)$, constructed in Theorem \ref{ThmSMfSCInfty}, will in general not be an inner $\Hom$ object in $\SMan$. In fact, already for purely even supermanifolds $S,X,Y$ (i.e. ordinary smooth manifolds), the required identification $C^\infty(S,C^\infty(X,Y)) \cong C^\infty(S \times X,Y))$ does not exist as mentioned above.
\end{enumerate}
Of course, one may also resort to a weaker notion of smoothness on spaces of the form $C^\infty(M,N)$ to circumvent these difficulties (cf. \cite{KM}, p.429/430). Remarkably, the approach of Alldridge  does not meet problem b) either (\cite{All}, Theorem 2.28). It would be interesting to study the underlying reason for this different behaviour.
\end{Rem}

\vspace{\baselineskip}

We finally introduce a slight extension of Definition \ref{DefSuperSmoothDom}. It will be relevant in chapter \ref{ChapSmooth} for the construction of mapping spaces with \emph{non-compact domain}, since these do in general not fit into the framework discussed so far. The following example (a simplification of \eqref{EqIdentSpace}) indicates the problem: 

\begin{Ex}\label{ExCounterChart}
Let $M,N$ be manifolds and assume that $M$ is not compact. Consider the vector bundle\footnote{Here, $\Lambda^{ev\geq 2}$ also denotes the trivial vector bundle over $N$ with fibre $\Lambda^{ev\geq 2}$ so that $\Lambda^{ev\geq 2} \otimes TN$ is just the tensor product of bundles.} $T_\Lambda := \Lambda^{ev\geq 2} \otimes TN \rightarrow N$. Setting $T(\Lambda):= C^\infty(M,T_\Lambda)$, we obtain a functor $\Gr \rightarrow \Set$ in the obvious way. We will see in chapter \ref{ChapSmooth} (Theorem \ref{ThmBundleStruct} and the preceding discussion) that $C^\infty(M,T_\Lambda)$ is a smooth manifold and its model space near $\sigma$ is given by 
\begin{align*}
\Gamma_c(M,f^\ast TN) \oplus  \Gamma_c(M,\Lambda^{\ev\geq 2} \otimes f^\ast TN) &\cong 
(\Lambda \otimes \Gamma_c(M,f^\ast TN))_\0,
\end{align*}
where $f \in C^\infty(M,N)$ is the map underlying $\sigma$. It is the set of $\Lambda$-points of the superrepresentable module $\qq{\Gamma_c(M,f^\ast TN)}$ where $\Gamma_c(M,f^\ast TN)$ is purely even. By \eqref{EqRestriction}, its open subfunctors are restrictions, i.e. of the form $\Lambda \mapsto V \times \Gamma_c(M,\Lambda^{\ev\geq 2} \otimes f^\ast TN) = \qq{\Gamma_c(f^\ast TN)}|_{V}(\Lambda)$ for some $c^\infty$-open set $V \subset C^\infty(M,f^\ast TN)$.\\
Functors of the form $T$ will arise in the description of spaces of morphisms. Assume it admits supercharts, i.e. there exists an isomorphism of functors $\alpha : \qq{\Gamma_c(f^\ast TN)}|_{V} \overset{\sim}{\rightarrow} U$ onto an open subfunctor $U \subset T$. $\qq{\Gamma_c(f^\ast TN)}|_{V}(\Lambda)$ contains elements of the form $(v,0)$ for $v \in V$ and $0 \in \Gamma_c(M,f^\ast TN \otimes \Lambda^{\ev\geq 2})$. Such elements are clearly invariant under the action of all $\rho \in \Hom_\Gr(\Lambda,\Lambda')$ and due to functoriality, so are $\alpha_{\Lambda}(v,0)$. Putting $\rho = \pr^\Lambda_\R : \Lambda \twoheadrightarrow \R$, it follows that $\alpha_{\Lambda}(v,0) \in C^\infty(M,TN\otimes\R)$. But due to the choice of the topology on $C^\infty(M,T_\Lambda)$ in \eqref{EqCharts2}, this implies that $U(\Lambda) \subset C^\infty_c(M,T_\Lambda)$. In particular, we see that it is impossible to cover $C^\infty(M,T_\Lambda)$ and obtain an atlas on $T$ with charts of the form \eqref{EqRestriction}, using the constructions in \eqref{EqCharts1}, \eqref{EqCharts2}. This argument moreover shows that the naive idea, to choose reference sections $\sigma_{0,\Lambda} \in C^\infty(M,T_\Lambda)$ and consider open chart neighbourhoods of compactly supported perturbations of $\sigma_{0,\Lambda}$ as in \eqref{EqCharts2} necessarily breaks functoriality. Intuitively, one needs open neighbourhoods in $C^\infty(M,T_\Lambda)$ without restrictions on the supports of sections in $T_\Lambda$.
\end{Ex}

To circumvent these difficulties, one may either change the topologies on mapping spaces or allow for more general superdomains. Since I am not aware of appropriate topologies on spaces like $C^\infty(M,T_\Lambda)$, we will generalise the notion of superdomain. Recall from Example \ref{ExLocAffine} that $\Gamma(f^\ast T_\Lambda)$ carries the structure of a locally affine space (and hence the structure of a smooth manifold) modelled on $\Gamma_c(f^\ast T_\Lambda)$. If $V \subset C^\infty(M,f^\ast TN)$ is open, then the set $\qq{\Gamma(f^\ast TN)}|_{V}(\Lambda)$ from Example \ref{ExCounterChart} is a disconnected, open subset of $\Gamma(f^\ast T_\Lambda)$ (where all spaces may be equipped with $c^\infty$-topologies). Its (trivial) tangent bundle is given by 
\begin{align}  \label{EqTanlocAff}
T(\qq{\Gamma(f^\ast TN)}|_{V}(\Lambda)) &= \qq{\Gamma(f^\ast TN)}|_{V}(\Lambda) \oplus \qq{\Gamma_c(f^\ast TN)}(\Lambda) 
\end{align}

Thus, the following definition is reasonable:

\begin{Def} \label{DefLocAffDom}
\begin{enumerate}
 \item Let $\qq{A}$, $\qq{W}$ be superrepresentable modules. $\qq{A}$ is called locally affine with model $\qq{W}$, if for each $\Lambda \in \Gr$, $\qq{A}(\Lambda)$ is a locally affine space modelled on $\qq{W}(\Lambda)$ in the sense of Example \ref{ExLocAffine}. Open subfunctors of $\qq{A}$ are called (convenient) locally affine superdomains. 
 \item A morphism $D_1 \rightarrow D_2$ in $\TGr$ between locally affine superdomains is called supersmooth, if it satisfies the two conditions from Definition \ref{DefSuperSmoothDom}.
\end{enumerate}
\end{Def}

Note that the locally affine superdomains are still precisely the restrictions of locally affine superrepresentable modules as in \eqref{EqRestriction}. In fact, the proof of Proposition 4.8 in \cite{Sa2} is still valid in this situation. Using this enlarged class of superdomains, we may finally set

\begin{Def} \label{DefLocAffSMf}
A locally affine supermanifold is a functor $X \in \Man^\Gr$ satisfying the conditions from Definition \ref{DefSMfFunc} with ``superdomain'' replaced by ``locally affine superdomain''.
\end{Def}

\begin{Rem}
Even though locally affine supermanifolds are modelled on locally affine spaces, note that this concepts still allows for reasonable tangent objects. It follows from \eqref{EqTanlocAff} and the $\Lambda^{\ev}$-linearity of chart transitions demanded in Definitions \ref{DefSuperSmoothDom} and \ref{DefSMfFunc}, respectively, that tangent spaces are still (convenient) superrepresentable modules and not just locally affine. In fact, the construction of tangent bundles in \cite{Sa1} chapter 3.8.2  carries over without difficulties. Loosely speaking, we only changed the class of allowed coordinate neighbourhoods in order to be able to use functor language. These neighbourhoods now may be disjoint unions of sets which, individually, are diffeomorphic to $c^\infty$-open subsets of a convenient vector spaces but whose unions only fit into the bigger, locally affine space.
\end{Rem}

 
\section{The algebraic Description of $\qq{\SC}^\infty(X,Y)$} \label{ChapAlg}

In this chapter, we will determine the geometric and algebraic structure of higher points of $\qq{\SC}^\infty(X,Y)$ where $X = (\tilde{X},\O)$ and $Y = (\tilde{Y},\cR)$ are supermanifolds of dimension $p|q$ and $p'|q'$ respectively, defined using either the BKL- or the equivalent Molotkov-Sachse picture. This has been done in \cite{Sa1}, chapter 7.2.1 and \cite{SaWo}, chapter 5.2 for invertible elements of $\qq{\SC}^\infty(X,X)$. A similar technique is also described in \cite{Hel1}. For simplicity, we will restrict to a skeleton of $\Gr$ by choosing the representative $\Lambda_n := \tL \R^n$ for each isomorphism class of finite-dimensional Grassmann algebras. Let $\PR_n := (\{0\}, \Lambda_n)$ denote the corresponding superpoint. We will expand the elements of $\qq{\SC}^\infty(X,Y)(\Lambda_n) \cong \Hom_\BKL(\PR_n \times X, Y)$ w.r.t. odd directions in $\PR_n$, i.e. w.r.t. a basis of $\Lambda_n$. The resulting coefficients turn out to be super differential operators of suitable degree and parity. We will first discuss the theory of linear differential operators and jet modules on supermanifolds using an algebraic approach. It starts with the definition of differential operators, rather than jet bundles, as the fundamental notion. General references for this approach are \cite{Kras1} chapter 1, \cite{Nes1} chapter 9 and \cite{Sad1} chapter 6.2 for superjet modules. A slightly different approach is discussed in \cite{CCF}, chapter 4.4\\

Let $A$ be a supercommutative $\R$-superalgebra\footnote{It is possible to use a general supercommutative ring instead of $\R$ (see the references quoted above) but for the present work, $\R$ is sufficient}. Let $P,Q$ be supermodules over $A$ or some other superring which will be indicated if necessary. When speaking about linear maps, we will indicate the underlying ring of scalars by a subscript. Later, $A$ will be the algebra of (local) superfunctions on some supermanifold and $P,Q$ the locally free modules of sections in super vector bundles.\\ 
Following the exposition in \cite{Kras1}, chapter 1 \textsection 1 and \cite{Sad2} section 3.2, we introduce a left and a right $A$-module structure on $\qq{\Hom}_\R(P,Q)$ as well as their commutator $\delta$ by
\begin{align} \label{EqDefModulHom}
(a\cdot \Delta)(p) &:= a\Delta(p) &
(\Delta \cdot a)(p) &:= \Delta(ap) & 
\delta_a \Delta := a\cdot \Delta - (-1)^{|a||\Delta|} \Delta\cdot a
\end{align}
Here $\Delta \in \qq{\Hom}_\R(P,Q)$, $a \in A$ and $p \in P$. Abstract differential operators (we will omit the prefix ``super'') are now defined by

\begin{Def} \label{DefDO}
A map $\Delta \in \qq{\Hom}_\R(P,Q)$ is called linear differential operator from $P$ to $Q$ of order $ \leq \ k \in \N_0$ and parity $|\Delta| \in \Z_2$ if the following identity is satisfied:
\begin{align*}
 \forall a_0,\ldots a_k \in A : \quad \delta_{a_0} \circ \cdots \circ \delta_{a_k} \Delta = 0
\end{align*}
The set of all differential operators of order $\leq k$ from $P$ to $Q$ is denoted by $\Diff^k(P,Q)$   
\end{Def}

It is clear that $\Diff^k(P,Q)$ inherits from $\qq{\Hom}_\R(P,Q)$ a left, a right and an $A$-bi-supermodule structure. Moreover, differential operators of order zero are just $\R$-linear maps, which supercommute with $A$, i.e. which are $A$-linear. Thus, we have $\Diff^0(P,Q) = \qq{\Hom}_A(P,Q)$. Differential operators show the expected behaviour under composition: If $R$ is another $A$-supermodule, then
\begin{align} \label{EqComposeDiffs}
\Delta \in \Diff^k(P,Q), \ \Delta' \in \Diff^k(Q,R) &\Longrightarrow \Delta' \circ \Delta \in \Diff^{k+l}(P,R)
\end{align}
The parity of $\Delta' \circ \Delta$ is given by $|\Delta'|+|\Delta|$ (cf. \cite{Kras1}, Prop 1.1.3, the rule for the parities is obtained as usual).
\begin{Ex} \label{Ex1DimDOs}
Let $P = Q = A$ and $\Delta \in \Diff^1(A,A)$. We decompose $\Delta = (\Delta - \Delta(1)) + \Delta(1)$ where $\Delta(1) \in A = \Diff^0(A,A)$ is considered to be a multiplication operator (cf. \cite{Sad2}, p.57/58). A simple calculation shows that $\Delta-\Delta(1)$ is a super derivation of parity $|\Delta|$ on $A$. For $A =\O(U)$, this says that first order operators can be decomposed into super vector fields and multiplication operators. In fact, this decomposition is a direct sum (\cite{Kras1} (1.1)). For higher order operators, the situation is more involved (see e.g. \cite{Kras1} Def. 1.1.9 ff).\\
Vector fields along morphisms as defined in chapter \ref{ChapPre} (above \eqref{EqPullBack}) fit into this framework as follows: Let $\Phi=(\Phi^\ast,\varphi) : X \rightarrow Y$ be a morphism of supermanifolds. For $V \subset \tilde{Y}$ open, let $P = A = \cR(V)$ and $Q := \O(\varphi^{-1}(V))_\Phi$ the $\cR(V)$-module defined above \eqref{EqExtScalars}. An element $\Delta \in \Diff^1(\cR(V),Q)$ without constant term satisfies
\begin{align*}
   \Delta(gg') = \Delta(g)\cdot g' + (-1)^{|g||\Delta|} g\cdot \Delta(g')
               = \Delta(g)\Phi(g') + (-1)^{|g||\Delta|} \Phi(g) \Delta(g')
\end{align*}
Hence, $\Delta$ satisfies the super Leibniz rule along $\Phi$ and elements of $\Diff^1(\cR(V),\O_\Phi(V))$ without constant term generalise the classical notion of vector fields along maps. This example will motivate Definition \ref{DefDoAlong}.
\end{Ex}

\vspace{\baselineskip}

Let $U \subset \tilde{X}$ be open and $A:= \O(U)$. The elements of $\Diff^k(\O(U),\O(U))$ are called linear differential operators of order $\leq k$, acting on superfunctions over $U$. Let $U$ be a coordinate neighbourhood with coordinate functions $(x^i,\theta^\alpha)$. For each multi-index $I = (I_\0,I_\1)$ s.t. $I_\0 = (i_1,\ldots,i_p)\in \N_0^p$ and $I_\1 = (\alpha_1,\ldots, \alpha_q) \in \Z_2^q$, we denote
\begin{align} \label{EqDefCoordDos}
 \frac{\del^{|I|}}{\del \xi^I} := 
 \frac{\del^{i_1}}{\del (x^1)^{i_1}} \circ\cdots\circ \frac{\del^{i_p}}{\del (x^p)^{i_p}} \circ\cdots\circ \frac{\del^{\alpha_1}}{\del (\theta^1)^{\alpha_1}} \circ\cdots\circ \frac{\del^{\alpha_q}}{\del (\theta^q)^{\alpha_q}}
\end{align}
By the Leibniz rule, $\del^{|I|}/\del\xi^I \in \Diff^{|I|}(\O(U),\O(U))$ and its parity is given by $|I_\1| \mod 2$. The local and global properties of differential operators on $X$ are given in the next theorem (cf. \cite{CCF} Proposition 4.4.11 or \cite{Kos} section 2.9 (see also \cite{FH} 5.2 for details of the proof which are omitted in the other references):

\begin{Thm} \label{ThmSMfDOs}
For any supermanifold $X = (\tilde{X},\O)$ of dimension $p|q$, the assignment
\begin{align*}
U &\mapsto \Diff^k(\O,\O)(U) := \Diff^k(\O(U),\O(U)) \subset \Hom_\R(\O(U),\O(U)) 
\end{align*}
defines a locally free sheaf of $\O$-modules over $X$. In coordinates, a local base is given by the coordinate differential operators $\del^{|I|}/\del\xi^I$ from \eqref{EqDefCoordDos} for $|I| = |I_\0|+|I_\1| \leq k$.  
\end{Thm}
  
\begin{Rem} \label{RemSMfDOs}
\begin{enumerate}
\item Theorem \ref{ThmSMfDOs} generalises to sections of super vector bundles $\E,\F$. This follows from the fact that locally, $\E \cong \O^{\mathrm{rk}(E)}$, $\F \cong \O^{\mathrm{rk}(F)}$ so that 
      $\Diff^k(\E,\F)$ inherits the properties from $\Diff^k(\O,\O)$.  
\item Theorem \ref{ThmSMfDOs} implicitly contains the statement, that operators $\Delta \in \Diff^k(\E,\F)(U)$ are local,  
       i.e. $\supp(\Delta(e)) \subset \supp(e)$ for all $e \in \E(U)$. 
 \end{enumerate}
\end{Rem}

The following definition generalises Example \ref{Ex1DimDOs} to higher order operators:

\begin{Def} \label{DefDoAlong}
Let $\Phi=(\Phi^\ast,\varphi) : X \rightarrow Y $ be a morphism of supermanifolds, $P$ a $\cR$-module on $\tilde{Y}$ and $Q$ an $\O$-module on $\tilde{X}$. A linear differential operator from $P$ to $Q$ along $\Phi$ of degree $\leq k$ over the open set $V \subset \tilde{Y}$ is an element of $\Diff^k(P,\varphi_\ast Q_\Phi)(V) = \Diff^k(P(V), \varphi_\ast Q(V)_\Phi)$. Differential operators along $\Phi$ over open sets $U \subset \tilde{X}$ are defined by $\varphi^\ast \Diff^k(P,\varphi_\ast Q_\Phi)(U)$. 
\end{Def}

Recall that the $\cR(V)$-module structure on $\varphi_\ast Q(V)_\Phi$ was defined by $r\cdot q = \Phi^\ast(r)q$ above \eqref{EqExtScalars}. Clearly, $\Phi^\ast$ is a differential operator of order zero along $\Phi$. Differential operators along morphisms have the following property:

\begin{Prop} \label{PropDOAlong}
Let $\Phi$, $P$ and $Q$ satisfy the assumptions of the preceding definition and assume that $P$ is locally free. Then we have for all $V \subset \tilde{Y}$:
\begin{align} \label{EqPropDOAlong}
 \Diff^k(P(V),\varphi_\ast Q(V)_\Phi) &\cong \varphi_\ast Q(V) \otimes_\Phi \Diff^k(P(V),\cR(V)) 
\end{align}
In particular, differential operators along morphisms also form a locally free sheaf over $\tilde{Y}$ and $\tilde{X}$, respectively.
\end{Prop}

\begin{Proof}{}
By Theorem \ref{ThmSMfDOs} and Remark \ref{RemSMfDOs}, the l.h.s. of \eqref{EqPropDOAlong} is a sheaf and moreover, by Proposition 4.6 in \cite{SR1}, the tensor product on the r.h.s. defines a sheaf, too. Hence, both sides define locally free sheaves and it is sufficient to prove the statement for open sets $V \subset \tilde{Y}$ such that $P(V)$ is a free module. Taking the pullback along $\varphi$, we get the corresponding statement for the operators defined over $U \subset \tilde{X}$.\\
Observe that there is an isomorphism of left $\cR(V)$-supermodules $\cR(V) \otimes_\Phi \varphi_\ast Q(V) \cong \varphi_\ast Q(V)_\Phi$ given by $r \otimes q \mapsto \Phi(r)q$ whose inverse is $q \mapsto 1 \otimes q$. Proposition \ref{PropJet} (b) yields
\begin{align*}
 \Diff^k(P(V),\varphi_\ast Q(V)_\Phi) 
&\cong \Diff^k(P(V), \cR(V) \otimes_\Phi \varphi_\ast Q(V)) \\
&\cong \qq{\Hom}_{\cR(V)}(\Jet^k P(V), \cR(V) \otimes_\Phi \varphi_\ast Q(V)) \\
&\cong \qq{\Hom}_{\cR(V)}(\Jet^k P(V), \cR(V)) \otimes_\Phi \varphi^\ast Q(V) \\
&\cong \varphi_\ast Q(V) \otimes_\Phi \Diff^k(P(V), \cR(V)) 
\end{align*}
In the third step, we used the compatibility of $\qq{\Hom}$ and $\otimes$ (see \cite{BouAlg1} II.4.2.2 (ii) which generalises to supermodules) taking into account that $P(V)$ is a free module of finite rank. 
\end{Proof}

\begin{Rem}
The proposition states that a differential operator along $\Phi$ is essentially given by an $\O$-linear combination of differential operators on $Y$ of the same order composed with the morphism 
       $\Phi$: In the case $P = \cR, Q = \O$,  we may identify $\Delta \in \Diff^k(\cR(V), \O(\varphi^{-1}(V))_\Phi)$ with $\sum_i f_i \otimes_\Phi \Delta_i$ for suitable $\Delta_i \in \Diff^k(\cR((V)),\cR(V))$ and $f_i \in \varphi_\ast\O(V)$. For $g \in \cR(V)$, the action of $\Delta$ on $g$ is given by (cf. \eqref{EqPullBack})
       \begin{align} \label{EqCoordPBDiff}
          \Delta(g) &= \tsum_i f_i \otimes_\Phi \Delta_i(g)  = f_i\Phi(\Delta_i(g)) 
       \end{align}
By Proposition \ref{PropDOAlong}, the sheaf of differential operators of order $\leq 1$ without constant term over $\tilde{X}$ acting on superfunctions is given by $\O \otimes_\Phi \varphi^\ast \T_Y = 
       \Phi^\ast \T_Y$. Hence, we recover Example \ref{Ex1DimDOs}.  
\end{Rem}

To obtain an analogue of \eqref{EqDOJetDiag} and \eqref{EqJetFactorProp}, consider for $b \in A$ the map $\delta^b \in \qq{\End}_\R(A \otimes_\R P)$ defined by $\delta^b(a \otimes p) = (ba)\otimes p - (-1)^{|b||a|} a \otimes (bp)$. Following \cite{Sad2} 3.2.6, the $k$-th superjet module and the natural jet map are defined by
\begin{align}\label{EqDefJetMod} 
 \Jet^k(P) &:= A \otimes_\R P / M_{k+1}  & &\text{for } M_{k+1}:= \mathrm{span}_A\{\delta^{a_0}\cdots\delta^{a_k}(x) | x \in A\otimes_\R P, a_i \in A \} \\
\jet^k : P  &\rightarrow \Jet^k(P), & &\mspace{-40mu}\jet^k(p) := 1 \otimes p + M_{k+1} \notag
\end{align}
$\Jet^k(P)$ inherits an $A$-module structure from $A \otimes_\R P$. The following properties can be found in \cite{Kras1} I.2 :

\begin{Prop}\label{PropJet}
 Let $P,Q$ be $A$-modules and $k \in \N_0$, then we have:
\begin{enumerate}
\item $\jet^k \in \Diff^k(P,\Jet^k(P))_\0$. Its image generates $\Jet^k(P)$.
\item Each $\Delta \in \Diff^k(P,Q)$ uniquely factors over $\jet^k$, i.e. we have natural bijections 
       \begin{align*}
        \Diff^k(P,Q) &\cong \qq{\Hom}(\Jet^k(P),Q) & \Delta = \mathrm{sb}(\Delta) \circ \jet^k &\mapsto \mathrm{sb}(\Delta)
       \end{align*}
\item For $k \leq l$, there is a unique map $\pi^{l,k} \in \Hom_A(\Jet^l(P),\Jet^k(P))$ such that $\jet^k = \pi^{k,l} \circ \jet^l$. For $k \leq l     
       \leq m$, we furthermore have $\pi^{m,l} \circ \pi^{l,k} = \pi^{m,k}$.
\end{enumerate}
\end{Prop}

\begin{Rem} \label{RemJetProblem}
In contrast to the algebraic construction of differential operators, which is equivalent to the common notion in local coordinates by Theorem \ref{ThmSMfDOs}, this is not true for the jet modules. In fact, to obtain the jet bundle in the sense of a locally free sheaf of supermodules from the module $\Jet^k(P)$ constructed above, a ``geometrisation procedure'' has to be applied. We will not discuss these details but refer to \textsection 1.3 of \cite{Kras1}. See also \cite{HRMM} for a different approach to superjet bundles.    
\end{Rem}

In case $P=A$, $A \otimes_\R A$ is an $A$-superalgebra by $(a \otimes b) \cdot (a' \otimes b' ) = (-1)^{|b||a'|} (aa')\otimes(bb')$. We have the following elementary lemma:

\begin{Lemma} \label{LemJetAlg}
For all $k \in \N_0$, the submodule $M_k \subset A \otimes_\R A$ from \eqref{EqDefJetMod} is a graded ideal. In particular, $\Jet^{k-1}(A) = A / M_k$ carries a natural structure of a $\R$-superalgebra and we have $\jet^k(ab) = \jet^k(a) \jet^k(b)$.
\end{Lemma}

Based on the formalism introduced so far, we now determine the structure of the higher points of $\qq{\SC}^\infty(X,Y)$. These results have been obtained in \cite{FH} with different proofs. Let $\Phi^{(n)} \in \qq{\SC}^\infty(X,Y)(\Lambda_n)$. By Definition \ref{DefSCInfty}, \eqref{EqAequiCat} and Proposition \ref{PropMorGlob}, we may identify
\begin{align} \label{IdentPhiN}
\Phi^{(n)} \in \SC^\infty(\PR_n \times X, Y) &\cong \Hom_\BKL(\PR_n \times X, Y) \cong \Hom_{\mathrm{SAlg}}(\cR(\tilde{Y}), \Lambda_n \otimes \O(\tilde{X}))  
\end{align}
From this point, we will focus on the algebras $\O := \O(\tilde{X})$, $\cR := \cR(\tilde{Y})$ of global sections, keeping in mind that the full sheaf morphisms can be reconstructed from the global data (cf. Proposition \ref{PropMorGlob}). Recall that $\Phi^{(n)\ast}$ is even by definition, i.e. preserves parity. Denoting the generators of $\Lambda_n$ by $\eta_1,\ldots,\eta_n$ and choosing $g \in \cR$, we can decompose $\Phi^{(n)\ast}$ as
\begin{align} \label{EqFineStr1}
 \Phi^{(n)\ast}(g) &= \sum_{I \in \Z_2^n} \eta^I \Phi^{(n)}_I(g) & \text{ with } \quad \Phi^{(n)}_I: \cR &\map \O  \quad \text{ for } I \in \Z_2^n
 \end{align}
We regard $\Phi^{(n)}_I$ as maps between global superfunctions but again, we may localize to obtain morphisms of sheaves of vector spaces. Since $\Phi^{(n)\ast}$ is even, we obviously have 
\begin{align*}
|\Phi^{(n)}_I| &= |I|\mod 2  & |\Phi^{(n)}_I(g)| &= (|I| \mod 2) + |g|
\end{align*}
  
By \eqref{IdentPhiN}, $\Phi^{(n)\ast}$ is a unital homomorphism of $\R$-superalgebras and this property determines the structure of the coefficients $\Phi^{(n)}_I$. We immediately obtain

\begin{Lemma} \label{LemSimpPropCoef}
 \begin{enumerate}
  \item $\Phi^{(n)}_\emptyset$ is a morphism of superalgebras $\cR \rightarrow \O$. The induced morphism of supermanifolds will also be denoted by $\Phi^{(n)}_\emptyset : X \rightarrow Y$.
  \item For $I \neq \emptyset$, we have $\Phi^{(n)}_I(1) = 0$. 
 \end{enumerate}
\end{Lemma}
  
Hence, the structure of $\Phi^{(n)}_\emptyset$ is completely understood. We will now show that the remaining coefficients are differential operators:

\begin{Thm} \label{ThmStructComp}
 Let $\Phi^{(n)} \in \qq{\SC}^\infty(X,Y)(\Lambda_n)$ and let $\Phi^{(n)}_I$ ($|I| > 0$) be the coefficients from \eqref{EqFineStr1}. Then $\Phi^{(n)}_{I} \in \qq{\Hom}_\R(\cR, \O$\raisebox{-1.5pt}{$\scriptscriptstyle \Phi^{(n)}_\emptyset$}$)$ is a linear differential operator along $\Phi^{n}_\emptyset$ of degree $\leq |I|$ and of parity $(|I| \mod 2)$ . In particular $\Phi^{(n)}_I$ are sections of a super vector bundle on $X$.    
\end{Thm}

\begin{Proof}{}
For $k \in \N_0$, let $\iota_{k+1}:\PR_k \hookrightarrow \PR_{k+1}$ be the morphism defined by $\iota_{k+1}^\ast(\eta^{k+1}) =0$. We rewrite \eqref{EqFineStr1} as
\begin{align} \label{EqRewriteMorph}
 \Phi^{(k+1)\ast} &= \sum_{I \in \Z_2^k} \eta^I \Phi^{(k+1)}_I + \sum_{I \in \Z_2^k} \eta^I\eta^{k+1} \Phi^{(k+1)}_{(I,k+1)} =: \iota_{k+1}^\ast \Phi^{(k+1)\ast} + D^{(k+1)}
\end{align}
For $f,g \in \cR$, the multiplicativity of $\Phi^{(k+1)\ast}$ implies 
\begin{align}
 \iota_{k+1}^\ast\Phi^{(k+1)\ast}(fg) &= \iota_{k+1}^\ast\Phi^{(k+1)\ast}(f)\iota_{k+1}^\ast\Phi^{(k+1)\ast}(g) \label{EqSplitKp1} \\
 D^{(k+1)}(fg) &= \iota_{k+1}^\ast\Phi^{(k+1)\ast}(f) D^{(k+1)}(g) + D^{(k+1)}(f) \iota_{k+1}^\ast\Phi^{(k+1)\ast}(g) \notag
\end{align}
$\iota_{k+1}^\ast\Phi^{(k+1)\ast} : \cR \rightarrow \Lambda_n \otimes \O \subset \Lambda_{k+1} \otimes \O$ defines the morphism\footnote{Of course, this also follows from the functoriality in Definition \ref{DefSCInfty}} $\Phi^{(k+1)}\iota_{k+1} : \PR_k \times X \rightarrow Y$ which we trivially extend to $\PR_{k+1} \times X$. Moreover, $D^{(k+1)}$ is an even derivation along $\Phi^{(k+1)}\iota_{k+1}$. We now prove the theorem by induction on the number of generators $n$ of the Grassmann algebra $\Lambda_n$. 

The statement in the case $n=1$ follows from \eqref{EqSplitKp1} (set $k=0$): Since $D^{(1)} = \eta^1 \Phi^{(1)}_1$ is an even derivation along $\Phi^{(1)}_\emptyset = \Phi^{(1)}\iota_1$, $\Phi^{(1)}_1$ is an odd differential operator of order $1$ along $\Phi^{(1)}_\emptyset$. To prove the step $n \rightarrow n+1$, we use \eqref{EqSplitKp1} with $k=n$. Since $D^{(n+1)}$ is a first order operator along $\Phi^{(n+1)}\iota_{n+1}$ without constant term, Proposition \ref{PropDOAlong} yields
\begin{align*}
 D^{(n+1)} = \sum_l h^l \otimes_{\Phi^{(n+1)}\iota_{n+1}} \Delta_l 
\end{align*}
for some coefficients $h^l = \sum_{K \in \Z_2^{n+1}} \eta^K h^l_K \in \Lambda_{n+1} \otimes \O$ and operators $\Delta_l \in \Diff^1(\cR,\cR)$. Since $D^{(n+1)}$ is even, we moreover have $|h^l| = |\Delta_l|$. Applying the induction hypothesis to $\Phi^{(n+1)}\iota_{n+1} : \PR_n \times X \rightarrow Y$, we know that its components $\Phi^{(n+1)}_I$ ($I \in \Z_2^{n}$) are differential operators of degree $\leq |I|$ and parity $|I| \mod 2$. As in \eqref{EqCoordPBDiff}, we have for $f \in \cR$
\begin{align*}
 \sum_{I \in \Z_2^{n}} \eta^{I}\eta^{n+1} \Phi^{(n+1)}_{(I,n+1)}(f) & = D^{(n+1)}(f) = \sum_{l,J \in \Z_2^{n}, K \in \Z_2^{n+1}} \eta^K h^l_K \eta^J (\Phi^{(n+1)}_J \circ \Delta_l)(f)  
\end{align*}
Comparing the $\eta$-decompositions of the left- and right hand side of this equation yields
\begin{align} \label{EqDecompD}
\Phi^{(n+1)}_{(I,n+1)} &= \sum_{l, K \sqcup J = (I,n+1)} (-1)^{\sign(K,J)} h^l_K \Phi^{(n+1)}_J \circ \Delta_l 
\end{align}
Since $\Phi^{(n+1)}_J$ and $\Delta_l$ have order $\leq |J|$ and $1$, respectively, it follows that $\Phi^{(n+1)}_{(I,n+1)}$ has order $\leq |I| + 1$ and it remains to prove that $\Phi^{(n+1)}_{(I,n+1)}$ has parity $(|I|+1) \mod 2$. Since we have $\Phi^{(n+1)}_J = |J| \mod 2$ by assumption and $|\Delta_l| = |h^l| = |h^l_K| + (|K| \mod 2)$ for each $K \in \Z_2^{n+1}$, we conclude from \eqref{EqDecompD} and the parity rule below \eqref{EqComposeDiffs} that
\begin{align*}
 |\Phi^{(n+1)}_{(I,n+1)}| &= |h^l_K| + |\Phi^{n+1}_{J}| + |\Delta_l| = (|K| + |J|) \mod 2.
\end{align*}
The sum in \eqref{EqDecompD} runs over all multi-indices $K,J$ satisfying $K \sqcup J = (I,n+1)$, hence $|\Phi^{(n+1)}_{(I,n+1)}| = (|K| + |J|) \mod 2 = |(I,n+1)| \mod 2$, which finishes the proof.
\end{Proof}

It is easy to see from this theorem, that $\Phi^{(n)\ast}$ is a differential operator of order $\leq n$. In fact, from \eqref{EqFineStr1} and Proposition \ref{PropJet}, we obtain
\begin{align*}
\Phi^{(n)\ast} = \sum_{I\in \Z_2^{n}} \eta^I \mathrm{sb}(\Phi^{(n)}_I) \circ \jet^{|I|} = (\sum_{I\in \Z_2^{n}} \eta^I \mathrm{sb}(\Phi^{(n)}_I) \circ \pi^{|I|,n}) \circ \jet^n \in \Hom_{\mathrm{SAlg}}(\cR,\Lambda_n \otimes \O)
\end{align*}
where $\mathrm{sb}(\Phi^{(n)}) =\sum_{I\in \Z_2^n} \eta^I \mathrm{sb}(\Phi^{(n)}_I) \circ \pi^{|I|,n} \in \Hom_{\mathrm{SAlg}}(\Jet^n(\cR),\Lambda_n \otimes \O)$ which proves the claim. Since $\jet^n$ is an algebra homomorphism  whose image generates $\Jet^n(\cR)$ by Lemma \ref{LemJetAlg}, we obtain: 

\begin{Cor} \label{CorStructComp2}
Every $\Phi^{(n)} \in \qq{\SC}^\infty(X,Y)(\Lambda_n)$ defines a differential operator of order $\leq n$ along $\Phi^{(n)}_\emptyset$.
Its total symbol is a homomorphism of superalgebras $\Jet^n(\cR) \rightarrow \Lambda_n \otimes \O$. Conversely, each $\Psi^{(n)} \in \Hom_{\mathrm{SAlg}}(\Jet^n(\cR), \Lambda_n \otimes \O)$ induces an $\Lambda_n$-point of $ \qq{\SC}^\infty(X,Y)$, defined by $\Psi^{(n)} \circ \jet^n \in \Hom_{\mathrm{SAlg}}(\cR, \Lambda_n \otimes \O)$. 
\end{Cor} 

\begin{Rem} \label{RemFurtherSimp}
Writing $\Phi^{(n)} = \Phi^{(n)}_\emptyset + \Phi^{(n)}_{\nil}$ w.r.t. the decomposition $\Lambda_n = \R \oplus \Lambda_n^{\nil}$, we may use the homomorphism property from Corollary \ref{CorStructComp2} to show that the $\Lambda_n$-points of $\qq{\SC}^\infty(X,Y)$ correspond bijectively to pairs $(\Phi^{(n)}_\emptyset, \Phi^{(n)}_{\nil})$, where $\Phi^{(n)}_{\nil} \in \T_Y \otimes_{\scriptscriptstyle \Phi^{(n)}_\emptyset} (\Lambda_n^{\nil} \otimes \O)$. This characterization does no longer contain a multiplicativity condition. We will not give a proof of this statement because we will not use it and the proof would require to deal with the subtleties of superjet bundles mentioned in Remark \ref{RemJetProblem}. A similar result (in a different setting) will be discussed in detail in Theorem \ref{ThmPointsSCInfty} below.         
\end{Rem}

\vspace{\baselineskip}

We next present a similar description of $\qq{\SC}^\infty(X,Y)(\Lambda_n)$ entirely in terms of the vector bundles $TX \cong T\tilde{X} \oplus E$ and $TY \cong T\tilde{Y} \oplus F$, making a \emph{non-canonical} choice for the Batchelor bundles $E$ and $F$ (cf. Theorem \ref{ThmBachelor} ff). We will use this result in chapter \ref{ChapSmooth} to define a natural smooth structure on the sets $\qq{\SC}(X,Y)(\Lambda_n)$.\\ 

Starting with $n=0$, recall that a morphism $\Phi: X \rightarrow Y$ is given by $\varphi \in C^\infty(\tilde{X},\tilde{Y})$ and a morphism $\Phi^\ast : \Gamma(\tL F^\ast) \rightarrow \Gamma(\tL E^\ast)$ acting on the superalgebras of global superfunctions. Since $\Phi^\ast$ preserves the $\Z_2$-grading, we have
\begin{align}
 \Phi^\ast(\Gamma(\tL^{\ev} F^\ast)) &\subset \Gamma(\tL^{\ev} E^\ast) & 
 \Phi^\ast(\Gamma(\tL^{\odd} F^\ast)) &\subset \Gamma(\tL^{\odd} E^\ast) \label{EqPresGrading}
\end{align}

We will prove that $\Phi^\ast$ is a differential operator $\Gamma(\tL F^\ast) \rightarrow \Gamma(\tL E^\ast)$ along $\varphi$. Recall that $\Gamma(\tL F^\ast)$ carries its usual structure as a $C^\infty_{\tilde{Y}}$-module and $\Gamma(\tL E^\ast)$ becomes an $C^\infty_{\tilde{Y}}$-module $\Gamma(\tL E^\ast)_\varphi$ by $g \cdot e := (\varphi^\ast g) e$. In this way, we can speak of differential operators acting on smooth sections of vector bundles along $\varphi$ as defined in Definition \ref{DefDoAlong}. Since the commutation rules appearing in this definition now only involve functions $f \in C^\infty(\tilde{Y})$, these operators are ordinary differential operators acting on sections of vector bundles. The super structure only enters in form of the multiplicativity of $\Phi^\ast$ and \eqref{EqPresGrading}.\\

Since being a differential operator is a local property, we may assume that the bundle $E$ is a trivial. Let $E^{(q)}$ denote the trivial vector bundle of rank $q$ over $\tilde{X}$ and let $\{e^1,\ldots,e^{q}\}$ be a frame for $E^{(q)\ast}$. In analogy to \eqref{EqFineStr1}, we have
\begin{align} \label{EqMorphDec}
\Phi^\ast(f) &= \sum_{I \in \Z_2^{q}} \Phi_I(f) e^I & \text{ for } f&\in \Gamma(\tL F^\ast)  
\end{align}
and the conditions \eqref{EqPresGrading} are clearly equivalent to $\Phi_I({\Gamma(\tL^{\ev (\odd)}F^\ast)}) = 0$ for $|I|$ odd (even). We will use the following lemma to deduce the properties of  $\Phi^\ast$ from those of the $\Phi_I$:

\begin{Lemma} \label{LemDOConstr}
Let $\Delta : \Gamma(\tL F^\ast) \rightarrow \Gamma(\tL E^\ast)$ be a differential operator of order $\leq k$ along $\varphi : \tilde{X} \rightarrow \tilde{Y} $ and $e \in \Gamma(\tL E^\ast)$. Then, the  map
\begin{align*}
\Delta\cdot e : \Gamma(\tL F^\ast) &\map \Gamma(\tL E^\ast), &
(\Delta\cdot e) (\sigma) &:= (-1)^{|e||\sigma|}(\Delta(\sigma))e 
\end{align*}
is also a differential operator of order $\leq k$. The same holds for maps of the form $e \cdot \Delta$. 
\end{Lemma}

The proof is a straightforward calculation using Definition \ref{DefDO} and \eqref{EqDefModulHom}, we omit the details. Recall that $\Diff^k(\Gamma(F), \varphi_\ast\Gamma(E)_\varphi))$ denotes the module of differential operators along $\varphi$ acting on sections of $F$ with values in sections of $E$. Theorem \ref{ThmStructComp} and Lemma \ref{LemDOConstr} imply that each $\Phi_I$ is a differential operator of order $\leq |I|$ along $\varphi$. Hence, $\Phi^\ast$ is a differential operator of order $q= \mathrm{rk}(E)$ since this is the maximal order appearing in \eqref{EqMorphDec}. However, this statement is not optimal and can be improved in the following way:

\begin{Thm} \label{ThmStructComp2}
Let $\Phi : X \rightarrow Y$ be a morphism of supermanifolds. Then, the coefficient functions $\Phi_I$ appearing in \eqref{EqMorphDec} are differential operators along $\varphi $ of order $\lfloor |I|/2 \rfloor$, where $\varphi$ is given by $\Phi_\emptyset = \varphi^\ast$. In particular, $\Phi^\ast$ is a differential operator of order $\lfloor q/2 \rfloor$. 
\end{Thm}

Note that the last statement makes sense independently of the existence of a global frame on $E$ and hence is valid on a general supermanifold.\\ 

\begin{Proof}{}
By Lemma \ref{LemDOConstr}, it is sufficient to show that the $\Phi_I$ are differential operators along $\varphi$ of order $\leq \lfloor |I|/2 \rfloor$. We will use an argument similar to the proof of Theorem \ref{ThmStructComp}. As noted before \eqref{EqMorphDec}, we can assume the $E = E^{(q)}$. We do induction on $q = \mathrm{rk}(E^{(q)})$.\\

The case $q=0$ is clear since $\Phi^\ast = \varphi^\ast$. To do the induction step $q \rightarrow q+1$, we decompose $\Phi^\ast =: \iota_{q+1}^\ast \Phi^\ast + e^{q+1} D :\Gamma(\tL F^\ast) \rightarrow \Gamma(\tL E^{(q+1)\ast})$ as in \eqref{EqRewriteMorph}. In the following, we identify $E^{(q)\ast}$ with a subbundle of $E^{(q+1)\ast}$ in the obvious way. As in the proof of Theorem \ref{ThmStructComp}, $\iota_{q+1}^\ast \Phi^\ast$ gives the morphism of supermanifolds $\Phi\iota_{q+1} : (\tilde{X}, \Gamma(\tL E^{(q)\ast})) \rightarrow Y$ and $D$ is an odd derivation along it. Choosing local coordinates $(\xi^k)$ on $Y$ , we obtain from Proposition \ref{PropDOAlong} and \eqref{EqDefCoordDos}: 
\begin{align} \label{EqDAnalyse}
D & = \sum_{k=1}^{p'+q'} h^k \otimes_{\Phi\iota_{q+1}} \tfrac{\del}{\del \xi^k}
       = \sum_{k=1}^{p'+q'} \sum_{I \in \Z_2^{q}} h^k e^I \Phi_I \circ \tfrac{\del }{\del \xi^k}
\end{align}
For the coefficients $h^k \in \Gamma(\tL E^{(q)\ast})$, we have $|h^k| = |\xi^k| + 1$ because $D$ is odd. By induction hypothesis applied to $\Phi\iota_{q+1}$, $\Phi_I$ is a differential operator along $\varphi$ of order $\leq \lfloor |I|/2 \rfloor$ for each $I \subset \Z_2^{q}$. By reordering the frame $\{e^i\}$, the argument actually applies to all $\Phi_I$ such that $|I| \leq q$. It remains to prove that $\Phi_{(1,\ldots,1)}$ is a differential operator of order $\leq \lfloor (q+1)/2\rfloor$, where $(1,\ldots,1) \in \Z_2^{q+1}$. Decomposing $h^k = \sum_{J \in \Z_2^{q}} h_{J}^k e^J$, we obtain the following expression from \eqref{EqDAnalyse}: 
\begin{align} \label{EqTopTerm}
\Phi_{(1,\ldots,1)} &= \sum_{k=1}^{p'+q'} \sum_{I + J = (1,\ldots,1) \in \Z_2^q} (-1)^r \sign(J,I) e^1\cdots e^{q+1} h^k_J \Phi_I \circ \tfrac{\del}{\del \xi^k}
\end{align}
We now have to distinguish two cases: First assume, that $q+1$ is even, i.e. $q+1=2s$. The maximal order of the differential operators $\Phi_I$ occurring in \eqref{EqTopTerm} is then $\lfloor q/2 \rfloor = s-1$ by hypothesis. Since $\tfrac{\del}{\del \xi^k}$ is a first order operator, $\Phi_{(1,\ldots,1)}$ is then a differential operator of order $s$ which equals $\lfloor (q+1)/2 \rfloor$.\\
Now assume that $q+1$ is odd, i.e. $q+1 = 2s+1$. All summands in \eqref{EqTopTerm} with $|I| \leq q-1$ yield contributions of order $\leq \lfloor (q-1)/2 \rfloor +1 = s = \lfloor (q+1)/2\rfloor$. If $|I| = q$, the only component of $h^k$ contributing to \eqref{EqTopTerm} is $h^k_{(0,\ldots,0)}$. Since $h^k_{(0,\ldots,0)}\neq 0$ implies that $h^k$ is even, the parity relation below \eqref{EqDAnalyse} shows that there may be terms of the form $\Phi_I(\del g / \del \theta^k)$ but none of the form $\Phi_I(\del g / \del x^k)$. $\tfrac{\del}{\del \theta^k}$ commutes with functions on $\tilde{Y}$ and hence is a differential operator of order $0$ on $\Gamma(\tL F^\ast)$. This implies that also for $|I|=q$, all the summands in \eqref{EqTopTerm} are of order $\leq \lfloor q/2 \rfloor = s$. Since $\lfloor (q+1) / 2 \rfloor = s$, this finishes the proof. 
\end{Proof}

Using the direct sum structure $\tL E^\ast = \bigoplus_k \tL^k E^\ast$, we introduce the following projections:
\begin{align*}
 \pr^k &: \tL E^\ast \smap \tL^k E^\ast &
 \pr^{\odd} &:= \sum_{k \in 2\N_0+1} \pr^k &
 \pr^{\ev, \geq 2} &:= \sum_{k \in 2\N} \pr^k
\end{align*}
Next to $\varphi \in C^\infty(\tilde{X},\tilde{Y})$, any morphism $\Phi : X \rightarrow Y$ gives rise to the $\R$-linear maps
\begin{align} \label{EqDefD0D1}
D_0 &:= \pr^{\ev \geq 2} \circ \Phi^\ast|_{C^\infty_{\tilde{}}} : C^\infty_{\tilde{Y}} \rightarrow  \Gamma(\tL^{\ev \geq 2} E^\ast) \\
D_1 &:= \pr^{\odd} \circ \Phi^\ast|_{\Gamma(F^\ast)}: \Gamma(F^\ast) \rightarrow  \Gamma(\tL^{\odd} E^\ast) \notag 
\end{align}

Since $C^\infty_{\tilde{Y}}$ and $\Gamma(F^\ast)$ generate $\Gamma(\tL F^\ast)$ multiplicatively, Theorem \ref{ThmStructComp2} now implies 

\begin{Cor}
Let $r := \lfloor q/2 \rfloor$. The morphisms $\Phi : X \rightarrow Y$ of supermanifolds are in one to one correspondence with triples $(\varphi, D_0, D_1)$ such that
\begin{align*}
\varphi & \in C^\infty(\tilde{X},\tilde{Y}) &  
    D_0 & \in \Diff^r(C^\infty_{\tilde{Y}},\varphi_\ast\Gamma(\tL^{\ev \geq 2} E^\ast)_\varphi) &
    D_1 & \in \Diff^r(\Gamma(F^\ast),\varphi_\ast\Gamma(\tL^{\odd} E^\ast)_\varphi),
\end{align*}
subject to the following relations:
\begin{align}
D_0(fg) &= \varphi^\ast f D_0(g) + D_0(f) \varphi^\ast g + D_0(f)D_0(g) \notag \\
D_1(f\sigma) &= \varphi^\ast f D_1(\sigma) + D_0(f)D_1(\sigma)  \label{EqCompBatDo}
\end{align}
\end{Cor}

Obviously neither $D_0$ nor $D_1$ can be induced by a bundle map $\tL F^\ast \rightarrow \tL E^\ast$. They must be higher order differential operators because the extra summands in \eqref{EqCompBatDo} spoil the $C^\infty_{\tilde{Y}}$-linearity.\\

For $\varphi \in C^\infty(\tilde{X},\tilde{Y})$, we have a canonical map 
\begin{align*}
\mathrm{I}^\varphi : \varphi^\ast \Jet^l(\tL  F^\ast) &\twoheadrightarrow \tL^0 E^\ast \subset \tL E^\ast, & 
(x, \jet^l_{\varphi(x)}(f)) &\mapsto \pr^0(f(x))
\end{align*}
Since every differential operator $\Gamma(\tL F^\ast) \rightarrow \Gamma(\tL E^\ast)$ of order $\leq l$ factors over the respective module $\Gamma(\Jet^l(\tL F^\ast))$ of jets, we obtain the following corollary: 

\begin{Cor} \label{CorDifferentChar}
Let $r := \lfloor q/2 \rfloor$. $\qq{\SC}^\infty(X,Y)(\R)$ is in bijection to the following sets: 
\begin{enumerate}
 \item Pairs $(\varphi, D)$, where $\varphi \in C^\infty(\tilde{X},\tilde{Y})$ and $D \in \Gamma(\varphi^\ast \Jet^r(\tL 
       F^\ast)) \rightarrow \Gamma(\tL^{\nil} E^\ast)$ such that $\mathrm{I}_\ast^\varphi+D$ is $C^\infty_{\tilde{X}}$-linear, multiplicative, unital and parity preserving.
 \item Pairs $(\varphi, H)$, where $\varphi \in C^\infty(\tilde{X},\tilde{Y})$ and $H : \varphi^\ast \Jet^{r}(\tL F^\ast) \rightarrow 
       \tL^{\nil} E^\ast$ is a morphism of vector bundles s.t.
       $\mathrm{I}^\varphi+ H$ is a unital, parity-preserving morphism of algebra bundles.
\end{enumerate}
\end{Cor}

\begin{Proof}{}
It is clear that (b) is just the reformulation of (a) using maps between bundles instead of maps between the corresponding spaces of sections. It remains to prove the first bijection. Since $\pr^0 \circ \Phi^\ast = \varphi^\ast$ defines a smooth map $\varphi$, setting $\widetilde{D}:= \pr^{\nil} \circ \Phi^\ast$ yields 
\begin{align} \label{EqDecomPhiAst}
\Phi^\ast = \varphi^\ast + \widetilde{D}
\end{align}
where $\widetilde{D} \in \Diff^{r}(\Gamma(\tL F^\ast),\Gamma(\tL^{\nil} E^\ast)_\varphi)$. Next, note that for any vector bundle $S \rightarrow \tilde{Y}$ we have an isomorphism of $\varphi_\ast C^\infty_{\tilde{X}}$-module-sheaves $\Gamma(S) \otimes_{C^\infty_{\tilde{Y}}} (\varphi_\ast C^\infty_{\tilde{X}})_\varphi \cong \varphi_\ast\Gamma(\varphi^\ast S)$ and in particular, $\Gamma(S) \otimes_{C^\infty(\tilde{Y})} C^\infty(\tilde{X})_\varphi \cong \Gamma(\varphi^\ast S)$. Together with \eqref{EqJetFactorProp}, this yields the following identification of spaces: 
\begin{align} \label{EqJetIso}
\Diff^{r}(\Gamma(\tL F^\ast),\Gamma(\tL E^\ast)_\varphi) 
&\cong \Gamma(\Jet^r(\tL F^\ast)) \otimes_{C^\infty_{\tilde{Y}}} \Gamma(\tL E^\ast)_\varphi \notag \\
&\cong (\Gamma(\Jet^r(\tL F^\ast)) \otimes_{C^\infty_{\tilde{Y}}}  ((C^\infty_{\tilde{X}})_\varphi \otimes_{C^\infty_{\tilde{X}}} \Gamma(\tL E^\ast)) \notag \\
&\cong \Gamma(\varphi^\ast \Jet^r(\tL F^\ast)) \otimes_{C^\infty_{\tilde{X}}} \Gamma(\tL E^\ast))  
\end{align}
Thus, we can identify $\widetilde{D}$ with an element $D \in \Gamma(\varphi^\ast \Jet^r(\tL F^\ast)) \otimes_{C^\infty_{\tilde{X}}} \Gamma(\tL^{\nil} E^\ast))$. Under the same identification, the first summand $\varphi^\ast$ in \eqref{EqDecomPhiAst} is mapped to $\mathrm{I}^\varphi_\ast$.\\
Since $\Phi^\ast$ is multiplicative, unital and parity-preserving, the same holds for $\mathrm{I}^\varphi_{\ast} + D$. Conversely, given $\varphi$ and $D$ satisfying the required properties, we can set $\iota^\ast \Phi^\ast = \varphi^\ast$ and use \eqref{EqJetIso} backwards to construct $\pr^{\nil} \circ \Phi^\ast = \widetilde{D}$ out of $D$ and $\varphi$.  $\varphi^\ast + \widetilde{D}$ then defines the required morphism $\Phi : X \rightarrow Y$. Both constructions are clearly inverse to each other.
\end{Proof}

To remove the multiplicativity condition from the characterisations in Corollary \ref{CorDifferentChar}, recall from Corollary \ref{CorIdentJet} that we may identify $\Jet^r(\tL F^\ast) \cong \Sym^{\leq r} T^\ast \tilde{Y} \otimes \tL F^\ast$ after the choice of connections on $F$ and $T \tilde{Y}$\footnote{We discuss the multiplicativity condition on the level of jets rather than differential operators, since $\Jet^r(\tL F^\ast)$ inherits a natural algebra structure from $\tL F^\ast$ whereas there is no counterpart on the corresponding space of differential operators.}. By the universal property of the symmetric algebra (see e.g. \cite{Eis1} A.2.3), we may identify the unital super algebra homomorphisms $\Sym(V^\ast) \rightarrow A$ for any super vector space $V$ and supercommutative super algebra $A$ in the following way:  
\begin{align} \label{EqSymHom}
\Hom_{\mathrm{SAlg}}(\Sym(V^\ast),A) &\overset{\sim}{\rightarrow} (V \otimes A)_\0, & \xi &\mapsto \xi|_V 
\end{align}
Choosing a basis $\{v_i\}$ of $V$ with dual basis $\{v^i\}$, the inverse of this bijection is given by multiplicative extension, i.e.  
\begin{align} \label{EqSymHomInverse}
\exp: (V \otimes A)_\0 &\rightarrow \Hom_{\mathrm{SAlg}}(\Sym(V^\ast),A), &  \tsum_i v_i \otimes \mu^i &\mapsto (\tsum_I c_I v^I \mapsto \sum_I c_I \mu^I)  
\end{align}
A more explicit expression for the inverse is obtained as follows: Since $V \otimes A \subset \Sym(V) \otimes A$, we may form $\exp(\mu) := \sum_{l=0}^\infty \tfrac{1}{l!} \mu^l$. Even though this sum is formally infinite, its contraction with elements of $\Sym(V^\ast)$ is clearly well-defined. It is straightforward to verify that for $\sum_I c_I v^I \in \Sym(V^\ast)$, we have
\begin{align*}
\langle \exp(\mu), \tsum_I c_I v^I \rangle &= \tsum_I c_I \mu^I 
\end{align*}
where $\langle, \rangle$ denotes the pairing $\Sym(V) \otimes \Sym(V^\ast) \otimes A \rightarrow A$. This justifies the name $\exp$ for the map in \eqref{EqSymHomInverse}. We can now prove the following result (see \cite{BK}, Prop. 6 for a related result obtained in a different fashion):

\begin{Thm} \label{ThmPointsSCInfty}
There are bijections
\begin{align*}
 \qq{\SC}^\infty(X,Y)(\R) 
 &\cong \{(\varphi, \sigma) \mid \varphi \in C^\infty(\tilde{X},\tilde{Y}),  \sigma \in \Gamma(\Hom_\R(\varphi^\ast (T^\ast Y), \tL^{\nil} E^\ast))\}  \\
 &\cong \{(\varphi, \tilde{\sigma}) \mid \varphi \in C^\infty(\tilde{X},\tilde{Y}), \tilde{\sigma} \in  \Gamma((\varphi^\ast TY \otimes \tL^{\nil} E^\ast )_\0\} \notag
\end{align*}
\end{Thm}

\begin{Proof}{}
It is clearly enough to prove the claim for a fixed underlying map $\varphi \in C^\infty(\tilde{X},\tilde{Y})$. Let $r := \lfloor q/2 \rfloor$. By Corollary \ref{CorDifferentChar} (b) and Corollary \ref{CorIdentJet} (applied with $A = \varphi^\ast \tL F^\ast$) we have to classify bundle maps $\psi : \varphi^\ast \Jet^{r}(\tL F^\ast) \rightarrow \tL^{\nil} E^\ast$ s.t. 
\begin{align*}
\mathrm{I}^\varphi + \psi : \varphi^\ast \Jet^{r}(\tL F^\ast) 
\cong \Sym^{r} T^\ast \tilde{Y} \otimes \tL F^\ast 
&\map \tL E^\ast 
\end{align*}
is a homomorphism of algebra bundles over $\tilde{X}$. Since the symmetric algebra of the $\Z_2$-graded vector bundle $T^\ast Y \cong T^\ast \tilde{Y} \oplus F^\ast$ can be identified with $\Sym(T^\ast\tilde{Y}) \otimes \tL F^\ast$ (\cite{Eis1}, Appendix A2.3), we obtain 
\begin{align}\label{EqAlgBunQuot}
 \Sym^{\leq r} T^\ast \tilde{Y} \otimes \tL F^\ast  &\cong 
 \Sym (T^\ast Y) / \bigl( \Sym^{> r} T^\ast \tilde{Y} \otimes \tL F^\ast \bigr)
\end{align}
By the universal property of the symmetric algebra, we have 
\begin{align} \label{EqUniSymm}
 \Hom_{\mathrm{SAlg}}(\Sym (T^\ast Y), \tL E^\ast) & \cong \Hom_\R(T^\ast Y, \tL E^\ast), & 
         \beta & \mapsto \beta|_{T^\ast Y},                
\end{align}
Denoting the underlying smooth map of $\Phi$ by $\tilde{\Phi}$, we therefore have an injective map
\begin{align} \label{EqMapInjective}
  \{\Phi \in \Hom_{\sss \BKL}(X,Y) \mid \tilde{\Phi} = \varphi\} &\rightarrow \Hom_\R(\varphi^\ast(T^\ast Y), \tL^{\nil} E^\ast) \\
 \mathrm{I}^\varphi + \psi  & \mapsto (\mathrm{I}^\varphi + \psi)|_{T^\ast Y} =  \psi|_{T^\ast Y} \notag
\end{align}
To prove that this map is surjective, note that every $\beta' \in \Hom_\R(\varphi^\ast(T^\ast Y), \tL^{\nil} E^\ast)$ gives rise to an algebra bundle homomorphism $\beta : \Sym (T^\ast Y) \rightarrow \tL E^\ast$ by \eqref{EqUniSymm}. It remains to show that $\beta$ descends to the quotient from \eqref{EqAlgBunQuot}, i.e. that $\beta(p) = 0 $ for all $p \in \Sym^{> r} T^\ast \tilde{Y} \otimes \tL F^\ast$. Let $p = dy^{i_1} \cdots dy^{i_k} \otimes f$ for some $f \in \tL F^\ast$ and $k \geq r + 1$. We compute the Grassmann-$\Z$-degree $\mathrm{deg}(\beta(p))$. Since $\beta'$ preserves parity, we have $\beta'(dx^i) \in \tL^{\ev \geq 2} E^\ast$. Thus 
\begin{align*}
 \mathrm{deg(\beta(p))} 
&\geq 2k \geq 2(\lfloor \tfrac{q}{2} \rfloor + 1) 
 = \begin{cases} q+2 & q \in 2\N\\ q+1 & q \in 2\N+1 \end{cases} \quad > q
\end{align*}
Since $q = \mathrm{rk}(E^\ast)$, it follows that $0 = \beta(p) \in \tL E^\ast$ and we conclude that the map \eqref{EqMapInjective} is bijective.\\
Identifying vector bundle morphisms with sections of the corresponding bundle of homomorphisms, the first bijection of the proposition now follows from $\eqref{EqMapInjective}$ when $\varphi$ varies over $C^\infty(\tilde{X},\tilde{Y})$. The second bijection is clear.
\end{Proof}

To describe higher points of $\qq{\SC}^\infty(X,Y)$, i.e. the sets $\qq{\SC}^\infty(X,Y)(\Lambda_n) = \Hom_{\sss \BKL}(\PR_n \times X, Y)$, we note that the sheaf of superfunctions on the product $\PR_n \times X$ (cf. Remark \ref{RemProducts}) is given by  
\begin{align*}
\Lambda_n \otimes_\R \Gamma(\tL E^\ast) &\cong \Gamma(\Lambda_{n} \otimes \tL E^\ast) \cong \Gamma(\tL(\R^n \oplus E^\ast))  
\end{align*}
where $\Lambda_{n}$ and $\R^n$ denote the vector spaces as well as the induced trivial bundles over $\tilde{X}$ built out of them. Thus, replacing $E$ by $\R^n \oplus E$, Theorem \ref{ThmPointsSCInfty} implies the following statement:

\begin{Cor} \label{CorPointsSCInfty}
We have bijections
\begin{align*}
\qq{\SC}^\infty(X,Y)(\Lambda_n)
  &\cong \{(\varphi, \sigma) \mid \varphi \in C^\infty(\tilde{X},\tilde{Y}),  \sigma \in \Gamma(\Hom_\R(\varphi^\ast (T^\ast Y), \tL^{\nil} (\R^n \oplus E^\ast)))\}  \\
  &\cong \{(\varphi, \tilde{\sigma}) \mid \varphi \in C^\infty(\tilde{X},\tilde{Y}), \tilde{\sigma} \in  \Gamma((\varphi^\ast TY \otimes \tL^{\nil} (\R^n \oplus E^\ast))_\0 )\} \notag
\end{align*}
\end{Cor}

This description suggests to interpret $\qq{\SC}^\infty(X,Y)(\Lambda_n)$ as an infinite dimensional bundle over $C^\infty(\tilde{X},\tilde{Y})$ with fibre at $\varphi$ given by $\Gamma((TY \otimes \tL^{\nil} (\R^n \oplus E^\ast))_\0)$. We will come back to this idea in the last section.


\section{Finite-dimensional supermanifolds revisited} \label{ChapFinDim}

We will now use Corollary \ref{CorPointsSCInfty} to obtain a description of finite-dimensional supermanifolds and their morphisms on the level of $\Lambda_n$-points. In particular, this provides us with a description of atlases on supermanifolds, that generalizes to the infinite-dimensional situation (cf. \eqref{EqPointwiseDefChart}).\\

A geometric description of the $\Lambda_n$-points of a finite-dimensional supermanifold $Y$ is readily obtained from Corollary \ref{CorPointsSCInfty}. Putting $X := \{0\}$, identifying $Y \cong (\tilde{Y}, \Gamma(\tL F^\ast))$ as before and noting that $C^\infty(\{0\},\tilde{Y}) \cong \tilde{Y}$, we get
\begin{align} \label{EqFinDimObjId}
\SC^\infty(\PR_n,Y) \cong  \Hom_\BKL(\PR_n,Y) &\cong \{(y,\nu) \mid y \in \tilde{Y}, \nu \in (T_y Y \otimes \Lambda_n^{\nil} )_\0 \} 
\cong (TY  \otimes \Lambda_n^{\nil})_\0
\end{align}
In particular, this set carries the structure of a finite-dimensional smooth manifold induced by the bundle $(TY \otimes \Lambda_n^{\nil} )_\0$, whose smooth structure is determined by $Y$. Since $\Hom_\BKL(\PR_n,Y) \cong \Hom_{\mathrm{SAlg}}(\Gamma(\tL F^\ast),\Lambda_n)$, each element $\nu \in (TY  \otimes \Lambda_n^{\nil})_\0$ induces a map $\hat{\nu} : \Gamma(\tL F^\ast) \rightarrow \Lambda_n$. By Corollary \ref{CorDifferentChar}, this map factors over $\jet^{r}$, where now $r := \lfloor n/2 \rfloor$, and we may identify $\jet^r(\sigma) \cong S^r(\sigma)$ for $\sigma \in \Gamma(\tL F^\ast)$. Recall that $\mathrm{dim}(Y) = p'|q'$. Using the bijection \eqref{EqSymHom}, \eqref{EqSymHomInverse} as in the proof of Theorem \ref{ThmPointsSCInfty}, we find the following explicit expression: 
\begin{align} \label{EqExplicitNPoint}
\hat{\nu}(\sigma) &= \langle \exp(\nu), S^r(\sigma) \rangle = \sum_{ I \in \N_0^{p'\!,r}, J \in \Z_2^{q'} }
                                                                       \tfrac{1}{I!} (D_I \sigma)|_{\tilde{\nu}}(w_J) \nu_\2^I \nu_\1^J
\end{align}
In the last step, we used the notation from Remark \ref{RemFlatTaylor} where $(w_{1},\ldots,w_{p'+q'})$ is an adapted basis of $T_{\tilde{\nu}}Y$ and for $J \in \Z_2^{q'}$, we set $w_J := w_{\scriptscriptstyle{p'+1}}^{\scriptscriptstyle{j_1}} \wedge \cdots \wedge w_{\scriptscriptstyle{p'+q'}}^{\scriptscriptstyle{j_{q'}}}$. $\tilde{\nu} \in \tilde{Y}$ denotes the basepoint of $\nu$ and $\nu^i_\2$ ($\nu^i_\1$) the even (odd) part of $\nu^i \in \Lambda_n^{\nil}$, defined by the decomposition $\nu = \tsum_i w_i \otimes \nu^i$. In the case that $Y = \R^{p'|q'}$, let $(x^i,\theta^\alpha)$ be the canonical coordinates. We may expand any $\sigma \in \Gamma(\tL F^\ast) \cong C^\infty(\R^{p'},\tL \R^{q'})$ as $\sigma(x) = \sum_{J \in \Z_2^{q'}} \sigma_J(x)\theta^J$. \eqref{EqExplicitNPoint} now takes the form 
\begin{align*}
\hat{\nu}(\sigma) &= \sum_{I \in \N_0^{p'\!,r}, J \in \Z_2^{q'}} \tfrac{1}{I!} \tfrac{\del^{|I|}\sigma_J}{\del x^I}\bigl|_{\tilde{\nu}} \nu_{\2}^I\nu_{\1}^J 
\end{align*}
This is exactly the result obtained in Proposition 3.13 of \cite{BCF1}.\\

Next, we describe the map $\Phi_\ast$ on $\Lambda_n$-points, induced by a morphism $\Phi \in \Hom_\BKL(X,Y)$ with underlying smooth map $\varphi \in C^\infty(\tilde{X},\tilde{Y})$. By definition, we have
\begin{align*}
\Phi_\ast : \Hom(\PR_n,X) &\rightarrow \Hom(\PR_n,Y), & \hat{\mu} &\mapsto \Phi_\ast \hat{\mu} = \hat{\mu} \circ \Phi^\ast
\end{align*}
where $\Phi^\ast: \cR  \rightarrow \O$ is the super algebra homomorphism defining $\Phi$. To derive an explicit expression using the characterisations of $\Lambda_n$-points from \eqref{EqFinDimObjId}, \eqref{EqExplicitNPoint}, we also identify $X \cong (\tilde{X}, \Gamma(\tL E^\ast))$. The image $\nu \in (TY \otimes \Lambda_n^{\nil})_\0$ of $\mu \in (TX \otimes \Lambda_n^{\nil})_\0$ under $\Phi$, defined by $\hat{\nu} := \hat{\mu} \circ \Phi^\ast$, is then uniquely characterized by
\begin{align} \label{EqComputeRho}
\langle \exp(\nu), S^r(\sigma) \rangle &= \hat{\nu}(\sigma) = \hat{\mu}(\Phi^\ast(\sigma)) = \langle \exp(\mu), S^r(\Phi^\ast(\sigma))\rangle 
\end{align}
For $\nu \in (T_{\tilde{\nu}} Y \otimes \Lambda_n^{\nil})_\0$, we clearly have $\nu = \sum_j w_j \otimes \langle \exp(\nu), w^j \rangle$. We may choose local functions $y^j$ ($j=1,\ldots,p'$) on $\tilde{Y}$ and local sections $y^\alpha$ ($\alpha = p'+1,\ldots,p'+q'$) in $F^\ast$ s.t. $dy^j (\tilde{\nu}) = w^j$, $y^j(\tilde{\nu}) = 0$, $y^\alpha(\tilde{\nu}) = w^\alpha$ and s.t. $dy^j, y^\alpha$ are parallel. We clearly have $S^r(y^j)(\tilde{\nu}) = w^j$ ($j=1,\ldots,p'+q'$) and thus obtain from \eqref{EqComputeRho}
\begin{align} \label{EqPhiAstSpecial}
&\Phi_\ast : (TX  \otimes \Lambda_n^{\nil})_\0 \rightarrow (TY  \otimes \Lambda_n^{\nil})_\0, \\
&\Phi_\ast(\mu) = \sum_{j=1}^{p'+q'} w_j \otimes \langle \exp(\mu), S^r(\Phi^\ast(y^j)) \rangle \notag
= \sum_{j=1}^{p'+q'} w_j \otimes \sum_{I \in \N_0^{p,r}} \sum_{J \in \Z_2^{q}} \tfrac{1}{I!} D_I(\Phi^\ast(y^j))|_{\tilde{\mu}}(v_J) \mu_{\2}^I \mu_\1^J
\end{align}
Here, the decomposition of $\mu$ is defined w.r.t an adapted basis $(v_1\ldots,v_{p+q})$ of $TX$ and the ``frame'' $\{y^j\}$ is defined in a neighbourhood of $\tilde{\nu} = \varphi(\tilde{\mu})$. Again, if we specialize to $X = \R^{p|q}$ and $Y = \R^{p'|q'}$, we reproduce the result in \cite{BCF1}, formula (3.3).\\

We will be mostly interested in the case that $\Phi$ is induced by a bundle morphism $\phi : E \rightarrow F$ over $\varphi \in C^\infty(\tilde{X},\tilde{Y})$. More precisely, $\Phi^\ast$ is given by the superalgebra homomorphism
\begin{align} \label{EqSpeciaMor}
\xymatrix{
\Gamma(\tilde{Y},\tL F^\ast) \ar[r]^-{\varphi^\ast} & \Gamma(\tilde{X},\varphi^\ast \tL F^\ast) \ar[r]^-{(\wedge \phi^\ast)_\ast} & \Gamma(\tilde{X}, \tL E^\ast) }
\end{align}  
where $\wedge \phi^\ast : \varphi^\ast \tL F^\ast \rightarrow \tL E^\ast$ is a bundle morphism over $\tilde{X}$. Note that $\ast$ denotes either pullback or fibrewise dualisation, depending on the context. Using the ``frame'' $\{y^j\}$ constructed above, we clearly have
$\Phi^\ast(y^j) = y^j \circ \varphi$ ($j=1,\ldots,p'$) and $\Phi^\ast(y^\alpha) = \phi^\ast (y^\alpha \circ \varphi)$ ($\alpha=p'+1,\ldots,p'+q'$).  Since $dy^j$ and $y^\alpha$ are parallel and $y^j(\varphi(\tilde{\mu})) = 0$, the chain rules from \eqref{EqVBHomChainRule} and \eqref{EqFaaDiBruno} yield
\begin{align*}
S^r(y^j \circ \varphi) (\tilde{\mu}) &= w^j \circ S^r(\varphi)(\tilde{\mu})   &
S^r(\phi^\ast (y^\alpha \circ \varphi))(\tilde{\mu}) &= w^\alpha \circ S^r(\phi)(\tilde{\mu})
\end{align*}
From \eqref{EqPhiAstSpecial}, we now obtain the explicit expression  
\begin{align} \label{EqPhiSimpTotal}
\Phi_\ast(\mu) &= \sum_{I \in \N_0^{p,r}} \tfrac{1}{I!} D_I \varphi|_{\tilde{\mu}} \otimes \mu_\2^I +
   \sum_{I \in \N_0^{p,r}} \sum_{\alpha = p+1}^{p+q}  \tfrac{1}{I!} D_I \phi|_{\tilde{\mu}}(v_\alpha) \otimes \mu_\2^I \mu_\1^\alpha \notag \\
   &=: S^r(\varphi)(\mu) + S^r(\phi)(\mu)
\end{align}
which is independent of any choice of frames. It is valid on all of $(TX \otimes \Lambda_n^{\nil})_\0$ and not just in a neighbourhood inside $\tilde{X}$.\\ 

The case of a general morphism $\Phi$ is slightly more involved since the homomorphism $\Phi^\ast$ is no longer induced by a bundle morphism in the sense of \eqref{EqSpeciaMor}. By Theorem \ref{ThmStructComp2} / Corollary \ref{CorDifferentChar}, $\Phi^\ast$ splits as
\begin{align} \label{EqSplitSymbMor}
\Phi^\ast : \xymatrix{ \Gamma(\tL F^\ast) \ar[r]^-{S^k} & \Gamma(\Sym^{\leq k} T^\ast Y) \ar[r]^-{\varphi^\ast} & \Gamma(\varphi^\ast\Sym^{\leq k} T^\ast Y) \ar[rr]^-{(\symb\Phi^\ast)_\ast} & & \Gamma(\tL E^\ast) }
\end{align}
where $k := \lfloor \mathrm{rk}(E) /2 \rfloor = \lfloor q /2 \rfloor$  and we used a connection on $F$ to identify jets and symmetric algebras. Since the symbol $\symb(\Phi^\ast)$ is a bundle morphism, we may write $\symb\Phi^\ast = (\symb\Phi)^\ast$ where (omitting canonical identifications) $\symb\Phi := (\symb\Phi^\ast)^\ast: \tL E \rightarrow \varphi^\ast \Sym^{\leq k} TY$ is a homomorphism of superalgebra bundles. A straightforward computation, again based on the chain rules, yields
\begin{align*}
S^r(\Phi^\ast(y^j))(\tilde{\mu}) 
&= \begin{cases}
  (w^j + w^j \circ S^r(\varphi)) \circ_r S^r(\symb \Phi)(\tilde{\mu}) & j = 1,\ldots,p'\\
   w^j \circ S^r(\symb\Phi)(\tilde{\mu}) & j = p'+1, \ldots p'+q'
  \end{cases}
\end{align*}
By \eqref{EqPhiAstSpecial}, we obtain the analogue to \eqref{EqPhiSimpTotal} for a general morphism:
\begin{align} \label{EqPhiAstGeneral}
\Phi_\ast(\mu) &= \sum_{j=1}^{p'+q'} w_j \otimes \langle \exp(\mu), w^j(S^r(\symb\Phi)(\tilde{\mu}))  \rangle 
                + \sum_{j=1}^{p'} w_j \otimes \langle \exp(\mu), w^j(S^r(\varphi) \circ_r S^r(\symb\Phi))(\tilde{\mu}) \rangle \notag \\
               &=: S^r(\symb\Phi)(\mu) + S^r(\varphi) \circ_r S^r(\symb\Phi) (\mu)
\end{align}
Again, one may use a frame on $TX$ to give an explicit expression as a polynomial in $\mu_\2$ and $\mu_\1$. If $\Phi$ is of the special type described in \eqref{EqSpeciaMor}, \eqref{EqPhiAstGeneral} clearly reduces to \eqref{EqPhiSimpTotal}. Summarizing, we have obtained the following results:

\begin{Prop} \label{PropFindimNPoints}
After the choice of a Batchelor bundle, the $\Lambda_n$-points of a finite-dimensional supermanifold $X$ can be identified with the manifold $(TX \otimes \Lambda_n^{\nil})_\0$. Under this identification, a 
a morphism $\Phi : X \rightarrow Y$ as defined in \eqref{EqSpeciaMor} induces the smooth map
\begin{align} \label{EqMorPnPoints}
\Phi_\ast : (TX \otimes \Lambda_n^{\nil})_\0 &\rightarrow (TY \otimes \Lambda_n^{\nil})_\0 &
\mu &\mapsto S^k(\varphi)(\mu) + S^k(\phi)(\mu)  
\end{align}
A general morphism is represented by \eqref{EqPhiAstGeneral} or the local expression from \eqref{EqPhiAstSpecial}. 
\end{Prop}

The following example shows s simple cases, where $\Phi$ is not of the type \eqref{EqSpeciaMor}:

\begin{Ex}
Let $X = \R^{p|2}$, $Y = (\tilde{Y}, \Gamma(\tL F^\ast))$ and let $(x^i,\theta^\alpha)$ denote the standard coordinates on $X$. A morphism $\Phi : \R^{p|2} \rightarrow Y$ of supermanifolds is uniquely determined on $g \in C^\infty(\tilde{Y})$ and $s \in \Gamma(F^\ast)$ by
\begin{align*}
\Phi^\ast(g) &= \varphi^\ast g + \xi(g)\theta^1\theta^2 \in C^\infty(\R^p,\R \oplus \tL^2(\R^2)^\ast) &
\Phi^\ast(s) &= H(s) \in C^\infty(\R^p,(\R^2)^\ast)
\end{align*}
where $\varphi \in C^\infty(\R^p,\tilde{Y})$, $\xi \in \Gamma(\varphi^\ast T\tilde{Y})$ and $H \in \Gamma(\Hom_\R(\varphi^\ast F^\ast, \R^p \times (\R^2)^\ast)))$. This is in accordance with the results from chapter \ref{ChapAlg}, since the action on odd sections (given by $H$) is a differential operator of order $0$ whereas the action on even sections is of order $\leq 1$ due to the presence of the vector field. It is clear that $\Phi^\ast$ can not be induced by a bundle morphism $\varphi^\ast \tL F^\ast \rightarrow \R^p \times (\R^2)^\ast$ unless $\xi=0$. By \eqref{EqSplitSymbMor}, $\Phi$ is uniquely determined by the smooth map $\varphi$ and the symbol $\varphi^\ast \Sym^{\leq 1}(T^\ast Y) \rightarrow \R^p \times \tL(\R^2)^\ast$ map which, being a multiplicative bundle morphism, is uniquely determined by the restriction to elements of degree one: 
\begin{align*}
\varphi^\ast T^\ast Y \cong \varphi^\ast (T^\ast \R \oplus F^\ast) &\rightarrow \R^p \times \tL(\R^2)^\ast &
(\alpha, s) &\mapsto \alpha(\xi)\theta^1\theta^2 + H(s)
\end{align*}
Note that the splitting $\Jet^1(\tL F^\ast) \cong \tL F^\ast \oplus T^\ast\tilde{Y} \otimes \tL F^\ast$ is defined without reference to a connection. However, the representation of the action of $\Phi$ on general $\Lambda_n$-points is based on the choice of connections. It is a polynomial map 
\begin{align*}
\Phi_\ast : (T\R^{p|2} \otimes \Lambda_n^{\nil})_\0 \cong (\Lambda_n^{\ev})^p \times (\Lambda_n^{\odd})^2 &\rightarrow (TY \otimes \Lambda_n^{\nil})_\0  
\end{align*}
whose explicit form, given by \eqref{EqPhiAstGeneral} or \eqref{EqPhiAstSpecial}, contains (covariant) derivatives from $\varphi, \xi$ and $H$.
\end{Ex}

\begin{Rem}
In \eqref{EqPhiAstSpecial}, we recovered what is often called the Grassmann- or superanalytic expansion of $\Phi$. This notion occurs in the Rogers-DeWitt-approach to supergeometry and in the physics literature. We will not discuss it here but refer to \cite{Rog1}, chapter 4, \cite{DeWitt}, chapter 1 or \cite{BBHR}, chapter III.2  and the references in these monographs. The (to my knowledge) first systematic study of the relation of morphism $\R^{p|q} \rightarrow \R^{p'|q'}$ and the induced maps on the corresponding $\Lambda$-points can be found in \cite{Vor}, see also \cite{BCF1} Proposition 4.13. 
Comparing the ``Taylor expansions'' with nilpotent arguments from \eqref{EqExplicitNPoint}, \eqref{EqPhiAstSpecial} with the construction of Weil functors in \cite{KM}, 31.5, it is clear that the latter concept is closely related to functorial formulation of supergeometry. Again, we will not comment this aspect but refer to \cite{BCF1} (chapters 3, 4) and \cite{All} (section 2.4) for a discussion of Weil functors / bundles in the context of supergeometry. 
\end{Rem}

We finally apply these results to charts and chart-transitions of the supermanifold $X$. Let $\phi : E|_U \rightarrow U \times \R^q$ be a local trivialization of its Batchelor bundle such there is a chart $\varphi : U \rightarrow V \subset \R^p$ of $\tilde{X}$. They induce vector the bundle morphism $(\varphi,\mathrm{Id}_{\R^q}) \circ \phi : E|_U \rightarrow V \times \R^q$. As in \eqref{EqSpeciaMor}, there is an induced morphism $\Phi : X|_U \rightarrow \R^{p|q}|_V$ which is clearly an isomorphism. Choosing two trivializations and charts, we obtain two charts $\Phi_i$ ($i=1,2$) and a chart transition 
\begin{align*}
\xymatrix{ R^{p|q}|_{V_{12}} \ar[r]^{\Phi_1^{-1}} & X_{U_1 \cap U_2} \ar[r]^{\Phi_2} & R^{p|q}|_{V_{21}}}
\end{align*}
where $V_{12} := V_1 \cap \varphi_1\varphi_2^{-1}(V_2)$ and $V_{21}:= V_2 \cap \varphi_2\varphi_1^{-1}(V_1)$. We thus get the following corollary to Proposition \ref{PropFindimNPoints}:
\begin{Cor}
The chart transition $\Phi_2\Phi_1^{-1}$ induces the following map on the level of $\Lambda_n$-points:
\begin{align} \label{EqChartTransFun}
(T\R^{p|q}|_{V_{12}} \otimes \Lambda_n^{\nil})_\0 \cong (\qq{\R}^{p|q}|_{V_{12}}(\Lambda_n)) &\longrightarrow  
(\qq{\R}^{p|q}|_{V_{21}}(\Lambda_n)) \cong (T\R^{p|q}|_{V_{21}} \otimes \Lambda_n^{\nil})_\0  \\
\nu \quad& \mapsto \quad \varphi(\tilde{\nu}) + S^r(\varphi_2\varphi_1^{-1})(\nu) + S^r(\phi_2\phi_1^{-1})(\nu) \notag
\end{align}
\end{Cor}
Again, we decomposed $\nu = \tilde{\nu} + \nu^{\nil} \in \R^{p|q} \oplus \R^{p|q} \otimes \Lambda_n^{\nil}$. $\tilde{\nu}$ represents the coordinates of the basepoint in the bundle picture. Note that the chart domains have been identified with restrictions of (finite-dimensional), representable supermodules in the sense of Definition \ref{EqRepSMod}. Being derived from the morphism $\Phi\Phi_1^{-1} \in \Hom_{\BKL}(\R^{p|q}|_{V_{12}}, \R^{p|q}|_{V_{21}})$, it follows that \eqref{EqChartTransFun} defines a supersmooth morphism $\qq{\R}^{p|q}|_{V_{12}} \rightarrow \qq{\R}^{p|q}|_{V_{21}}$ (cf. Chapter  \ref{ChapSmooth} for a similar argument in the infinite-dimensional setting). This observation is not new, in fact, Voronov already established in \cite{Vor} that both points of view are in fact equivalent.


\section{The supersmooth structure of $\qq{SC}^\infty(X,Y)$} \label{ChapSmooth}

We will now equip $\qq{SC}^\infty(X,Y)$ with the structure of a smooth supermanifold. Following the discussion in chapter \ref{ChapPre}, we have to equip the sets 
\begin{align} \label{EqSetSmoothStruct}
\qq{SC}^\infty(X,Y)(\Lambda_n) &\cong \{(f,\sigma) \mid f \in C^\infty(\tilde{X},\tilde{Y}), \sigma \in \Gamma((f^\ast TY \otimes \tL^{\nil}(\R^n \oplus E^\ast))_\0)
\end{align}
from Corollary \ref{CorPointsSCInfty} with structures of (infinite-dimensional) smooth manifolds. This will be done in the following way: For any finite-dimensional vector bundle $\pi: S \rightarrow N$ and smooth map $f \in C^\infty(M,N)$, we have the pullback diagram
\begin{align} \label{DiagPullBackSPull}
 \xymatrix{
     f^\ast S \ar[rr]^{\pi^\ast f} \ar[d]^{f^\ast\pi} & & S \ar[d]^{\pi}  \\
       M \ar[rr]^f \ar@/^1pc/[u]^s \ar[urr]^{\tilde{s}} & & N     
  }
\end{align}
where $s \in \Gamma(f^\ast S)$, $\tilde{s} := \pi_S^\ast f \circ s \in C^\infty(M,S)$ and $\pi^\ast f$ is the canonical map. It is well known that $\pi^\ast f$ is smooth, a linear isomorphism on fibres and induces a bijection 
\begin{align} \label{EqNatPropPB}
\Gamma(f^\ast S) &\overset{\sim}{\rightarrow} \{\tilde{s} \in  C^\infty(M,S) \mid \pi \circ s' = f\} =: C^\infty(M,S)_f  
\end{align}
Hence, \eqref{EqSetSmoothStruct} and \eqref{EqNatPropPB} suggest to view the set $\qq{SC}^\infty(X,Y)(\Lambda_n)$ as a subset of the space of maps $C^\infty(\tilde{X}\times \tilde{Y}, TY \boxtimes \tL^{\nil}(\R^n \oplus E^\ast))$. We will discuss the details after Remark \ref{RemMichor} and first study the smooth structure of $C^\infty(M,S)$. We note that infinite-dimensional mapping spaces in the sense of \cite{Michoralt} have also been used in \cite{BK} (chapter 3) in a different, non-functorial way to describe certain morphism spaces.\\ 

Given a Riemannian metric on $N$ as well as a fibre metric and a compatible connection $\LC$ on $S$, there is an induced Riemannian metric $G$ on $S$ defined in \eqref{EqDefMetA}. By the construction described in Theorem \ref{ThmCInftyMN}, $C^\infty(M,S)$ becomes an infinite-dimensional smooth manifold. Using the splitting in horizontal and vertical subbundles from \eqref{EqDecompTS} and \eqref{EqIdentHV}, the model space near $\sigma : M \rightarrow S$ can be identified with
\begin{align} \label{EqSplitModelSpace}
 \Gamma_c(\sigma^\ast TS) &\cong \Gamma_c(M, \sigma^\ast (\mathcal{H}(S) \oplus \mathcal{V}(S))) 
 \cong \Gamma_c(f^\ast TN \oplus f^\ast S) \cong \Gamma_c(f^\ast TN) \oplus \Gamma_c(f^\ast S)
\end{align}
Here, we have set $f := \pi \circ \sigma$. Note that the first factor is just the model space of the chart $u_f : U_f \rightarrow V_f$ from \eqref{EqCharts2}. It is well known (see e.g. \cite{KM} 42.13) that the bundle projection induces the smooth map
\begin{align*}
 \pi_\ast : C^\infty(M,S) &\rightarrow C^\infty(M,N), \enspace s \mapsto \pi \circ s
\end{align*}
Together with \eqref{EqNatPropPB}, this suggest to view $C^\infty(M,S)$ as a bundle over $C^\infty(M,N)$, whose fibre at $f$ is given by $C^\infty(M,S)_f \cong \Gamma(f^\ast S)$. This is reasonable because any map $g \in U_f$ is homotopic to $f$ by \eqref{EqCharts1}, \eqref{EqCharts2} and thus, $\Gamma(f^\ast S) \cong \Gamma(g^\ast S)$; in fact, any map $g$ in the connected component of $f$ in $C^\infty(M,N)$ is smoothly homotopic to $f$ (cf. \cite{Woc}, Corollary 14 and \cite{KM2}). If we regard $\Gamma(f^\ast S)$ as a locally affine manifold modelled on $\Gamma_c(f^\ast S)$ as in Example \ref{ExLocAffine} , the decomposition in \eqref{EqSplitModelSpace} corresponds to the directions along the base and the fibres, respectively. In fact, the charts defined relative to $G$ on $C^\infty(M,S)$ induce local trivializations in a natural way: For $f \in C^\infty(M,N)$, consider the $c^\infty$-open subsets 
\begin{align*}
 \mathbb{U}_f &:= \pi^{-1}(U_{f}) \subset C^\infty(M,S) & 
 \mathbb{V}_f &:= V_f \times \Gamma(f^\ast S) \subset \Gamma_c(f^\ast TN) \times \Gamma(f^\ast S) 
\end{align*}
$\Gamma_c(f^\ast TN) \times \Gamma(f^\ast S)$ carries the locally affine topology induced by $c^\infty(\Gamma_c(f^\ast TN) \times  \Gamma_c(f^\ast S))$, which is strictly finer than the one of $c^\infty(\Gamma_c(f^\ast TN)) \times c^\infty(\Gamma_c(f^\ast S))$, see Remark \ref{RemCInftyTop}. In particular, the characterisation of smooth curves given in Lemma \ref{LemSmoothCurvCInftyMN} applies to $\mathbb{U}_f$ as well as $\mathbb{V}_f$. The sets $\mathbb{U}_f$ cover $C^\infty(M,S)$ and they are unions of entire fibres of the projection $\pi_\ast : C^\infty(M,S) \rightarrow C^\infty(M,N)$. In analogy to \eqref{EqCharts1}, we set  
\begin{align}\label{EqDefLocTriv}
 \mathbf{u}_f : \mathbb{U}_f &\rightarrow \mathbb{V}_f  &
 \mathbf{u}_f(\sigma) &:= (x \mapsto (x, (\pi^{TS},\exp^G)^{-1}\circ (0_f,\sigma)(x)))
\end{align}
where $0_f \in C^\infty(M,S)$ is the map $x \mapsto 0_{f(x)} \in S_{f(x)}$. It is well-defined by Lemma \ref{LemConvNghdA}. Moreover, its restriction to $U_{0_f} \subset C^\infty(M,S)$ is precisely the chart $u_{0_f} : U_{0_f} \rightarrow V_f \times \Gamma_c(f^\ast S)$ near $0_f$ as defined in \eqref{EqCharts1}. Using \eqref{EqExpMapBun} (with $a = 0_{f(x)}$), we find the following representation in terms of the parallel transport:
\begin{align} \label{EqLocTrivPt}
\mathbf{u}_{f}(\sigma) &= (u_f(\pi_\ast\sigma), P_{\pi_\ast\sigma, f}(\sigma)) \\
\mathbf{u}_{f}^{-1} (\xi,\zeta) &= P_{f,u_f^{-1}(\xi)}(\zeta) \notag
\end{align}
Here, $P_{\pi_\ast\sigma, f}(\sigma)$ denotes the map $x \mapsto P_{\pi_\ast\sigma(x), f(x)}(\sigma(x))$ etc.. \eqref{EqLocTrivPt} clearly shows that $\mathbf{u}_f$ is bijective and commutes with the projections onto $C^\infty(M,N)$ and $V_f$, respectively. In fact, the maps provide a local trivialization and we have the following theorem:

\begin{Thm} \label{ThmBundleStruct}
$\pi_\ast : C^\infty(M,S) \rightarrow C^\infty(M,N)$ is a smooth infinite dimensional fibre bundle. Its typical fibre $\Gamma(f^\ast S)$ carries the structure of a locally affine manifold as described in Example \ref{ExLocAffine}. Two different choices of connections on $S$ yield isomorphic bundles.
\end{Thm}

\begin{Proof}{}
To see that the local trivializations $\mathbf{u}_f$ are in fact diffeomorphisms, we show that they map smooth curves to smooth curves (see also Remark \ref{RemChartsSmooth}). For $\gamma \in C^\infty(\R, \mathbb{U}_f)$, we have the induced map $\gamma^\wedge \in C^\infty(\R \times M, S)$. By \eqref{EqLocTrivPt},
\begin{align} \label{EqLocTrivCurve}
(\mathbf{u}_f{}_\ast \gamma)^\wedge(t,x) &= (u_f(\pi(\gamma^\wedge(t,x))), P_{\pi(\gamma^\wedge(t,x)), f(x)}(\gamma^\wedge(t,x))))
\end{align}
which is clearly smooth in $(t,x) \in \R \times M$ by the usual smoothness properties of $\exp$ and $P$ on finite-dimensional manifolds. Let $[a,b] \subset \R$ be a compact interval. Since $\gamma \in C^\infty(\R, \mathbb{U}_f)$, there is compact set $K \subset M$ s.t. for each $x \in M \setminus K$, the map $[a,b] \ni t \mapsto \gamma^\wedge(t,x)$ is constant. By \eqref{EqLocTrivCurve}, the same is then true for $[a,b] \ni t \mapsto (\mathbf{u}_f{}\ast\gamma)^\wedge(t,x)$. Thus, by Lemma \ref{LemSmoothCurvCInftyMN}, we have $\mathbf{u}_f{}_\ast\gamma \in C^\infty(\R, \mathbb{V}_f)$ which implies $\mathbf{u}_f \in C^\infty(\mathbb{U}_f,\mathbb{V}_f)$. The proof of smoothness of $\mathbf{u}^{-1}_f$ is identical, since the criterion \ref{LemSmoothCurvCInftyMN} for smoothness of curves applies to $\mathbb{V}_f$ as well due to the choice of topology discussed above \eqref{EqDefLocTriv}.\\
To see that the bundle structure does not depend on the choice of the metrics on the base and the fibre as well as of the connection $\LC$, first note that the smooth structure on $C^\infty(M,S)$ does not depend on them. In fact, these objects only enter the construction of $G$ and the smooth structure itself is independent of $G$ by Theorem \ref{ThmCInftyMN}. To prove that the fibre bundle structures induced by two different connections $\LC,\LC'$ are isomorphic, note that both resulting trivializations are defined on the same neighbourhoods $\mathbb{U}_f$, $\mathbb{V}_f$. By \eqref{EqLocTrivPt}, the transition maps have the form $(\xi,\zeta) \mapsto (\xi, P'$\raisebox{-1pt}{$\scriptstyle {u_f^{-1}(\xi),f}$} $ \circ P$\raisebox{-1pt}{$\scriptstyle {f,u_f^{-1}(\xi)}$} $(\zeta))$ where $P'$ is the parallel transport associated to $\LC'$. They are bijective, fibrewise linear and smooth as composites of smooth maps. 
\end{Proof}

\begin{Rem} \label{RemChartsSmooth}
The charts on $C^\infty(M,S)$ (as defined in \eqref{EqCharts1} using the metric $G$) and the local trivializations from \eqref{EqDefLocTriv} are of course closely related. In fact, for $f \in C^\infty(M,N)$, $\sigma_0 \in \mathbb{U}_f$ and $s_0 := \mathbf{u}_f(\sigma_0)$, the chart $u_{\sigma_0} : U_{\sigma_0} \rightarrow V_{\sigma_0} \cong V_f \times \Gamma_c(f^\ast S)$ is simply given by $\mathbf{u}_f{}|_{U_{\sigma_0}} - (0,s_0)$. Since smoothness is a local property and the sets $U_{\sigma_0} \subset \mathbb{U}_f$, $V_{\sigma_0} \subset \mathbb{V}_f$ are open by definition of the locally affine topology, this provides an alternative proof of the smoothness of $\mathbf{u}_f$ and $\mathbf{u}_f^{-1}$. 
\end{Rem}

\begin{Rem}
The fibre bundle from Theorem \ref{ThmBundleStruct} is not a vector bundle in the sense of \cite{KM} chapter 29, even though the trivializations in \eqref{EqLocTrivPt} are fibrewise linear. In fact, the fibre $\Gamma(f^\ast S)$, equipped with the structure constructed above, is only a locally affine space (cf. Example \ref{ExLocAffine}). Moreover, it is interesting to note that the bijections \eqref{EqLocTrivPt} do not provide a vector bundle atlas if $M$ is not compact and the spaces $\Gamma(f^\ast S)$ are equipped with the convenient vector space structure described in chapter \ref{ChapPre}. In that case, $\mathbf{u}_f^{-1}$ fails to be smooth. This follows from the fact that smooth curves $\gamma \in C^\infty(\R,\Gamma(f^\ast S))$ bijectively correspond to smooth maps $\R \times M \rightarrow S$ s.t. $\gamma(t,x) \in S_{f(x)}$, without any assumptions on supports! Thus, the curve $(\mathbf{u}_f)^{-1}_\ast\gamma$ in general does not satisfy the condition in Lemma \ref{LemSmoothCurvCInftyMN} which was essential in the proof of Theorem \ref{ThmBundleStruct}.
\end{Rem}

\begin{Rem} \label{RemMichor}
It is possible to obtain the bundle structure on $C^\infty(M,S)$ in a different way, which was suggested to me by P. Michor on MathOverFlow. We will describe it only very briefly omitting technical details. The parallel transport $P$ induced by the connection $\nabla$ can be considered as a map
\begin{align} \label{EqExDefPT}
 P_{t_0} : C^\infty(\R,N) \times_{ev_{t_0},\pi} S \times \R &\rightarrow S, &
 P(\gamma, s_0, t) &:= s(t).
\end{align}
Here, $s : \R \rightarrow S$ is the unique solution of the differential equation $\nabla^{\gamma^\ast S}_t s = 0$ along the curve $\gamma : \R \rightarrow N$ satisfying the initial condition $s(t_0) = s_0$. By proving that it maps smooth curves to smooth curves, it can be shown that this map is in fact smooth, either in the sense of infinite-dimensional manifolds or in the (more general) sense of Fr\"{o}hlicher spaces (cf. \cite{KM} chapter 23). After a suitable modification of the chart $u_f$ from \eqref{EqCharts1}, one may consider the family of curves $\gamma_{f,g}(t) := u_f^{-1}(t \cdot u_f(g))$ which is smooth in $t \in \R$ and $g \in U_f$. It connects $f$ and $g$ through a family of geodesics and the induced map $\overline{\gamma}_{f,g} : M \rightarrow C^\infty(\R,N)$ sending each $M$ to this geodesic can be shown to be smooth, too. Composing the  parallel transport from \eqref{EqExDefPT} with these curves (or rather there inverses $\overline{\gamma}^{-1}$ defined by going backwards) provides us with the following alternative description of the local trivializations $\mathbf{u}_f$ from \eqref{EqDefLocTriv}:
\begin{align*}
\mathbb{U}_f \ni s &\mapsto (\pi \circ s, P_0 \circ (\overline{\gamma}^{-1}_{f,\pi_\ast(s)}{},s,1) \in \mathbb{V}_f
\end{align*}
Comparison with \eqref{EqLocTrivPt} shows, that the underlying geometric construction is the same as before. This approach does not rely on a Riemannian metric on the total space $S$ and it is now easier to change the manifold structure on $C^\infty(M,S)$ turning it into a vector bundle. Equipping $\Gamma(f^\ast S)$ (no compact supports !) with the convenient structure discussed after Definition \ref{DefSmoothMap}, we obtain a coarser topology on the spaces $\mathbb{V}_f$. Since chart transitions are still diffeomorphisms and fibrewise linear, $C^\infty(M,S) \rightarrow C^\infty(M,N)$ becomes a smooth vector bundle but the topology on its total space is again coarser than the one used in Theorem \ref{ThmBundleStruct}. We will not use this structure in the remaining part of this publication.       
\end{Rem}

\vspace{\baselineskip}

We will now use Theorem \ref{ThmBundleStruct} to equip $\qq{SC}^\infty(X,Y)(\Lambda_n)$ with the structure of an infinite dimensional manifold. Recall that, by definition, 
\begin{align*}
\pi^{n} : T^n &:= (TY \boxtimes \tL^{\nil}(\R^n \oplus E^\ast))_\0 \rightarrow \tilde{X} \times \tilde{Y}
\end{align*}
is the vector bundle over $\tilde{X} \times \tilde{Y}$ whose fibre at $(x,y)$ is $T^n_{x,y} = (T_y Y \otimes \tL^{\nil}(\R^n \oplus E^\ast_x))_\0$. Hence, $C^\infty(\tilde{X},T^n)$ carries a smooth structure, which is a smooth affine bundle over $C^\infty(\tilde{X},\tilde{X}\times\tilde{Y})$ by Theorem \ref{ThmBundleStruct}. However, comparing with \eqref{EqSetSmoothStruct}, $C^\infty(\tilde{X},T^n)$ is too large because its base is $C^\infty(\tilde{X},\tilde{X}\times \tilde{Y})$ rather than $C^\infty(\tilde{X},\tilde{Y})$ and it contains sections along arbitrary maps into the product. To restrict it, consider the projections
\begin{align*}
\xymatrix{ \pr^n_{\tilde{X}/\tilde{Y}} :  T^n \ar[r]^-{\pi^{n}} & \tilde{X} \times \tilde{Y} \ar[r]^-{\pr_{\tilde{X}/\tilde{Y}}} & \tilde{X} / \tilde{Y}} 
\end{align*}
By \cite{KM}, Theorem 42.20, $\Gamma(\pr_{\tilde{X}}) = \{(\mathrm{Id}_{\tilde{X}}, g) \mid g \in C^\infty({\tilde{X}},{\tilde{Y}})\}$ is a splitting submanifold of $C^\infty({\tilde{X}},{\tilde{X}} \times {\tilde{Y}})$ and the same is true for $\Gamma(\pr^n_{\scriptscriptstyle \tilde{X}}) \subset C^\infty({\tilde{X}}, T^n)$. Clearly, $\Gamma(\pr_{\scriptscriptstyle\tilde{X}}^n) \rightarrow \Gamma(\pr_{\scriptscriptstyle\tilde{X}})$ is the restriction of the bundle $C^\infty(\tilde{X}, T^n) \rightarrow C^\infty({\tilde{X}}, {\tilde{X}} \nolinebreak \times \nolinebreak {\tilde{Y}})$ to $\Gamma(\pr_{\tilde{X}})$. In particular, it inherits the structure of an affine bundle and its local trivializations are given by restricting those of $C^\infty({\tilde{X}},T^n)$ in \eqref{EqDefLocTriv}, \eqref{EqLocTrivPt}. Moreover, a point in this space is an element of $\Gamma((f,\mathrm{Id}_{\tilde{X}})^\ast T^n) \cong \Gamma(f^\ast TY \otimes \tL^{\nil}(\R^n \oplus E^\ast))_\0)$ for some $f \in C^\infty({\tilde{X}},{\tilde{Y}})$ so by \eqref{EqSetSmoothStruct}, we have a bijection
\begin{align} \label{EqIdentSpace}
\qq{\SC}^\infty(X,Y)(\Lambda_n) &\cong \Gamma(\pr^{n}_{\tilde{X}}) \subset C^\infty({\tilde{X}}, T^n)
\end{align}
Finally, $\Gamma(\pr_{\tilde{X}})$ is diffeomorphic to $C^\infty({\tilde{X}},{\tilde{Y}})$ since the map 
$I(f) := (\mathrm{Id}_{\tilde{X}},f)$ is in fact a diffeomorphism. This is easily verified by checking that $I$ and $I^{-1}$ map smooth curves to smooth curves, which follows from Lemma \ref{LemSmoothCurvCInftyMN} and the characterisation of smooth maps into submanifolds given in \cite{KM} (27.11).\\

To check functoriality, let $\rho \in \Hom_\Gr(\Lambda_n,\Lambda_m)$. Clearly, $\rho(\Lambda_n^{\nil})\subset \Lambda_m^{\nil}$. By Definition \ref{DefSCInfty}, we have $\qq{SC}^\infty(X,Y)(\rho)\Phi = \Phi \circ (\rho^\ast \times \mathrm{Id}_{\tilde{X}})$, which on the level of superfunctions translates to $(\qq{SC}^\infty(X,Y)(\rho)\Phi)^\ast = (\rho \otimes \mathrm{Id}_{\O}) \circ \Phi^\ast : \Lambda_n \otimes \Gamma(\tL E^\ast) \rightarrow \Lambda_m \otimes \Gamma(\tL E^\ast)$. Under the identifications \eqref{EqIdentSpace}, this map becomes
\begin{align} \label{EqRhoAction}
\qq{SC}^\infty(X,Y)(\rho) : \Gamma(\pr^{n}_{\tilde{X}}) &\rightarrow \Gamma(\pr^{m}_{\tilde{X}}), &
\sigma &\mapsto (\mathrm{Id}_{TY} \otimes \rho \otimes \mathrm{Id}_{E^\ast})_\ast \sigma
\end{align}
Thus, it is the restriction of the smooth map $(\mathrm{Id}_{TY} \otimes \rho \otimes \mathrm{Id}_{E^\ast})_\ast : C^\infty({\tilde{X}}, T^n) \rightarrow C^\infty({\tilde{X}}, T^m)$ to the submanifold $\Gamma(\pr^{n}_{\tilde{X}})$, which is well-defined and again smooth. Summarizing, we have shown:

\begin{Prop} \label{PropSCInfMan}
$\qq{SC}^\infty(X,Y)$ becomes a functor from the category of Grassmann algebras into the category of manifolds (in the sense of chapter \ref{ChapPre}) by \eqref{EqIdentSpace} and \eqref{EqRhoAction}.  
\end{Prop}

\begin{Rem}
In case the supermanifold $X$ is compact, it follows from Remark \ref{RemFrechetMf} that $\qq{SC}^\infty(X,Y)$ is actually a functor $\Gr \rightarrow \Man_F$ where $\Man_F$ denotes the category of Fr\'{e}chet-manifolds.
\end{Rem}

In order to establish that the functor $\qq{SC}^\infty(X,Y)$ carries the structure of a convenient supermanifold, it remains to provide a supersmooth atlas. By \eqref{EqSplitModelSpace}, the model space for  $C^\infty(\tilde{X},T^n)$ at $\sigma'$ with underlying map $(g,f) := \pi^{n}\circ\sigma' \in C^\infty(\tilde{X}, \tilde{X} \times \tilde{Y})$ is given by  
 \begin{align*}
 \Gamma_c(\tilde{X}, (g,f)^\ast T(\tilde{X} \times \tilde{Y})) \oplus \Gamma_c(\tilde{X},(f^\ast TY \boxtimes g^\ast\tL^{\nil}(\R^n \oplus E^\ast))_\0).
 \end{align*}
Restricting to $\Gamma(\pr^n_{\tilde{X}})$, the first summand reduces to $\Gamma_c(f^\ast T\tilde{Y})$ and we have $g = \mathrm{Id}_{\tilde{X}}$. Hence, the model space for the manifold $\qq{SC}^\infty(X,Y)(\Lambda_n)$ near $\sigma$ with underlying map $f := \pr^n_{\tilde{Y}}\circ\sigma \in C^\infty(\tilde{X},\tilde{Y})$ is given by
\begin{align} \label{EqModelSpacePr}
\Gamma_c(f^\ast T\tilde{Y}) \oplus \Gamma_c((f^\ast TY \otimes \tL^{\nil}(\R^n \oplus E^\ast))_\0)
&\cong (\Lambda_n \otimes \Gamma_c(f^\ast TY \otimes \tL E^\ast))_\0 
\end{align}
We see that it is given by the $\Lambda_n$-points of the superrepresentable module $\qq{\Gamma_c(f^\ast TY \otimes \tL E)}$ (cf. \ref{EqRepSMod}), which depends on the basepoint $f \in C^\infty(\tilde{X},\tilde{Y})$ . Note that the contribution $\Gamma_c(f^\ast T\tilde{Y})$ arising from the base $C^\infty(\tilde{X},\tilde{Y})$ is crucial to obtain superrepresentability. To define charts, we define subfunctors $\qq{\mathbb{U}}_f \subset \qq{SC}^\infty$, $\qq{\mathbb{V}}_f \subset \qq{\Gamma(f^\ast TY \otimes \tL E^\ast)}$ at $f \in C^\infty(\tilde{X},\tilde{Y})$ by
\begin{align} \label{EqDefSuperChart}
\qq{\mathbb{U}}_f(\Lambda_n) &:=  (\pr^n_{\tilde{X}})^{-1}(U_f) \\
\qq{\mathbb{V}}_f(\Lambda_n) &:=  V_f \times \Gamma((f^\ast TY \otimes \tL^{\nil}(\R^n \oplus E^\ast))_\0) \notag
\end{align}

\begin{Rem} \label{RemOpenSF}
Both subfunctors are restrictions (cf. \eqref{EqRestriction}) to open sets of $C^\infty(\tilde{X},\tilde{Y})$ and $\Gamma(f^\ast T \tilde{Y})$ respectively, hence open. For $\qq{\mathbb{U}}_f$, this is obvious and for $\qq{\mathbb{V}}_f$, it follows from the fact that $\Gamma(f^\ast T\tilde{Y})$ is equipped with the locally affine topology and hence $\Gamma_c(f^\ast T\tilde{Y})$ as well as $V_f$ are open subsets of $\Gamma(f^\ast T\tilde{Y})$. In particular, $\qq{\mathbb{V}}_f$ is a locally affine convenient superdomain in the sense of Definition \ref{DefLocAffDom}.
\end{Rem}

\begin{Rem}
Similar to the argument in Example \ref{ExCounterChart}, the functors $\qq{\mathbb{V}}_f$ cannot be realized as convenient superdomains. Moreover, we can not simply equip $\Gamma(\pr^n_{\tilde{X}})$ with the vector bundle structure described Remark \ref{RemMichor}, because the different topologies on base and fibres lead to model spaces, which are not superrepresentable. However, one may use a vector space topology on the summand $\Gamma(f^\ast F \otimes \Lambda^{\odd}(\R^n \oplus E^\ast))$ of the fibre. Applications (e.g. in field theory) may indicate, which structure is most useful.  
\end{Rem}

The definition of the charts itself is motivated by results from chapter \ref{ChapFinDim}. Recall that the situation discussed there corresponds to $X = \{0\}$. Thus, a reasonable definition of charts should be obtained by extending \eqref{EqMorPnPoints} to arbitrary $X$ of dimension $p|q$. Let $f \in C^\infty(\tilde{X},\tilde{Y})$ and $x \in \tilde{X}$. Using the open sets $U,V$ from \eqref{EqCharts0}, we obtain open sets $V_{f(x)} := V \cap T_{f(x)}\tilde{Y} \subset T_{f(x)}\tilde{Y}$ and $U_{f(x)} := \{y \in \tilde{Y} \mid (f(x),y) \in U \} \subset \tilde{Y}$. \eqref{EqLocTrivPt} suggests to consider
\begin{align} 
\varphi_{f(x)} : U_{f(x)} &\rightarrow V_{f(x)},  & \varphi_{f(x)}(y) &:= \exp_{f(x)}^{-1}(y) \label{EqPointwiseDefChart} \\
\phi_{f(x)}: F|_{U_{f(x)}} &\rightarrow V_{f(x)} \times F_{f(x)}, & \phi_{f(x)}(\zeta) &:= (\varphi_{f(x)}(\pi(\zeta)), P_{\pi(\zeta),f(x)}(\zeta) ) \notag
\end{align}
Clearly, for each $x$ and $f$, $\phi_{f(x)}$ is a vector bundle isomorphism over $\varphi_{f(x)}$ and we may apply the construction discussed in chapter \ref{ChapFinDim} to obtain a local chart $\Phi_{f(x)}$ of $Y$ near $f(x)$. Putting $r:=\lfloor n/2 \rfloor$, Proposition \ref{PropFindimNPoints} shows that $\Phi_{f(x)}$ maps a $\Lambda_n$-point $\mu \in (T(Y|_{U_{f(x)}}) \otimes \Lambda_n^{\nil})_\0$ of $Y$ to 
\begin{align*}
\Phi_{f(x)}{_\ast}(\mu) &= S^r(\varphi_{f(x)})(\mu) + S^r(\phi_{f(x)})(\mu) \in  V_{f(x)} \times ((T_{f(x)}\tilde{Y} \oplus F_{f(x)})\otimes \Lambda_n^{\nil})_\0
\end{align*}
Note that $V_{f(x)} \times (T_{f(x)}\tilde{Y} \oplus F_{f(x)})$ is simply the tangent bundle of the supermanifold defined by the trivial Batchelor bundle $V_{f(x)} \times F_{f(x)} \rightarrow V_{f(x)}$.\\ 

The charts on the infinite-dimensional spaces can now defined by using this formula for each $x \in \tilde{X}$ and extending the arguments of the maps $S^r(\varphi_{f(x)})$, $S^r(\phi_{f(x)})$ from $(T(Y|_{U_{f(x)}}) \otimes \Lambda_n^{\nil})_\0$ to $(T(Y|_{U_{f(x)}}) \otimes \tL^{\nil}(\R^n \oplus E_x))_\0$: 
\begin{align} \label{EqDefSAtlas}
\Phi_{f,n} : \qq{\mathbb{U}}_f(\Lambda_n) &\rightarrow \qq{\mathbb{V}}_f(\Lambda_n) &
\Phi_{f,n}(\sigma) &:= (x \mapsto S^r(\varphi_{f(x)})(\sigma(x))  + S^r(\phi_{f(x)})(\sigma(x)) )
\end{align}
where now $r := \lfloor (n+\mathrm{rk}(E))/2\rfloor$, $\sigma(x) \in (T_{g(x)}Y \otimes \tL^{\nil}(\R^n \oplus E^\ast_x))_\0$ and $g := \pr^n_{\tilde{Y}} \sigma \in C^\infty(\tilde{X},\tilde{Y})$ is the smooth map underlying $\sigma$. To prove smoothness and functoriality properties, recall from \eqref{EqPhiSimpTotal} that 
\begin{align}
S^r(\varphi_{f(x)})(\sigma(x)) &= \sum_{I \in \N_0^{p,r}} \tfrac{1}{I!} D_I \varphi_{f(x)}|_{\tilde{\sigma}(x)} \otimes \sigma(x)_\2^I \label{EqChartLocal1} \\
S^r(\phi_{f(x)})(\sigma(x)) &=  \sum_{I \in \N_0^{p,r}} \sum_{\alpha = p+1}^{p+q} \tfrac{1}{I!} D_I \phi_{f(x)}|_{\tilde{\sigma}(x)}(v_\alpha) \sigma(x)_\2^I \sigma(x)_\1^\alpha \label{EqChartLocal2}
\end{align}
From this, we can conclude

\begin{Lemma} \label{LemPropCharts}
The maps $\{\Phi_{f,n}\}_{n \in \N_0}$ from \eqref{EqDefSAtlas} define a natural isomorphism $\Phi_f : \qq{\mathbb{U}}_f \rightarrow \qq{\mathbb{V}}_f$ in $\MGr$.
\end{Lemma}
\begin{Proof}{}
To see that $\Phi_{f,n}$ is bijective for all $n \in \N_0$, observe that the maps \eqref{EqPointwiseDefChart} are bijective for each $x \in \tilde{X}$. Since $S^r$ is functorial by Lemma \ref{LemmaFuncSK}, the inverse to $\Phi_{f,n}$ is then simply given by
\begin{align*}
\Phi_{f,n}^{-1}: \qq{\mathbb{V}}_f(\Lambda_n) &\rightarrow \qq{\mathbb{U}}_f(\Lambda_n) &
\Phi_{f,n}^{-1}(\tau) &:= (x \mapsto S^r(\varphi^{-1}_{f(x)})(\tau(x))  + S^r(\phi^{-1}_{f(x)})(\tau(x)) ) 
\end{align*}
Next, let $\rho \in \Hom_\Gr(\Lambda_n,\Lambda_m)$.  From \eqref{EqRhoAction} and the fact that $\rho$ is an algebra homomorphism, we deduce 
\begin{align*}
\Phi_f(\qq{\mathbb{U}}_f(\rho)\sigma)(x) &= \sum_{I \in \N_0^{p,r}} \tfrac{1}{I!} \Bigl[ 
                                         D_I \varphi_{f(x)}|_{\tilde{\sigma}(x)} \otimes \rho(\sigma(x)_\2^I) + \!\!
                                         \sum_{\alpha = p+1}^{p+q} (D_I \phi_{f(x)}|_{\tilde{\sigma}(x)})(v_\alpha) \rho(\sigma(x)_\2^I \sigma(x)_\1^\alpha) \Bigr]
\end{align*}
where $\rho$ acts on $\tL(\R^n \oplus E^\ast) \cong \Lambda_n \otimes \tL E^\ast$ by $\rho \otimes \mathrm{Id}_{\bigwedge E^\ast}$. By \eqref{EqRepSMod}, the right hand side equals $\qq{\mathbb{V}}(\rho)\Phi_f(\sigma)(x)$ which proves that $\{\Phi_{f,n}\}_{n \in \N_0}$ is a natural transformation.\\
It remains to show that $\Phi_{f,n}$ and $\Phi_{f,n}^{-1}$ are smooth. We essentially use the argument which was already used in the proof of Theorem \ref{ThmBundleStruct}: For $\gamma \in C^\infty(\R,\qq{\mathbb{U}}(\Lambda_n))$, \eqref{EqChartLocal1} and \eqref{EqChartLocal2} imply that $(\Phi_{f,n}{}_\ast\sigma)^\wedge(t,x)$ is smooth in $t$ and $x$. Moreover, the dependence of $\varphi_{f(x)}$, $\phi_{f(x)}$ on $x \in \tilde{X}$ given in \eqref{EqPointwiseDefChart} implies, that the condition stated in Lemma \eqref{LemSmoothCurvCInftyMN} is preserved by $\Phi_{f,n}{}_\ast$. Thus, $\Phi_{f,n}$ map is smooth for all $n \in \N_0$. The same argument can be applied to $\Phi_{f,n}^{-1}$.  
\end{Proof}

We will now state and prove the main result of this article.
 
\begin{Thm} \label{ThmSMfSCInfty}
$\qq{SC}^\infty(X,Y)$ is a smooth, infinite-dimensional convenient supermanifold. A supersmooth atlas is provided by $\mathcal{A} := \{\Phi_f^{-1} : \qq{\mathbb{V}}_f \rightarrow \qq{\mathbb{U}}_f\}_{f \in C^\infty(\tilde{X},\tilde{Y})}$ as defined in \eqref{EqDefSAtlas}.
\end{Thm}

\begin{Proof}{}
First of all, it follows from \eqref{EqDefSuperChart} that $\{\qq{\mathbb{U}}_f\}$\raisebox{-2pt}{$\scriptscriptstyle {f\in C^\infty(\tilde{X},\tilde{Y})}$} covers $\qq{SC}^\infty(X,Y)$ because the $U_f$ cover $C^\infty(\tilde{X},\tilde{Y})$. Thus, $\mathcal{A}$ is an open covering by locally affine superdomains by Remark \ref{RemOpenSF} and Lemma \ref{LemPropCharts}.\\
It remains to check the second condition in Definition \ref{DefSMfFunc}. By Remark \ref{RemSuperLin}, it is equivalent to check that differentials of the chart transitions 
\begin{align*}
\Psi := \Phi_{g,n} \circ \Phi^{-1}_{f,n} : \qq{\mathbb{V}}_f|_{V_{fg}}(\Lambda_n) &\overset{\sim}{\longrightarrow} \qq{\mathbb{V}}_g|_{V_{fg}}(\Lambda_n) 
\end{align*}
are $\Lambda_{n}^{\mathrm{ev}}$-linear for all $n \in \N_0$ and $f,g \in C^\infty(\tilde{X},\tilde{Y})$. Here, we set $V_{f,g} := V_f \cap V_g \subset C^\infty(\tilde{X},\tilde{Y})$. We first observe that, by Theorem \ref{ThmCInftyMN}, it is enough to verify this linearity for arbitrary but fixed $x \in \tilde{X}$.  To simplify notation, let $\varphi_x := \varphi_{g(x)}\circ\varphi_{f(x)}^{-1}|_{V_{fg,x}}$  and $\phi_x := \phi_{g(x)}\circ\phi_{f(x)}^{-1}|_{V_{fg,x}}$ with the notation from \eqref{EqPointwiseDefChart}. By Remark \ref{RemFunkSK} and Lemma \ref{LemmaFuncSK}, the chart transition at $x$ is given by
\begin{align} \label{EqExplicitCoordChange}
\Psi_x:  V_{f(x)} \times (T_{f(x)}Y \otimes \tL^{\nil}(\R^n \oplus E_x^\ast))_\0 
      &\rightarrow  V_{g(x)} \times (T_{g(x)}Y \otimes \tL^{\nil}(\R^n \oplus E_x^\ast))_\0 \\
\Psi_x(\kappa) &= S^r(\varphi_x)(\kappa) + S^r(\phi_x)(\kappa)
\end{align}
Note that both, the domain and the range of $\Psi_x$, carry the flat geometry induced by $T_{f(x)}Y$ and $T_{g(x)}Y$, respectively. In particular, we may use ordinary finite-dimensional calculus during the subsequent calculations (cf. Remark \ref{RemFlatTaylor}). For $\kappa \in V_{f(x)} \times (T_{f(x)}Y \otimes \tL^{\nil}(\R^n \oplus E_x^\ast))_\0$ and $\tau \in (T_{f(x)}Y \otimes \tL(\R^n \oplus E_x^\ast))_\0$, we thus obtain
\begin{align} \label{EqCompDifferential}
d S^r(\varphi_x)|_{\kappa}(\tau)
&= \sum_{I \in \N_0^{p,r}} \tfrac{1}{I!} D_{\tilde{\tau}}D_I\varphi_x|_{\tilde{\kappa}} \kappa_\2^I 
 + \tfrac{1}{I!} D_I\varphi_x|_{\tilde{\kappa}} \sum_{l=1}^p i_l \tau_\2^l \kappa_\2^{I-e_l}  \\
d S^r(\phi_x)|_\kappa(\tau) 
&= \sum_{I \in \N_0^{p,r}} \sum_{\alpha=p+1}^{p+q} \tfrac{1}{I!} D_{\tilde{\tau}}D_I\phi_x|_{\tilde{\kappa}}(v_\alpha) \kappa_\2^I \kappa_\1^\alpha 
 + \tfrac{1}{I!} D_I\phi_x|_{\tilde{\kappa}}(v_\alpha) 
   (\sum_{l=1}^p i_l \tau_\2^l \kappa_\2^{I-e_l} \kappa_\1^\alpha + \kappa_\2^I \tau_\1^\alpha) \notag
\end{align}
The decomposition $\tau = \tau_\0 + \tau_\1 \in T_{f(x)}\tilde{Y} \otimes \tL(\R^n \oplus E_x^\ast)_\0 \oplus F_{f(x)} \otimes \tL(\R^n \oplus E_x^\ast)_\1$ is clearly preserved under multiplication by $\lambda \in \Lambda_{n}^{\ev}$ which allows us to discuss the cases $\tau = \tau_\0$, $\tau = \tau_\1$ separately. Since $\tau_\0$ furthermore splits as $\tau_\0 = \tilde{\tau} + \tau_\2$ and assuming $\tau = \tau_\1$ implies $\tilde{\tau} = 0$, it follows immediately from \eqref{EqCompDifferential} that $d\Psi_x|_\kappa(\lambda \tau_\1) = \lambda d\Psi_x|_\kappa(\tau_\1)$.\\
To discuss the case $\tau = \tau_\0$, we first observe that smoothness of $\Psi_x$ implies $\R$-linearity of $d\Psi_x$. Thus, we may assume $\lambda \in \Lambda_n^{ev \geq 2}$ which implies $\widetilde{\lambda \tau} = 0$ and $(\lambda\tau)_\2 = \lambda\tau = \lambda(\tilde{\tau}+\tau_\2)$. Inserting it into $dS^r(\varphi_x)$ yields
\begin{align*}
 dS^r(\varphi_x)|_{\kappa}(\lambda\tau) &= \sum_{I\in \N_0^{p,r}} \tfrac{1}{I!} D_I\varphi_x|_{\tilde{\kappa}} \sum_{l=1}^p i_l \lambda(\tilde{\tau}+\tau_\2)^l \kappa_\2^{I-e_l}
\end{align*}
Comparing this expression with \eqref{EqCompDifferential}, it follows that $\lambda dS^r(\varphi_x)|_{\kappa}(\tau) = dS^r(\varphi_x)|_\kappa(\lambda\tau)$ is equivalent to
\begin{align} \label{EqEquivdS}
\sum_{I \in \N_0^{p,r}} \tfrac{1}{I!} D_I\varphi_x|_{\tilde{\kappa}} \sum_{l=1}^p i_l \lambda \tilde{\tau}^l \kappa_\2^{I-e_l} 
&= \lambda \sum_{I \in \N_0^{p,r}} \tfrac{1}{I!} D_{\tilde{\tau}}D_I\varphi_x|_{\tilde{\kappa}} \kappa_\2^I
\end{align}
Next, we rearrange the double sum on the left hand side: 
\begin{align} \label{EqCompLambdaLin}
\sum_{l=1}^p \sum_{I \in \N_0^{p,r}} \tfrac{i_l}{I!} D_I\varphi_x|_{\tilde{\kappa}} \lambda \tilde{\tau}^l \kappa_\2^{I-e_l} &=
\lambda \sum_{l=1}^p \sum_{J \in \N_0^{p,r-1}} \tfrac{1}{J!} \tilde{\tau}^l D_{v_l} D_J\varphi_x|_{\tilde{\kappa}} \kappa_\2^J =
\lambda\sum_{J \in \N_0^{p,r}} \tfrac{1}{J!} D_{\tilde{\tau}} D_J\varphi_x|_{\tilde{\kappa}} \kappa_\2^J
\end{align}
Here, we performed the following steps: For each fixed $l=1,\ldots,p$, we can relabel the multi-index by setting $J := I- e_l$ because summands with $i_l = 0$ do not contribute. Observing that $i_l / I! = 1/J!$ yields the first equality. The second equality follows from the fact, that the Grassmann degree of $\lambda$ is $\geq 2$, hence $\lambda \kappa_\2^J = 0$ for $|J| = r = \lfloor (n + rk(E))/2 \rfloor$. Hence, we have proven \eqref{EqEquivdS}, i.e. $\Lambda_{n}^{\ev}$-linearity of $dS^r(\varphi_x)$.\\ 
Finally, the $\Lambda_n^{\mathrm{ev}}$-linearity of $dS^r(\phi_x)$ is obtained using a completely analogous argument. This follows from \eqref{EqCompDifferential}, because for $\tau_\1 = 0$, the expressions for $dS^r(\varphi_x)$ and $dS^r(\phi_x)$ have the same form. Summarizing, we have shown $d\Psi_x|_\kappa(\lambda\tau) = \lambda d\Psi_x|_\kappa(\tau)$ which concludes the proof.
\end{Proof}

\begin{Rem}
The argument below equation \eqref{EqCompLambdaLin} shows that it is necessary to use a suitable jet-extension of $\varphi_{f(x)}$ and $\phi_{f(x)}$ for the construction of the supercharts. Taking into account only their $0$-th order jets would result in a covering by morphisms in $\MGr$ which do not form a supersmooth atlas, even in the case when $\tilde{X}$ is compact.
\end{Rem}

\vspace{\baselineskip}

In the course of the discussion of $\qq{\SC}^\infty(X,Y)$, we used non-canonical isomorphisms
\begin{align} \label{EqChoiceIso}
\Phi_E : (\tilde{X}, \Gamma(\tL E^\ast)) &\overset{\sim}{\longrightarrow} X &
\Phi_F : (\tilde{Y}, \Gamma(\tL F^\ast)) &\overset{\sim}{\longrightarrow} Y 
\end{align}
of supermanifolds, where $E$ and $F$ are Batchelor bundles for $X$ and $Y$, respectively. Moreover, we chose Riemannian metrics on $\tilde{X}$, $\tilde{Y}$ and connections on $E$, $F$ to define the atlases for $C^\infty(\tilde{X}, T^n)$. Even though the resulting smooth structure is independent of all these choices by Theorem \ref{ThmBundleStruct}, the supersmooth structure may still depend on the connections on $F$ and $T\tilde{Y}$, because these enter the definition of the charts in \eqref{EqDefSAtlas} through $S^r$. In contrast, there is no additional dependence on the metric on $\tilde{X}$ and the connection on $E$. This is consistent with the observation that only the former connections implicitly enter the description of the model space for $SC^\infty(X,Y)(\Lambda_n) \cong \Gamma(\pr^{n}_{\tilde{X}})$ in \eqref{EqModelSpacePr}. In fact, the contribution $\Gamma_c(f^\ast T\tilde{Y})$ in \eqref{EqModelSpacePr} results from the identification $T\tilde{Y} \cong \mathcal{H}(TY)$ which depends on the connection on $F$ as well as the Levi-Civita connection on $T\tilde{Y}$.\\ 

The following result shows that also the supersmooth structure on $\qq{SC}^\infty(X,Y)$ is independent of these choices:

\begin{Prop}
The supersmooth structures on $\qq{\SC}^\infty(X,Y)$, induced by different identifications $\Phi_E$, $\Phi_F$ and different connections on $F$ are superdiffeomorphic, i.e. isomorphic in $\SMan$.
\end{Prop}

\begin{Proof}{}
Let $\Phi_{E_i}$ and $\Phi_{F_i}$ ($i=1,2$) two different choices for the isomorphisms in \eqref{EqChoiceIso} and $\LC^i$ connections on $T\tilde{Y} \oplus F_i \cong TY$. We clearly have $E_1  \cong E_2$ and $F_1 \cong F_2$. Denote the two resulting supermanifolds by $\qq{\SC}(X,Y)^{(i)}$ and the corresponding manifolds of sections from \eqref{EqIdentSpace} by $\Gamma(\mathrm{pr}_{\tilde{X}}^{\scriptscriptstyle k,(i)})$. Abbreviating $\Psi_{E_1E_2} := (\Phi_{E_1})^{-1}\Phi_{E_2}$ and $\Psi_{F_2F_1} := (\Phi_{F_2})^{-1}\Phi_{F_1}$, we get an induced isomorphism on the level of function algebras for each $n \in \N_0$: 
\begin{align} \label{EqChangeBat}
\Hom(\Gamma(\tL F_1^\ast), \Lambda_n \otimes \Gamma(\tL E_1^\ast)) &\cong \Hom(\Gamma(\tL F_2^\ast), \Lambda_n \otimes \Gamma(\tL E_2^\ast)) \\
\hat{\sigma} &\mapsto (\mathrm{Id}_{\Lambda_n} \otimes \Psi^\ast_{E_1E_2}) \circ \hat{\sigma} \circ \Psi^\ast_{F_2F_1} \notag
\end{align}
We will discuss the dependence on $\Phi_F$ and $\Phi_E$ separately. First let $\Phi_{E_1} = \Phi_{E_2}$ and $\psi \in C^\infty(\tilde{Y},\tilde{Y})$ be the diffeomorphism underlying $\Psi_{F_1F_2}$. Using the identification \eqref{EqIdentSpace}, Proposition \ref{PropFindimNPoints} and \eqref{EqPhiAstGeneral} yield the map on $\Lambda_n$-points induced by \eqref{EqChangeBat}:  
\begin{align*}
\Psi_n: \Gamma(\pr^{n,(1)}_{\tilde{X}}) &\cong \qq{\SC}(X,Y)^{(1)}(\Lambda_n) \rightarrow \Gamma(\pr^{n,(2)}_{\tilde{X}}) \cong \qq{\SC}(X,Y)^{(2)}(\Lambda_n) \\  
\Psi_n(\sigma) &:= \bigl( x \mapsto S^r(\symb \Psi_{F_2F_1}) \circ_r S^r(\psi)(\sigma(x))+ S^r(\symb \Psi_{F_2F_1})(\sigma(x)) \bigr) 
\end{align*}
The arguments given in the proof of Lemma \ref{LemPropCharts} establish that $\{\Psi_n\}_{n \in \N_0}$ defines an isomorphism in $\MGr$ and it remains to check its supersmoothness (cf. Remark \ref{RemSMfMorphism}). As in the proof of Theorem \ref{ThmSMfSCInfty}, it is sufficient to check the $\Lambda_n^{\ev}$-linearity of the differential of the coordinate representatives of $\Psi_n$ at each fixed $x \in \tilde{X}$. Using \eqref{EqPhiAstSpecial}, it is not difficult to adapt the computations following \eqref{EqCompDifferential} to this more general situation, we omit the details\footnote{Alternatively, one may use the results from \cite{BCF1}, Theorem 4.2, to conclude that the map induced by the morphism of supermanifolds $\Psi_{F_1F_2}$ has the desired properties.}. Hence, $\{\Psi_n\}_{n \in \N_0}$ provides a superdiffeomorphism $\qq{\SC}(X,Y)^{(1)} \overset{\sim}{\longrightarrow} \qq{\SC}(X,Y)^{(2)}$. Note that we have in particular shown that the supersmooth structure does not depend on the choice of the connection on $F$.\\
In a second step, assume $\Phi_{F_1} = \Phi_{F_2}$ and let $\psi \in C^\infty(\tilde{X},\tilde{X})$ now denote the diffeomorphism underlying $\Psi_{E_1E_2}$. According to \eqref{EqSplitSymbMor}, we may split $\Psi_{E_1E_2} = \symb \Psi_{E_1E_2}^\ast \circ \psi^\ast \circ S^k$ putting $k := \lfloor \mathrm{rk}(E)/2 \rfloor$. Recall that $\qq{SC}^\infty(X,Y)^{(i)}(\Lambda_n)$ is a bundle over $C^\infty(\tilde{X},\tilde{Y})$ (cf. \eqref{EqIdentSpace}). The map from \eqref{EqChangeBat} translates to the map $\Psi_n$ given on the fibre at $f \in C^\infty(\tilde{X},\tilde{Y})$ of this bundle by the following diagram:
\begin{align} \label{EqDecompPsiE}
\xymatrix{ \Gamma((f^\ast TY \otimes \tL^{\nil}(\R^n \oplus E_1^\ast))_\0)  \ar[r]^-{\psi^\ast \circ S^{k}} \ar[rd]_{\Psi_n}
                    & \Gamma(\psi^\ast(f^\ast TY \otimes \Sym^{\leq k}(\R^n \oplus T^\ast X)^{\mathrm{nil}})_\0) \ar[d]^-{(\symb\Psi_{E_1E_2}^\ast)_\ast} \\
                    &  \Gamma(\psi^\ast (f^\ast TY \otimes \tL^{\nil}(\R^n \oplus E_2^\ast))_\0) } 
\end{align}
This representation extends to all of $\qq{SC}^\infty(X,Y)^{(i)}(\Lambda_n)$. Each $\Psi_n$ is smooth by Lemma \ref{LemSmoothCurvCInftyMN} since the maps in \eqref{EqDecompPsiE} preserve the conditions given in that lemma. By \eqref{EqChangeBat}, we have $\qq{SC}^\infty(X,Y)^{(2)}(\rho)\circ \Psi_n = \Psi_m \circ \qq{SC}^\infty(X,Y)^{(1)}(\rho)$ for $\rho \in \Hom_\Gr(\Lambda_n,\Lambda_m)$. Hence, $\{\Psi_n\}_{n \in \N_0}$ defines an isomorphism in $\MGr$ whose inverse is induced by $\Psi_{E_2E_1}$. To prove supersmoothness, let $\Psi$ denote a coordinate transformation from $f$  to $g$ (cf. proof of Theorem \ref{ThmSMfSCInfty}). For $\tau \in \qq{\mathbb{V}}_f(\Lambda_n)$, let $\hat{\tau}$ denote the associated homomorphism of superalgebras. Its image under \eqref{EqChangeBat} in coordinates around $g$, evaluated at some $x \in \tilde{X}$, is then given by  
\begin{align} \label{EqMapPsiEE}
[\mathrm{ev}_x \circ (\mathrm{Id_{\Lambda_n}} \otimes \Psi_{E_1E_2}^\ast) \circ \hat{\tau}] \circ \Psi^\ast 
\in \Hom_{\mathrm{SAlg}}(\Gamma(\tL(V_{g(x)} \times F^\ast_{g(x)})), \tL(\R^n \oplus E^\ast_{2,x}))
\end{align}
Let $\Psi_{E_1E_2}(\tau)(x)$ be the element in $V_{f(x)} \times (T_{f(x)}Y \otimes \tL^{\nil} (\R^n \oplus E^\ast_{2,x}))$ defined by the homomorphism $[\mathrm{ev}_x \circ (\mathrm{Id}_{\Lambda_n} \otimes \Psi_{E_1E_2}^\ast) \circ \hat{\tau}]$ as described above \eqref{EqExplicitNPoint}. Similar to \eqref{EqExplicitCoordChange}, the map on $\Lambda_n$-points corresponding to \eqref{EqMapPsiEE} reads
\begin{align*}
\qq{\mathbb{V}}_f(\Lambda_n) &\longrightarrow \qq{\mathbb{V}}_g(\Lambda_n), & \tau &\mapsto \bigl(x \mapsto S^r(\varphi_x)(\Psi_{E_1E_2}(\tau)(x)) + S^r(\phi_x)(\Psi_{E_1E_2}(\tau)(x))\bigr)  
\end{align*}
Since $\Psi_{E_1E_2}^\ast$ acts trivially on $\Lambda_n$ by \eqref{EqChangeBat}, it is not hard to verify that $\Psi_{E_1E_2}(\kappa + \lambda\tau) = \Psi_{E_1E_2}(\kappa) + \lambda\Psi_{E_1E_2}(\tau)$ for all $\lambda \in \Lambda_n^{\ev}$ and $\kappa \in \qq{\mathbb{V}}_f(\Lambda_n), \tau \in \qq{\Gamma_c(f^\ast TY \otimes \tL E)}(\Lambda_n)$. Thus, we may again adapt the calculations following \eqref{EqCompDifferential} by replacing $\tau$ by $\Psi_{E_1E_2}(\tau)$ to obtain supersmoothness of $\{\Psi_n\}$. Hence, we have shown that different identifications $\Phi_{E_i}$ lead to isomorphic supermanifolds $\qq{SC}^\infty(X,Y)^{(i)}$.
\end{Proof}

This result completes the construction of the supermanifold $\qq{SC}^\infty(X,Y)$. Further details (e.g. tangent structures, supersmoothness of canonical maps as compositions etc.) as well as the construction of ``component fields'' will be discussed in a separate publication.


\end{document}